\documentclass[hidelinks,onefignum,onetabnum]{siamart250211}

\usepackage{lipsum}
\usepackage{booktabs}
\usepackage{amsfonts}
\usepackage{graphicx}
\usepackage{epstopdf}
\usepackage{caption}
\usepackage{hyperref}
\usepackage{makecell}
\usepackage{algorithmic}
\usepackage{wavelet}
\usepackage[font=footnotesize]{subcaption}
\usepackage{verbatim}
\usepackage{amsopn}
\usepackage{listings}
\usepackage{tikz}
\usepackage[normalem]{ulem}

\DeclareFontFamily{U}{mathx}{}
\DeclareFontShape{U}{mathx}{m}{n}{<-> mathx10}{}
\DeclareSymbolFont{mathx}{U}{mathx}{m}{n}
\DeclareMathAccent{\widehat}{0}{mathx}{"70}
\DeclareMathAccent{\widecheck}{0}{mathx}{"71}

\usepackage{hyperref}
\usepackage[nameinlink]{cleveref}

\usepackage[top=0.7in,bottom=0.7in,left=0.9in,right=0.9in]{geometry}

\allowdisplaybreaks

\ifpdf
  \DeclareGraphicsExtensions{.eps,.pdf,.png,.jpg}
\else
  \DeclareGraphicsExtensions{.eps}
\fi

\DeclareMathOperator* {\Span} {span}

\newtheorem{example}{Example} 
\newsiamremark{remark}{Remark}
\newsiamremark{hypothesis}{Hypothesis}
\crefname{hypothesis}{Hypothesis}{Hypotheses}
\newsiamthm{claim}{Claim}
\newsiamremark{fact}{Fact}
\crefname{fact}{Fact}{Facts}
\crefname{subsection}{Section}{Sections}
\Crefname{subsection}{Section}{Sections}

\crefname{example}{Example}{Examples}
\Crefname{example}{Example}{Examples}

\newcommand {\textB} [1] {\textcolor{blue}{#1}}
\newcommand {\ii} {\mathbf{i}} 

\newcommand{\cut}[1]{ \textcolor{red}{} } 
\renewcommand{\SS}{\mathcal{S}} 
\newcommand{\spt}{\mathbf{c}^*} 
\newcommand{\sptt}{\mathbf{d}^*}
\newcommand{\bpt}{\mathbf{b}^*} 
\newcommand{\LN}{\mathsf{\Pi}}

\begin{document}

\headers{Analysis of high-order, compact, and symmetric FDMs}{Q. Feng, B. Han, M. Michelle, and J. Sim}

\title{
Convergent Compact Finite Difference Methods with Symmetric Matrices and Highest Possible Orders for Variable-Coefficient Elliptic Equations\thanks{A preprint.
\funding{Research supported in part by Natural Sciences and Engineering Research Council (NSERC) of Canada under grants RGPIN-2024-04991, RGPIN 2026-05814, DGECR 2026-00371, the University of Alberta start-up grant, and the Digital Research Alliance of Canada.}}}

\author{Qiwei Feng\thanks{Mathematics Department, King Fahd University of Petroleum and Minerals (KFUPM), Dhahran, Eastern Province, Saudi Arabia
  (\email{qiwei.feng@kfupm.edu.sa,qfeng@ualberta.ca}).}
\and Bin Han\thanks{Department of Mathematical and Statistical Sciences, University of Alberta, Edmonton, Alberta, Canada
  (\email{\{bhan,mmichell,jiwoon2\}@ualberta.ca}).
  \protect\\
\textbf{Data Availability:}
No external datasets were used in this study. The numerical results are reported in the article.
\protect\\
\textbf{Conflict of interest:}
The authors declare that they have no competing interests.}
\and Michelle Michelle\footnotemark[3]
\and Jiwoon Sim\footnotemark[3]}






\maketitle

\begin{abstract}
Symmetric positive definite (SPD) linear systems from a discretization scheme enable efficient iterative solvers such as the conjugate gradient method.
High-order schemes producing SPD matrices are highly desirable for solving $d$-dimensional (dD) PDEs. For the dD elliptic equation $-\nabla \cdot(a\nabla u)=f$ with a \emph{variable coefficient} $a$, we seek compact (i.e., $3^d$-point) finite difference methods (FDMs) that produce symmetric linear systems and attain the \emph{highest possible consistency orders}. We are particularly interested in schemes that also produce SPD matrices, allowing the use of efficient iterative solvers and facilitating convergence analysis. The standard 2nd-order compact FDM for this elliptic problem is well known from textbooks and produces a symmetric matrix by applying 1D centered flux discretization separately in each coordinate direction.
However, to the best of our knowledge, a general construction of higher-order compact $3^d$-point FDMs with symmetric discretization matrices for variable coefficients in arbitrary dimensions has remained open. This paper provides a comprehensive
study of the highest consistency orders attainable by compact FDMs with symmetric discretization matrices for the variable-coefficient elliptic equation on a uniform grid in a
dD hypercube.For $d=1$, we construct compact FDMs that attain arbitrarily high consistency orders, yield tridiagonal SPD matrices,
and converge at the corresponding orders in the $\infty$-norm and the weighted $2$-norm.
For $d\ge2$, we prove that the highest consistency order attainable by compact FDMs with symmetric discretization matrices is $6$
when $\tfrac{|\nabla a|^2}{a^2}-\tfrac{2\Delta a}{a}$ is constant (i.e., $a^{1/2}$ is an eigenfunction of the Laplace operator $\Delta$),
and 4 when this quantity is nonconstant.
For $d=2,3,4$ and $M=4,6$ (requiring $a^{1/2}$ to be an eigenfunction of $\Delta$ for $M=6$), we further prove that our proposed compact schemes of consistency order $M$ produce SPD matrices
and achieve $M$th-order convergence in both the weighted $2$-norm and the $\infty$-norm. 
To facilitate implementation and reproducibility, we provide explicit stencil coefficients for the proposed convergent schemes in dimensions $1\le d\le 4$, together with a general construction procedure. Numerical experiments demonstrate the convergence and performance of the proposed schemes.
\end{abstract}




\section{Introduction and motivation}
\label{sec:intro}
    We consider the following $d$-dimensional (dD) elliptic equation with variable coefficients and the Dirichlet boundary condition:
    \begin{equation}
    \label{eq:PDE}
            \mathcal{L} u := -\nabla \cdot (a \nabla u) = f \quad \mbox{in } \Omega \quad \mbox{with}\quad
            u = g \quad \mbox{on } \partial \Omega,
    \end{equation}
    where $a,f,g$ are sufficiently smooth functions and $a(\mathbf{x}) \ge a_0 >0$ for all $\mathbf{x} \in \Omega$. Such an equation models many physical phenomena such as fluid flow in heterogeneous media \cite{B72book}, steady-state heat conduction in media where the thermal conductivity varies, electrostatic potential from a charge distribution \cite{GD07}, and semiconductors \cite{S84book}. In this paper we are interested in developing high-order compact finite difference methods (FDMs) to solve the elliptic equation \eqref{eq:PDE} such that the linear systems are symmetric positive definite (SPD) and the consistency order is as high as possible. Following the convention of FDMs, the domain $\Omega$ is assumed to be a dD hypercube. For simplicity of discussion, we further assume $\Omega=(0,1)^d$ and mesh size $h=1/N$, where $N\ge2$ is an integer. To find an approximate solution, $u_h$, to \eqref{eq:PDE} using FDMs with a given grid size $h>0$, we can use the 2nd-order FDM readily available in standard textbooks (e.g., \cite{JS14book,L07book})
    \begin{align*}
        & \sum_{i=1}^d \left[ - a(\spt + \tfrac{h}{2} \vec{e}_i) u_h(\spt + h \vec{e}_i)
        + (a(\spt + \tfrac{h}{2} \vec{e}_i) + a(\spt - \tfrac{h}{2} \vec{e}_i)) u_h(\spt) \right. \\
        & \qquad \left. - a(\spt - \tfrac{h}{2} \vec{e}_i) u_h(\spt - h \vec{e}_i) \right] = h^2 f(\spt), \quad \text{for all } \spt \in \Omega_h,
    \end{align*}
    where $\Omega_h := \Omega \cap (h \mathbb{Z}^d)$ and each $\vec{e}_i$ is the $i$th standard basis vector of $\R^d$. This scheme is easily obtained by applying the 1D centered flux discretization separately in each coordinate direction, which automatically makes the discretization matrix symmetric.

    SPD linear systems from a discretization scheme enable efficient iterative solvers such as the conjugate gradient method (CGM). High-order schemes producing SPD matrices are highly desirable for solving dD elliptic equation \eqref{eq:PDE}. Two natural questions arise:
    \begin{itemize}
    \item[(Q1)] \emph{What is the highest consistency order attainable by a compact $3^d$-point FDM for \eqref{eq:PDE} with variable $a$ that produces a symmetric linear system?}
    \item[(Q2)] \emph{Furthermore, do schemes attaining this highest consistency order produce SPD matrices?}
    \end{itemize}
    It is striking that, despite decades of research on FDMs for \eqref{eq:PDE}, the fundamental question (Q1) remains unresolved in arbitrary dimensions, to the best of our knowledge. Our extensive literature review, including detailed searches of online databases such as MathSciNet, identified no characterization of the highest consistency orders attainable by compact schemes that produce symmetric linear systems for general variable coefficients in arbitrary dimensions. To address the gaps in the literature, we aim to
    \begin{itemize}
        \item[(1)] fully characterize the highest consistency order attainable by compact $3^d$-point FDMs on uniform grids for \eqref{eq:PDE} with variable $a$ that produce symmetric linear systems,
        \item[(2)] explain how such schemes can be constructed, and
        \item[(3)] prove the SPD property and convergence, whenever possible.
    \end{itemize}

    As a consequence of our results established in this paper, our answer to (Q1) is arbitrarily high orders for $d=1$, order 6 for $d\ge 2$ and constant $\tfrac{|\nabla a|^2}{a^2}-\tfrac{2\Delta a}{a}$ (i.e., $a^{1/2}$ is an eigenfunction of the Laplace operator $\Delta$, because $\tfrac{|\nabla a|^2}{a^2}-\tfrac{2\Delta a}{a}=\lambda$ is equivalent to $\Delta (a^{1/2})= - \tfrac{\lambda}{4} a^{1/2}$),
    and
    order 4 for $d\ge 2$ when this quantity is nonconstant.
    Moreover, we prove that our constructed schemes produce SPD matrices for $1\le d\le4$, thereby answering (Q2) affirmatively. An additional condition on $a$ is required for the 4th-order scheme when $d=4$. FDMs that produce SPD matrices are of particular interest because they allow efficient matrix storage, enable the use of many fast solvers, and facilitate proofs of convergence rates. From now on, we refer to FDMs that produce symmetric linear systems as \emph{\textbf{symmetric FDMs}} (see \cref{def}).


    \subsection{Challenges of high-order, compact, and symmetric FDMs}
    The above 2nd-order FDM is easy to construct: applying the one-dimensional centered flux discretization separately in each coordinate direction immediately produces a symmetric linear system. However, this simple construction does not directly extend to higher-order FDMs on the same compact $3^d$-point stencils. Compactness limits the available degrees of freedom, high-order consistency imposes numerous local Taylor conditions, and matrix symmetry requires compatibility between neighboring stencils. These requirements must be satisfied simultaneously. The compatibility required by the matrix symmetry is explained below.

    An FDM for the PDE \eqref{eq:PDE} seeks the approximate solution $u_h : \Omega_h \rightarrow \R$ such that
    \[
    \sum_{\sptt \in \Omega_h} C_{p}(\spt) u_h(\sptt) = f_h(\spt),
    \quad p:=h^{-1}(\sptt - \spt),
    \quad \text{for all } \spt \in \Omega_h,
    \]
    where each $C_{p}(\spt)$ is a known stencil coefficient centered at the base point $\spt$, and the grid function $f_h : \Omega_h \rightarrow \R$ is also known. The previous relation can be written as a linear system, where all stencil coefficients are stored in a matrix $A$. Let $\dot{\iota} (\spt)$, $\spt \in \Omega_h$, be a map from a grid point $\spt$ to its row index of the matrix $A$ under a given ordering of the grid points so that the matrix entry $A_{\dot{\iota}(\spt),\dot{\iota}(\sptt)}$ corresponds to the pair of grid points $\spt,\sptt \in \Omega_h$. The above relation gives
    \[
    A_{\dot{\iota}(\spt), \dot{\iota}(\sptt)} = C_{p}(\spt) \quad \text{and} \quad
    A_{\dot{\iota}(\sptt), \dot{\iota}(\spt)} = C_{-p}(\sptt).
    \]
    Now, the finite difference discretization matrix $A$ is symmetric iff for any $\spt,\sptt \in \Omega_h$, 
    \[
    A_{\dot{\iota}(\spt), \dot{\iota}(\sptt)} = A_{\dot{\iota}(\sptt), \dot{\iota}(\spt)} \quad \text{equivalently} \quad
    C_{p}(\spt) = C_{-p}(\spt + ph).
    \]
    This means that all stencil coefficients have to satisfy specific compatibility conditions so that the finite difference matrix $A$ is symmetric, which can be challenging when the stencil coefficients depend on their centers or base points. 



    Even on $\Omega = (0,1)^d$, characterizing and constructing compact symmetric FDMs of the highest possible consistency orders requires substantial technical analysis. Extensions to irregular domains or other PDEs serve as important future research directions. We anticipate some key techniques in this paper to carry over. 

    \subsection{Related work}
    When $d=1$, arbitrarily high-order, compact, but nonsymmetric FDMs can be constructed for \eqref{eq:PDE} with piecewise smooth $a$ \cite{HanMichelleWong}. This generalizes the 1D result of \cite{SDKS13} for $a=1$.

    When $d=2$, 6th-order compact symmetric FDMs for \eqref{eq:PDE} with $a=1$ were constructed in \cite{D20,SDKS13,SWSZ22,zhai2014new}. For $a=1$, \cite{SDKS13} also showed that the consistency order of 9-point compact FDMs is at most 4 on quasi-uniform grids and at most 3 on nonuniform grids. A 4th-order compact symmetric FDM for the special coefficient $a(x,y)=1/x$ was constructed in \cite{ICGI00}. For general variable $a$, \cite{feng2024sixth,han2025convergent} established that the highest consistency order attainable by compact, not necessarily symmetric FDMs is 6 on rectangular and curved domains. Second-order compact symmetric FDMs for irregular domains were developed in \cite{GFCK02}, where symmetry is preserved through linear extrapolation to ghost values near the boundary.

    When $d=3$ and $a=1$, 6th-order compact FDMs were developed in \cite{D20,spotz1996high,zhai2013family}, while \cite{FHLP26} presented a unified family of 4th-order compact schemes with convergence analysis. More generally, \cite{D20} proved that the highest consistency order attainable by compact symmetric FDMs for the constant-coefficient Poisson equation on uniform grids with periodic boundary conditions is 6 for all $d\ge2$. Arbitrarily high orders can be attained in any dimension by increasing the stencil size.

    For $d=2,3$, 4th-order compact schemes were developed in \cite{BTT11} for variable-coefficient Helmholtz equations, including \eqref{eq:PDE} as a special case, and in \cite{TFL21} for variable-coefficient diffusion equations. Fourth-order compact schemes for variable-diffusion equations in arbitrary dimensions were also developed in \cite{GL20}. These schemes generally do not produce symmetric matrices.

    Building on the summation-by-parts operators for second derivatives with variable coefficients developed in \cite{M12}, \cite{KEW21} proposed high-order symmetric but noncompact FDMs for variable-coefficient elliptic equations. Their interior stencils have up to 6th-order consistency, while global convergence of up to 5th-order was observed numerically.



    \subsection{Main contributions}
    We characterize the highest consistency orders attainable by compact $3^d$-point FDMs with symmetric discretization matrices for \eqref{eq:PDE} on uniform grids in $\Omega=(0,1)^d$. For $d=1$, arbitrarily high consistency orders are attainable. For $d\ge2$, the highest attainable consistency order is 4 when $\tfrac{|\nabla a|^2}{a^2}-\tfrac{2\Delta a}{a}$ is nonconstant on $\Omega$, and 6 when it is constant. The precise assumptions and results are stated in \cref{thm:1D,thm:max_order}.

    The quantity $\tfrac{|\nabla a|^2}{a^2}-\tfrac{2\Delta a}{a}$ also connects this characterization to Helmholtz-type equations and quantum mechanics. Under the transformation $v:=\sqrt{a} u$, \eqref{eq:PDE} becomes the Helmholtz-type equation $-\Delta v-\tfrac{1}{4}(\tfrac{|\nabla a|^2}{a^2}-\tfrac{2\Delta a}{a})v=f/\sqrt{a}$. If $a$ is interpreted as a positive probability density, then $\tfrac{|\nabla a|^2}{a^2}-\tfrac{2\Delta a}{a}$ is proportional to the Bohm quantum potential \cite[eqn (8)]{B52}. If $\tfrac{|\nabla a|^2}{a^2}-\tfrac{2\Delta a}{a}$ is constant in $\Omega$ (i.e., $\tfrac{|\nabla a|^2}{a^2}-\tfrac{2\Delta a}{a}=\lambda$ for a constant $\lambda$), then it is equivalent to 
    $\Delta (a^{1/2})= - \tfrac{\lambda}{4} a^{1/2}$. We can further say that
    $a^{1/2}$ is an eigenfunction of the Laplace operator $\Delta$.

    For $1\le d\le4$, we prove that the proposed schemes produce SPD matrices for sufficiently small grid sizes and converge at their respective consistency orders in both the weighted $2$-norm and the $\infty$-norm. The 4th-order result for $d=4$ requires an additional condition on $a$. These results are stated in \cref{SPD,SPD:4,SPD:sym:constant}.

    We provide a general construction procedure for compact symmetric FDMs and explicit stencil coefficients for the proposed convergent schemes in dimensions $1\le d\le4$. For $d\ge3$, our 4th-order schemes are linear combinations of schemes involving derivatives of the problem data along either one or two coordinate directions. This structure supports their construction in arbitrary dimensions. Numerical experiments demonstrate the convergence and performance of the proposed schemes.

    \subsection{Organization}
    We organize this paper in order of increasing technicality. For the reader's convenience, we list the notation and definitions in \Cref{notdef}. The main results on the highest attainable consistency orders in arbitrary dimensions are stated in \Cref{sec:main:thms}. We present explicit compact symmetric FDMs in dimensions $1\le d\le4$ and state their SPD and convergence results in \Cref{sec:stencils}. Numerical experiments verifying the convergence rates are presented in \Cref{sec:exp}. The construction procedure is described in \Cref{sec:construct:dD}, and the proofs and supporting details are provided in \Cref{sec:details}. Finally, we present concluding remarks in \Cref{sec:conclusion}.


\section{Notation and definitions}
    \label{notdef}
    For the sake of clarity, we list the notations and definitions that we use throughout the paper:
    \begin{itemize}
    \item $\mathbb{N}_0$ stands for the set of all nonnegative integers  $\N \cup \{0\}$.
    \item The standard unit vector basis of $\R^d$ is denoted by $\vec{e}_j$ with $1 \le j \le d$, while the $d$-dimensional zero vector is denoted by $\vec{0}$. The symbol $\rightarrow$ on top of a variable indicates that it is a $d$-dimensional real vector.
    \item Given a vector $\vec{v} \in \R^{N}$ with $\vec{v}_i$ denoting its $i$-th entry, define $\| \vec{v} \|_2^2 := \sum_{i=1}^{N} |\vec{v}_i|^2$ and $\| \vec{v} \|_\infty := \max_{1\le i \le N} |\vec{v}_i|$.
    \item Let $\bm{k} \in \N^{d}_0$ such that $\bm{k} := (k_1,\ldots,k_d)$. The sum of the components of a multi-index $\bm{k}$ is given by $|\bm{k}| := \sum_{j = 1}^d k_j$. For $\bm{\ell} \in \N^{d}_0$, we say that $\bm{k} \le \bm{\ell}$ if $k_j \le \ell_j$ for all $1 \le j \le d$. Meanwhile, if $\vec{0} \le \bm{\ell} \le \bm{k}$, the binomial coefficient is defined as $\binom{\bm{k}}{\bm{\ell}} := \prod_{j = 1}^d \binom{k_j}{l_j}$.
    \item For a smooth function $v$, $\partial^{\bm{k}} v$ stands for its $\bm{k}$-th partial derivative.
    \item We use $\spt \in \Omega$ to denote the center (base) point of a finite difference stencil.
    \item The notation $\bo(h^M)$ with various subscripts refers to a function that is bounded by $C h^M$ as $h \to 0^+$, where the constant $C$ only depends on the expressions and their derivatives in the subscript, and $C$ remains bounded if its dependencies are bounded.
    \item $\td(k)$, $k \in \Z$, is the sequence such that $\td(0) := 1$ and $\td(k) := 0$ for $k \neq 0$.
    \end{itemize}

    Throughout this paper, unless otherwise stated, we assume that we have a uniform grid on the domain $\Omega:=(0,1)^d$ with mesh size $h = N^{-1}$, $N \in \N$, and define $\Omega_h := \Omega \cap (h\Z^d)$. The discussion at the start of \cref{sec:intro} motivates the following definition of a compact, symmetric, and $M$-th order consistent finite difference scheme for \eqref{eq:PDE}.


    \begin{definition}
        \label{def}
        Consider a finite difference scheme $\mathcal{L}_h u_h = f_h$ for \eqref{eq:PDE} on $\Omega_h$, where $u_h$ and $f_h$ are grid functions on $\Omega_h$  (i.e., $u_h, f_h : \Omega_h \rightarrow \R$), $f_h$ depends on $a$, $f$, as well as their derivatives, and the discretization operator $\mathcal{L}_h$ has the form
        \begin{equation}
            \label{eq:FDM_form}
            \mathcal{L}_h u_h(\spt) = \sum_{p \in \SS} C_p(\spt) u_h(\spt + ph) \quad \text{for all} \quad \spt \in \Omega_h,
        \end{equation}
        where $C_p(\spt) \in \R$ and $ \SS = \{-1,0,1\}^d$. Such a scheme is called compact.
        \begin{itemize}

            \item[(a)] The finite difference scheme $\mathcal{L}_h u_h = f_h$ is also symmetric if
            \be \label{sym:cond}
            C_p(\spt) = C_{-p} (\spt + ph) \quad \forall \spt, \spt + ph \in \Omega_h, \; p \in \SS.
            \ee

            \item[(b)] Suppose there exists $\Lambda(h) \neq 0$ dependent only on $a$ and $h$ such that $f_h (\spt) = \Lambda(h) (f(\spt) + \so(1))$ uniformly in $\spt$. The finite difference scheme $\mathcal{L}_h u_h = f_h$ is also $M$th-order consistent (to the operator $\mathcal{L}$) for some $M \in \N$ if
            \begin{equation}
                \label{eq:consistency_def}
                \sup_{\spt \in \Omega_h}  \left| \mathcal{L}_h u(\spt) - f_h(\spt) \right| = \bo (\Lambda(h) h^M),
                \quad \forall u \in C^\infty (\overline{\Omega}),
            \end{equation}
            where $f_h$ is viewed as the result of a well-defined mapping on $a$ and $u$.
        \end{itemize}
    \end{definition}

\begin{figure}[htbp]
        \centering
        \includegraphics[width=0.3\textwidth]{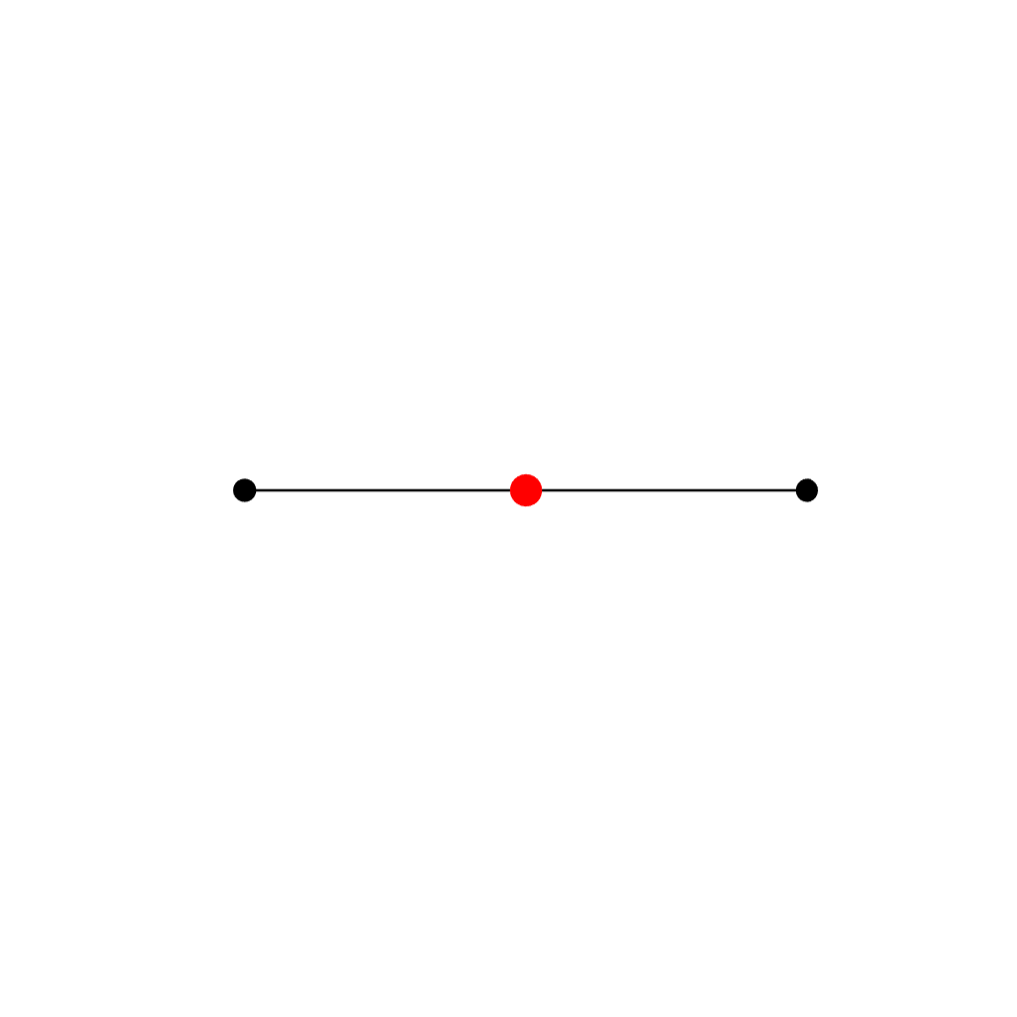}
        \includegraphics[width=0.3\textwidth]{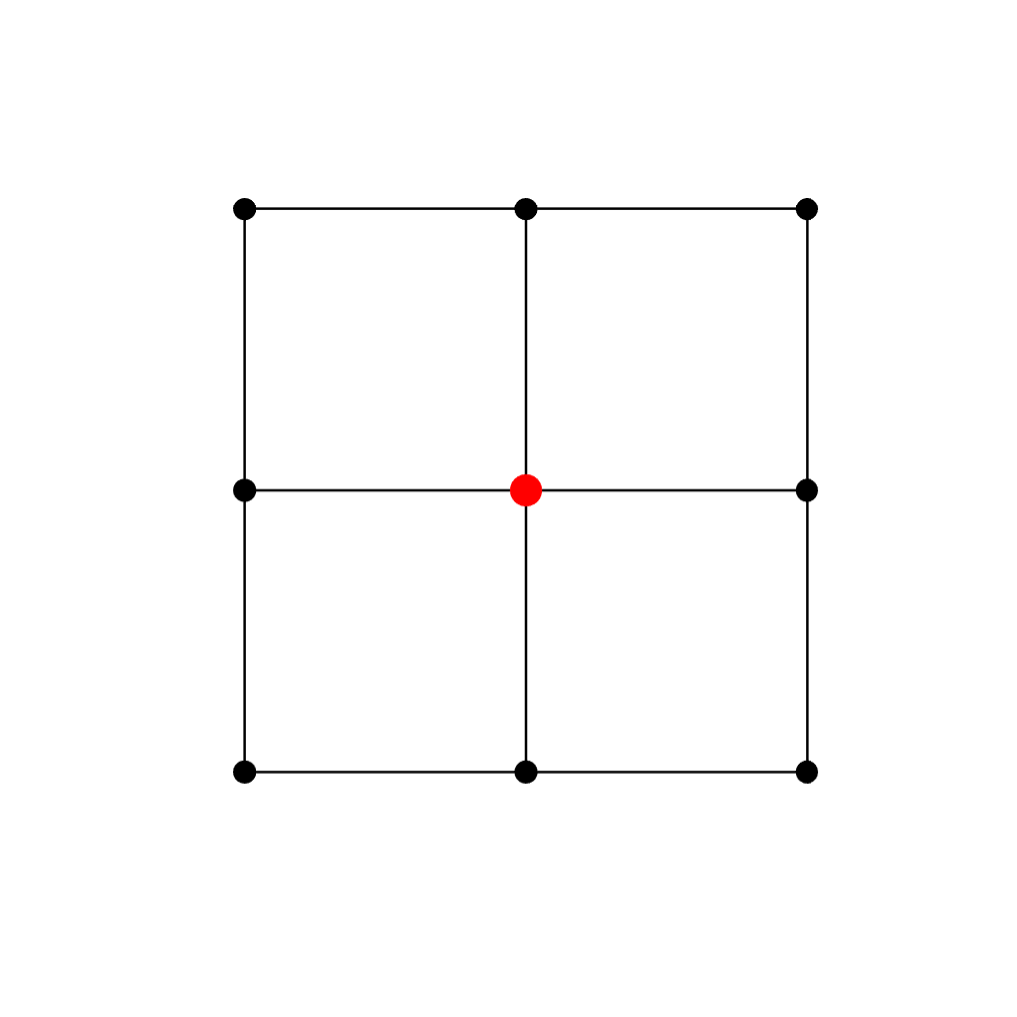}
        \includegraphics[width=0.3\textwidth]{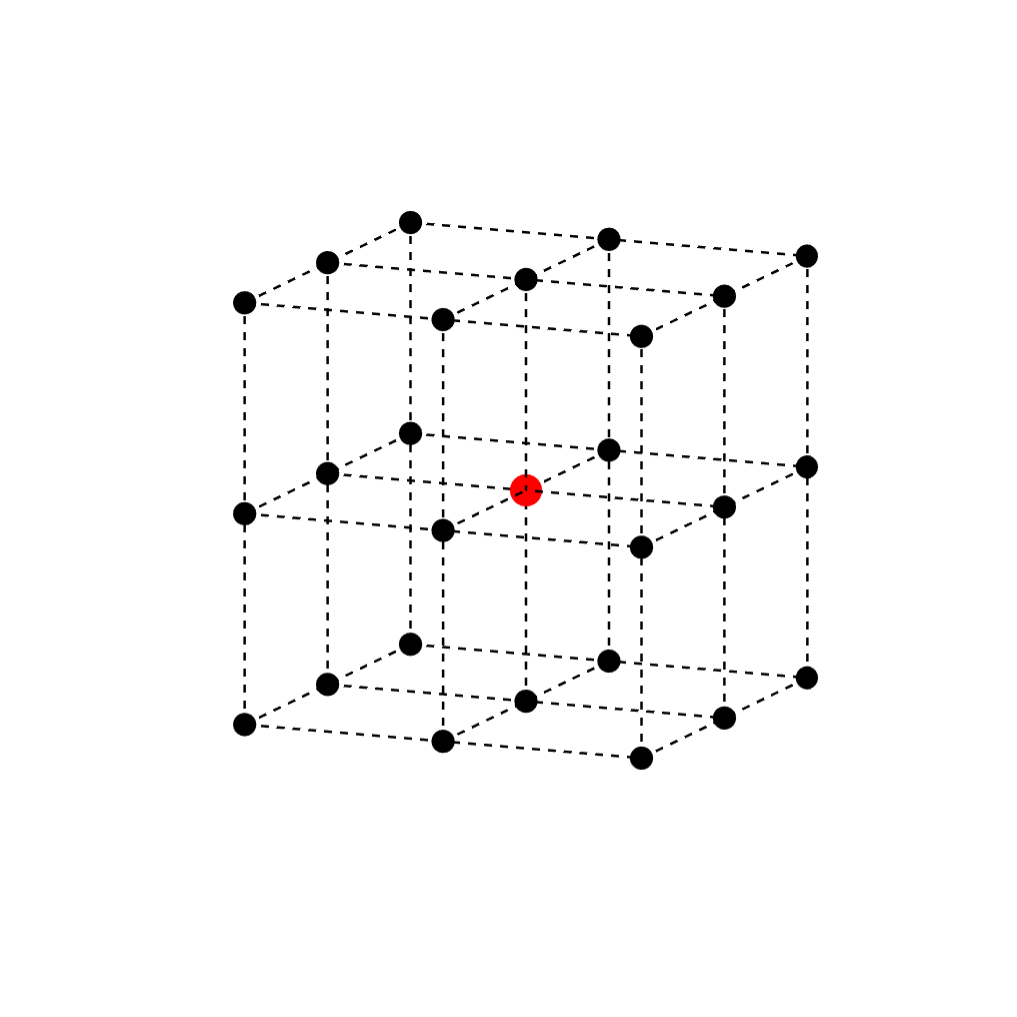}
        \caption{Compact finite difference stencils for $d=1,2,3$.}
        \label{fig:stencil:size}
    \end{figure}

    We emphasize again that, throughout this paper, our compact finite difference stencils are supported on the $3^d$-point stencil to solve \eqref{eq:PDE} with $\Omega = (0,1)^d$. See \cref{fig:stencil:size} for a visualization of the stencils, where $d=1,2,3$.

    Condition \eqref{eq:consistency_def} means the scheme produces a relative error of order $h^M$, instead of an absolute error. This prevents the trivial multiplication of $\so(1)$ factors in the stencil coefficients $C_p(\spt)$, and together with the stability estimates, allows the consistency order to yield to the corresponding convergence rate.

    \section{Maximum consistency orders}
    \label{sec:main:thms}
    In this section, we present the first two main theoretical results concerning the highest attainable consistency orders for compact and symmetric FDMs. To enhance readability, the proofs and detailed derivations are deferred to \cref{sec:details}. The first theorem states that a 1D compact, symmetric FDM can achieve arbitrarily high consistency order. The center point of a 1D finite difference stencil is a scalar, and thus we denote it by $c^{*}$.

    \begin{theorem} \label{thm:1D}
        Let $\Omega=(0,1)$, $M$ be a positive integer, $0<h<1$, and $c^* \in \Omega$ with $(c^* -h, c^* + h) \subset \Omega$. Assume that $a$ and $f$ are sufficiently smooth functions. There is a compact 3-point, symmetric, and $M$th-order consistent finite difference scheme for \eqref{eq:PDE} given by
        \begin{align*}
        \mathcal{L}_{h} u(c^*) & := C_{-1}(c^*)  u(c^* -h) + C_0(c^*)u(c^*) + C_{1}(c^*) u(c^* +h)  \\
        & = \sum_{\ell =0}^{M-1} d_{\ell}(h) h^{\ell + 2} f^{(\ell)}(c^*) + \bo_{a,f}(h^{M+2}),
        \end{align*}
        where $C_{-1}$, $C_0$, $C_1$, and $d_{\ell}$ are defined as
        \be \label{c:1D}
        \begin{aligned}
        & C_1(c^*) = -\mathcal{E}^{M}(c^* + \tfrac{h}{2}), \quad
        C_{-1}(c^*) = -\mathcal{E}^{M}(c^* - \tfrac{h}{2}), \quad
        C_0(c^*) = -C_1(c^*) - C_{-1}(c^*), \\ 
        & d_{\ell}(h) =
        -\tfrac{2((-1)^{\ell+1}-1)}{2^{\ell+2} (\ell+1)!}  + \sum_{k = 0}^{\ell} \tfrac{C_{-1}(c^*) (-1)^{\ell - k+1} \mathcal{G}_{k}(c^{*}-\tfrac{h}{2}) + C_{1}(c^*) \mathcal{G}_{k}(c^{*}+\tfrac{h}{2})}{2^{\ell+2} (\ell-k)!},
        \end{aligned}
        \ee
        such that
        \begin{align}
        \nonumber
        & \mathcal{E}^{M}(\cdot) := 2a(\cdot) \left(2 + \sum_{j=2}^{M+1} \tfrac{E_{j,1}(\cdot)}{j!} \Big(\tfrac{h}{2} \Big)^{j-1} (1+(-1)^{j-1}) \right)^{-1} + \bo_{a}(h^{M+1}),\\
        \label{Gk}
        &
        \mathcal{G}_{k}(\cdot) := \sum_{j= k + 2}^{M + 1} \tfrac{F_{j,k}(\cdot)}{j!} \Big(\tfrac{h}{2}\Big)^{j-k-2} (1+(-1)^{j-1}), 
        \end{align}
        and the quantities $E_{j,1}$, $F_{j,k}$ can be computed using these recursion relations
        \begin{align*}
            & E_{2,1} = -\tfrac{a'}{a}, \quad
            E_{j+1,1} = E'_{j,1}-\tfrac{a'}{a}E_{j,1}, \\
            & F_{2,0} = - \tfrac{1}{a}, \quad
            F_{j,-1} = - \tfrac{E_{j,1}}{a}, \quad
            F_{j+1,k} = F'_{j,k}+F_{j,k-1}.
        \end{align*}
    \end{theorem}
    \begin{proof}
        See \cref{sec:pf:1D}.
    \end{proof}

    The second theorem states the highest consistency orders that can be achieved by $d$-dimensional compact nonsymmetric and symmetric FDMs, where $d\ge 2$.

    \begin{theorem}
        \label{thm:max_order}
        Let $\Omega = (0,1)^d$ with $d\ge 2$, $M \in \N$, and $h > 0$.  Assume that $a$, $f$, and $g$ are sufficiently smooth functions. Consider the compact finite difference scheme $\mathcal{L}_h u_h = f_h$ for \eqref{eq:PDE} on $\Omega_h$ given by $\mathcal{L}_h u_h (\spt) := \sum_{p \in \SS} C_p (\spt) u_h(\spt + ph)$ with $\spt \in \Omega_h$, $C_p (\spt) := \sum_{k = 0}^{M + 1} c_{p, k}(\spt) h^k + \bo_{a} (h^{M + 2})$, $c_{p, k} = \bo_{a} (1)$ for $p \in \SS$, and
        $f_h(\spt):=h^2 f(\spt) + \sum_{k=1}^{M-1} f_{h,k}(\spt) h^{k+2} + \bo_{a,f}(h^{M+2})$, $f_{h,k} = \bo_{a,f}(1)$.
        Then, the following statements hold:

        \begin{itemize}
        \item[(a)] For any $d \ge 2$, if $\frac{|\nabla a|^2}{a^2} - \frac{2\Delta a}{a}$ is not constant on $\Omega$ (i.e., $a^{1/2}$ is not an eigenfunction of $\Delta$), then the maximum attainable consistency order of a compact $3^d$-point symmetric scheme is 4 and no higher.
        \item[(b)] For any $d\ge2$, if $\frac{|\nabla a|^2}{a^2}-\frac{2\Delta a}{a}$ is constant on $\Omega$ (i.e., $a^{1/2}$ is an eigenfunction of $\Delta$), then the maximum attainable consistency order of a compact $3^d$-point symmetric scheme is 6, which is also the maximum attainable when matrix symmetry is not required.
        \end{itemize}
    \end{theorem}

    \begin{proof}
        The proof establishes both the upper bounds and their attainability. The upper bound in item (a) is proved in \Cref{Qa:noncnst}, and schemes attaining this bound are constructed in \Cref{sec:basis,sec:any_dimension}. Item (b), together with the corresponding stencil constructions, is proved in \Cref{lem:sym:6th}. For schemes not required to be symmetric, the upper bound is proved in \Cref{lem:max_order_nonsymmetric}. Its attainability follows from the 6th-order basis stencils in \Cref{sec:basis}, extended to higher dimensions by adapting the construction procedure in \Cref{sec:any_dimension}.
    \end{proof}

    \section{Explicit finite difference stencils, SPD property, and convergence rates}
        \label{sec:stencils}
        For reproducibility purposes, we provide the explicit forms of our compact, high-order, and symmetric finite difference stencils.

        \subsection{Explicit stencils on $\Omega=(0,1)$}
        Given \cref{thm:1D}, a compact, 12th-order, and symmetric 1D finite difference scheme for \eqref{eq:PDE} can be obtained by taking
        \begin{equation}
        \label{E12}
        \begin{aligned}
        \mathcal{E}^{12}(c^{*} + \tfrac{h}{2}) & = a(c^{*}+\tfrac{h}{2})  + w_1(c^{*}+\tfrac{h}{2}) h^2 + w_2(c^{*}+\tfrac{h}{2}) h^4 + w_3(c^{*}+\tfrac{h}{2}) h^6 \\
        & \quad + w_4(c^{*}+\tfrac{h}{2}) h^8 + w_5(c^{*}+\tfrac{h}{2}) h^{10},
        \end{aligned}
        \end{equation}
        where
        \begin{align*}
            w_1(x) & := q_1 \tfrac{(a^{(1)}(x))^2}{a(x)} + q_2 a^{(2)}(x), \\
            w_2(x) & := q_3\tfrac{(a^{(1)}(x))^4}{(a(x))^3 } + q_4\tfrac{a^{(2)}(x)(a^{(1)}(x))^2}{(a(x))^2 } + q_{5}\tfrac{(a^{(2)}(x))^2}{a(x)} + q_{6}\tfrac{a^{(3)}(x)a^{(1)}(x)}{a(x)}
            + q_{7} a^{(4)}(x),\\
            w_3(x) & := q_{8}\tfrac{(a^{(1)}(x))^6}{(a(x))^5 } + q_{9}\tfrac{a^{(2)}(x)(a^{(1)}(x))^4}{(a(x))^4 } + q_{10}\tfrac{(a^{(2)}(x))^2(a^{(1)}(x))^2}{(a(x))^3}+ q_{11}\tfrac{(a^{(2)}(x))^3}{(a(x))^2 } \\
            & \quad + q_{12}\tfrac{a^{(3)}(x)(a^{(1)}(x))^3}{(a(x))^3 } + q_{13}\tfrac{a^{(3)}(x)a^{(2)}(x)a^{(1)}(x)}{(a(x))^2 } + q_{14}\tfrac{(a^{(3)}(x))^2}{a(x) }  \\
            & \quad + q_{15}\tfrac{a^{(4)}(x)(a^{(1)}(x))^2}{(a(x))^2} + q_{16}\tfrac{a^{(4)}(x)a^{(2)}(x)}{a(x) }+ q_{17}\tfrac{a^{(5)}(x)a^{(1)}(x)}{a(x)} +q_{18}a^{(6)}(x),\\
            w_4(x) & := q_{19}\tfrac{(a^{(1)}(x))^8}{(a(x))^7} + q_{20}\tfrac{ a^{(2)}(x)(a^{(1)}(x))^6}{(a(x))^6} + q_{21}\tfrac{ (a^{(2)}(x))^2(a^{(1)}(x))^4}{(a(x))^5} \\
            & \quad + q_{22}\tfrac{ (a^{(2)}(x))^3(a^{(1)}(x))^2}{(a(x))^4} + q_{23}\tfrac{ (a^{(2)}(x))^4}{(a(x))^3}  + q_{24}\tfrac{ a^{(3)}(x)(a^{(1)}(x))^5}{(a(x))^5} \\
            & \quad + q_{25}\tfrac{ a^{(3)}(x)a^{(2)}(x)(a^{(1)}(x))^3}{(a(x))^4} + q_{26}\tfrac{ a^{(3)}(x)(a^{(2)}(x))^2a^{(1)}(x)}{(a(x))^3} + q_{27}\tfrac{ (a^{(3)}(x))^2(a^{(1)}(x))^2}{(a(x))^3} \\
            & \quad  + q_{28} \tfrac{ (a^{(3)}(x))^2a^{(2)}(x)}{(a(x))^2} + q_{29} \tfrac{ a^{(4)}(x)(a^{(1)}(x))^4}{(a(x))^4}  +q_{30} \tfrac{ a^{(4)}(x)a^{(2)}(x)(a^{(1)}(x))^2}{(a(x))^3}  \\
            & \quad + q_{31}  \tfrac{ a^{(4)}(x)(a^{(2)}(x))^2}{(a(x))^2} + q_{32} \tfrac{ a^{(4)}(x)a^{(3)}(x)a^{(1)}(x)}{(a(x))^2} + q_{33} \tfrac{ (a^{(4)}(x))^2}{a(x)} \\
            & \quad + q_{34} \tfrac{ a^{(5)}(x)(a^{(1)}(x))^3}{(a(x))^3} + q_{35} \tfrac{ a^{(5)}(x)a^{(2)}(x)a^{(1)}(x)}{(a(x))^2} + q_{36} \tfrac{ a^{(5)}(x)a^{(3)}(x)}{a(x)} \\
            & \quad + q_{37} \tfrac{a^{(6)}(x)(a^{(1)}(x))^2}{(a(x))^2} + q_{38} \tfrac{ a^{(6)}(x)a^{(2)}(x)}{a(x)}  + q_{39} \tfrac{ a^{(7)}(x)a^{(1)}(x)}{a(x)} + q_{40} a^{(8)}(x),\\
            w_5(x) & := q_{41} \tfrac{ (a^{(1)}(x))^{10}}{(a(x))^9} + q_{42} \tfrac{ a^{(2)}(x)(a^{(1)}(x))^{8}}{(a(x))^8}  + q_{43} \tfrac{ (a^{(2)}(x))^{2}(a^{(1)}(x))^{6}}{(a(x))^7} \\
            & \quad +q_{44} \tfrac{ (a^{(2)}(x))^{3}(a^{(1)}(x))^{4}}{(a(x))^6} +q_{45} \tfrac{ (a^{(2)}(x))^{4}(a^{(1)}(x))^{2}}{(a(x))^5} +q_{46} \tfrac{ (a^{(2)}(x))^{5}}{(a(x))^4} \\
            & \quad + q_{47} \tfrac{ a^{(3)}(x)(a^{(1)}(x))^{7}}{(a(x))^7} + q_{48} \tfrac{ a^{(3)}(x)a^{(2)}(x)(a^{(1)}(x))^{5}}{(a(x))^6} + q_{49} \tfrac{ a^{(3)}(x)(a^{(2)}(x))^{2}(a^{(1)}(x))^{3}}{(a(x))^5} \\
            & \quad  + q_{50} \tfrac{ a^{(3)}(x)(a^{(2)}(x))^{3}a^{(1)}(x)}{(a(x))^4} + q_{51} \tfrac{ (a^{(3)}(x))^{2}(a^{(1)}(x))^{4}}{(a(x))^5} + q_{52} \tfrac{ (a^{(3)}(x))^{2}a^{(2)}(x)(a^{(1)}(x))^{2}}{(a(x))^4} \\
            & \quad + q_{53} \tfrac{ (a^{(3)}(x))^{2}(a^{(2)}(x))^{2}}{(a(x))^3} + q_{54} \tfrac{ (a^{(3)}(x))^{3}a^{(1)}(x)}{(a(x))^3} + q_{55} \tfrac{ a^{(4)}(x)(a^{(1)}(x))^{6}}{(a(x))^6} \\
            & \quad + q_{56} \tfrac{ a^{(4)}(x)a^{(2)}(x)(a^{(1)}(x))^{4}}{(a(x))^5} + q_{57} \tfrac{ a^{(4)}(x)(a^{(2)}(x))^{2}(a^{(1)}(x))^{2}}{(a(x))^4} + q_{58} \tfrac{ a^{(4)}(x)(a^{(2)}(x))^{3}}{(a(x))^3} \\
            & \quad + q_{59} \tfrac{ a^{(4)}(x)a^{(3)}(x)(a^{(1)}(x))^{3}}{(a(x))^4} + q_{60} \tfrac{ a^{(4)}(x)a^{(3)}(x)a^{(2)}(x)a^{(1)}(x)}{(a(x))^3} \\
            & \quad + q_{61} \tfrac{ a^{(4)}(x)(a^{(3)}(x))^{2}}{(a(x))^2} + q_{62} \tfrac{ (a^{(4)}(x))^{2}(a^{(1)}(x))^{2}}{(a(x))^3} + q_{63} \tfrac{ (a^{(4)}(x))^{2}a^{(2)}(x)}{(a(x))^2} \\
            & \quad + q_{64} \tfrac{ a^{(5)}(x)(a^{(1)}(x))^{5}}{(a(x))^5} + q_{65} \tfrac{ a^{(5)}(x)a^{(2)}(x)(a^{(1)}(x))^{3}}{(a(x))^4} + q_{66} \tfrac{ a^{(5)}(x)(a^{(2)}(x))^{2}a^{(1)}(x)}{(a(x))^3} \\
            & \quad  + q_{67} \tfrac{ a^{(5)}(x)a^{(3)}(x)(a^{(1)}(x))^{2}}{(a(x))^3} + q_{68} \tfrac{ a^{(5)}(x)a^{(3)}(x)a^{(2)}(x)}{(a(x))^2} + q_{69} \tfrac{ a^{(5)}(x)a^{(4)}(x)a^{(1)}(x)}{(a(x))^2} \\
            & \quad + q_{70} \tfrac{ (a^{(5)}(x))^{2}}{a(x)} + q_{71} \tfrac{ a^{(6)}(x)(a^{(1)}(x))^{4}}{(a(x))^4} + q_{72} \tfrac{ a^{(6)}(x)a^{(2)}(x)(a^{(1)}(x))^{2}}{(a(x))^3} \\
            & \quad + q_{73} \tfrac{ a^{(6)}(x)(a^{(2)}(x))^{2}}{(a(x))^2} + q_{74} \tfrac{ a^{(6)}(x)a^{(3)}(x)a^{(1)}(x)}{(a(x))^2}  + q_{75} \tfrac{ a^{(6)}(x)a^{(4)}(x)}{a(x)} \\
            & \quad + q_{76} \tfrac{ a^{(7)}(x)(a^{(1)}(x))^{3}}{(a(x))^3} + q_{77} \tfrac{ a^{(7)}(x)a^{(2)}(x)a^{(1)}(x)}{(a(x))^2} + q_{78} \tfrac{ a^{(7)}(x)a^{(3)}(x)}{a(x)} \\
            & \quad + q_{79} \tfrac{ a^{(8)}(x)(a^{(1)}(x))^{2}}{(a(x))^2}  + q_{80} \tfrac{ a^{(8)}(x)a^{(2)}(x)}{a(x)} + q_{81} \tfrac{ a^{(9)}(x)a^{(1)}(x)}{a(x)} + q_{82} a^{(10)}(x),
        \end{align*}
        whose coefficients take the following values
        \begin{align*}
        &q_{1} = \tfrac{ -1}{12},\quad q_{2} = \tfrac{ 1}{24},\quad q_{3} = \tfrac{ -1}{180},\quad q_{4} = \tfrac{ 17}{1440},\quad q_{5} = \tfrac{ -1}{720},\quad q_{6} = \tfrac{ -1}{240},\quad q_{7} = \tfrac{ 1}{1920},\\
        & q_{8} = \tfrac{ -11}{15120}, \quad q_{9} = \tfrac{ 23}{10080},\quad q_{10} = \tfrac{ -137}{80640},\quad q_{11} = \tfrac{ 11}{120960},\quad q_{12} = \tfrac{ -1}{1260},\quad q_{13} = \tfrac{ 31}{40320},\\
        & q_{14} = \tfrac{ -1}{16128},\quad q_{15} = \tfrac{ 31}{161280}, \quad q_{16} = \tfrac{ -1}{20160},\quad q_{17} = \tfrac{ -1}{26880},\quad q_{18} = \tfrac{ 1}{322560},\\
        & q_{19} = \tfrac{ -107}{907200}, \quad q_{20} = \tfrac{ 887}{1814400},\quad q_{21} = \tfrac{ -377}{604800}, \quad q_{22} = \tfrac{ 989}{4147200},\quad q_{23} = \tfrac{ -107}{14515200},\\
        & q_{24} = \tfrac{ -17}{100800}, \quad q_{25} = \tfrac{ 13}{37800}, \quad q_{26} = \tfrac{ -193}{1612800},\quad q_{27} = \tfrac{ -1}{22400}, \quad q_{28} = \tfrac{ 5}{387072},\\
        & q_{29} = \tfrac{ 101}{2419200}, \quad q_{30} = \tfrac{ -197}{3225600},\quad q_{31} = \tfrac{ 17}{3225600},\quad q_{32} = \tfrac{ 19}{1382400},\quad q_{33} = \tfrac{ -1}{2073600},\\
        & q_{34} = \tfrac{ -1}{120960}, \quad q_{35} = \tfrac{ 1}{129024},\quad q_{36} = \tfrac{ -1}{829440},\quad q_{37} = \tfrac{ 1}{774144},\quad q_{38} = \tfrac{ -1}{2903040}, \\
        & q_{39} = \tfrac{ -1}{5806080}, \quad q_{40} = \tfrac{ 1}{92897280},\quad q_{41} = \tfrac{ -2549}{119750400},\quad q_{42} = \tfrac{ 26263}{239500800},\\
        & q_{43} = \tfrac{ -94043}{479001600}, \quad q_{44} = \tfrac{ 67339}{479001600}, \quad q_{45} = \tfrac{ -35971}{1094860800},\quad q_{46} = \tfrac{ 2549}{3832012800},\\
        & q_{47} = \tfrac{ -751}{19958400}, \quad q_{48} = \tfrac{ 1319}{11404800},\quad q_{49} = \tfrac{ -1921}{19958400}, \quad q_{50} = \tfrac{ 22313}{1277337600},\\
        & q_{51} = \tfrac{ -379}{22809600}, \quad q_{52} = \tfrac{ 39}{1971200}, \quad q_{53} = \tfrac{ -1087}{510935040}, \quad q_{54} = \tfrac{ -1}{887040},\quad q_{55} = \tfrac{ 37}{3942400},\\
        & q_{56} = \tfrac{ -541}{22809600}, \quad q_{57} = \tfrac{ 7643}{567705600}, \quad q_{58} = \tfrac{ -751}{1277337600},\quad q_{59} = \tfrac{ 521}{79833600},\\
        & q_{60} = \tfrac{ -5743}{1277337600}, \quad q_{61} = \tfrac{ 83}{340623360}, \quad q_{62} = \tfrac{ -17669}{30656102400},\quad q_{63} = \tfrac{ 101}{958003200},\\
        & q_{64} = \tfrac{ -299}{159667200}, \quad q_{65} = \tfrac{ 601}{159667200}, \quad q_{66} = \tfrac{ -1087}{851558400},\quad q_{67} = \tfrac{ -29}{29937600},\\
        & q_{68} = \tfrac{ 59}{218972160}, \quad q_{69} = \tfrac{ 83}{567705600}, \quad q_{70} = \tfrac{ -1}{162201600}, \quad q_{71} = \tfrac{ 193}{638668800},\\
        & q_{72} = \tfrac{ -3349}{7664025600}, \quad q_{73} = \tfrac{ 299}{7664025600}, \quad q_{74} = \tfrac{ 83}{851558400},\quad q_{75} = \tfrac{ -1}{141926400},\\
        & q_{76} = \tfrac{ -1}{23950080}, \quad q_{77} = \tfrac{ 59}{1532805120}, \quad q_{78} = \tfrac{ -1}{170311680},\quad q_{79} = \tfrac{ 59}{12262440960},\\
        & q_{80} = \tfrac{ -1}{766402560}, \quad q_{81} = \tfrac{ -1}{2043740160}, \quad q_{82} = \tfrac{ 1}{40874803200}.
        \end{align*}
        Other $2n$-th order schemes with $1\le n \le 5$ can be recovered by simply neglecting all terms involving $h^{2k}$ with $k \ge n$ in $\mathcal{E}^{12}(c^* + \tfrac{h}{2})$. The foregoing expression is necessarily lengthy because the scheme has 12th-order consistency.

        \subsection{Explicit stencils on $\Omega=(0,1)^d$ with $2 \le d \le 4$ and nonconstant $|\nabla a |^2 / a^2 - 2\Delta a / a$}
        The $d$-dimensional compact, 4th-order, and symmetric finite difference stencils with $d=2,3,4$ take the following form
        \begin{equation}
        \label{stencil:explicit:uf}
        \mathcal{L}_h u (\spt) := \sum_{p \in \mathcal{S}} C_p(\spt) u(\spt + ph)  = f_h(\spt) + \bo_{a,f}(h^{6}),
        \end{equation}
        where $\mathcal{S}:=\{-1,0,1\}^d$. When $d=2$, the stencil coefficients $C_{p}(\spt)$ and the right-hand term $f_h(\spt)$ in \eqref{stencil:explicit:uf} take the following forms
        %
        %
        \begin{align}
            \nonumber
            & C_{p}(\spt) := \\
            \nonumber
            & \quad
            \left[ 
            \begin{aligned}
                & \begin{aligned}
                    & - \tfrac{a(\spt+ \tfrac{p h}{2})}{6} + \tfrac{a_x(\spt+ \tfrac{p h}{2}) a_y(\spt+ \tfrac{p h}{2})}{24a(\spt+ \tfrac{p h}{2})}h^2 \\
                    & \quad - \tfrac{4a_{xy}(\spt+ \tfrac{p h}{2})-a_{xx}(\spt+ \tfrac{p h}{2})}{48}h^2,
                \end{aligned}
                && p = (\pm 1,\pm 1),\\
                & - \tfrac{a(\spt+ \tfrac{p h}{2})}{6} - \left[\tfrac{a_x(\spt+ \tfrac{p h}{2}) a_y(\spt+ \tfrac{p h}{2})}{24a(\spt+ \tfrac{p h}{2})}-\tfrac{a_{xx}(\spt+ \tfrac{p h}{2})}{48}\right]h^2, && p = (\pm 1,\mp 1),\\
        %
                & \begin{aligned}
                    & - \tfrac{2a(\spt+ \tfrac{p h}{2})}{3} + \tfrac{a_x^2(\spt+ \tfrac{p h}{2})}{12a(\spt+ \tfrac{p h}{2})}h^2 - \tfrac{1}{24}[ 2a_{xx}(\spt+ \tfrac{p h}{2}) \\
                    & \quad -2a_{xy}(\spt+ \tfrac{p h}{2}) -a_{yy}(\spt+ \tfrac{p h}{2}) ]h^2,
                \end{aligned}
                && p = (\pm 1,0),\\
                & \begin{aligned}
                    & -\tfrac{2a(\spt+ \tfrac{p h}{2})}{3} + \tfrac{a_y^2(\spt+ \tfrac{p h}{2})}{12a(\spt+ \tfrac{p h}{2})}h^2 - \tfrac{1}{24}[ a_{yy}(\spt+ \tfrac{p h}{2}) \\
                    & \quad -2a_{xy}(\spt+ \tfrac{p h}{2}) ]h^2,
                \end{aligned}
                && p=(0,\pm 1),\\
                & -\sum_{p \in \mathcal{S} \backslash \{(0,0)\}} C_p (\spt), && p = (0,0),
            \end{aligned}
            \right. \\
            \nonumber
            & f_h(\spt) := f(\spt) h^2 - \tfrac{1}{12a^2(\spt)}[  a(\spt)(f(\spt)\Delta a(\spt)  +\nabla a(\spt) \cdot \nabla f(\spt)) \quad \\
            \label{2d:explicit}
            & \quad -a^2(\spt)\Delta f(\spt) - \|\nabla a(\spt)\|^2_2 f(\spt)  ]h^4.
        \end{align}

        Meanwhile, when $d=3$, the stencil coefficients $C_{p}(\spt)$ and the right-hand term $f_h(\spt)$ in \eqref{stencil:explicit:uf} take the following forms
        \begin{align*}
            & C_{p}(\spt) := \\
            & \quad
            \left\lceil
            \begin{aligned}
                & -\tfrac{7a(\spt+ \tfrac{p h}{2})}{15}, && p = (\pm 1,0,0), \\
                & -\tfrac{a(\spt+ \tfrac{p h}{2})}{10}, \hspace{7.1cm} && p=(0,\pm 1,\pm 1),\\
                & -\tfrac{a(\spt+ \tfrac{p h}{2})}{30} + \tfrac{a_y^2(\spt+ \tfrac{p h}{2})}{24 a(\spt+ \tfrac{p h}{2})} h^2, \hspace{4.5cm} && p=(\pm 1,\mp 1,\pm 1), \\
                & - \tfrac{a(\spt+ \tfrac{p h}{2})}{10} - \tfrac{a_y(\spt+ \tfrac{p h}{2}) a_z(\spt+ \tfrac{p h}{2})}{12a(\spt+ \tfrac{p h}{2})}h^2, && p = (0,\pm 1,\mp 1),
              \end{aligned}
            \right.
        \end{align*}
        \begin{align}
            \nonumber
            & \quad
            \left\lfloor
            \begin{aligned}
                & - \tfrac{a(\spt+ \tfrac{p h}{2})}{10} - \tfrac{a_x(\spt+ \tfrac{p h}{2}) a_z(\spt+ \tfrac{p h}{2})}{12a(\spt+ \tfrac{p h}{2})} h^2, && p=(\pm 1,0,\mp 1),\\
                & - \tfrac{a(\spt+ \tfrac{p h}{2})}{30} - \tfrac{a_y^2(\spt+ \tfrac{p h}{2}) - a_x^2(\spt+ \tfrac{p h}{2})}{24 a(\spt+ \tfrac{p h}{2})}h^2, && p=(\pm 1,\pm 1,\mp 1),\\
                & - \tfrac{7 a(\spt+ \tfrac{p h}{2})}{15} - \tfrac{a_z(\spt+ \tfrac{p h}{2})}{12a(\spt+ \tfrac{p h}{2})} [a_x(\spt+ \tfrac{p h}{2}) - a_y(\spt+ \tfrac{p h}{2})  ]h^2, && p=(0,\pm 1,0),\\
                & \begin{aligned}
                    & - \tfrac{a(\spt+ \tfrac{p h}{2})}{30} - \tfrac{1}{24}[a_{yy}(\spt+ \tfrac{p h}{2}) +a_{xz}(\spt+ \tfrac{p h}{2}) \\
                    &  \quad +a_{yz}(\spt+ \tfrac{p h}{2})-a_{xx}(\spt+ \tfrac{p h}{2})]h^2,
                    \end{aligned}
                    && p=(\pm 1,\pm 1, \pm 1),\\
                 & \begin{aligned}
                    & - \tfrac{a(\spt+ \tfrac{p h}{2})}{10} - \tfrac{a_y^2(\spt+ \tfrac{p h}{2}) -a_x^2(\spt+ \tfrac{p h}{2})}{12a(\spt+ \tfrac{p h}{2})}h^2 \\
                    & \quad - \tfrac{1}{12}[a_{xx}(\spt+ \tfrac{p h}{2})-a_{yy}(\spt+ \tfrac{p h}{2})]h^2,
                    \end{aligned}
                    && p = (\pm 1,0 ,\pm 1),\\
                & \begin{aligned}
                & - \tfrac{7a(\spt+ \tfrac{p h}{2})}{15} - \tfrac{1}{12a(\spt+ \tfrac{p h}{2})} [2a_x^2(\spt+ \tfrac{p h}{2}) -a_y^2(\spt+ \tfrac{p h}{2}) \\
                & \quad -a_z^2(\spt+ \tfrac{p h}{2}) -a_x(\spt+ \tfrac{p h}{2}) a_z(\spt+ \tfrac{p h}{2}) \\
                & \quad -a_y(\spt+ \tfrac{p h}{2}) a_z(\spt+ \tfrac{p h}{2})]h^2 - \tfrac{1}{24}[a_{yy}(\spt+ \tfrac{p h}{2})\\
                & \quad  +a_{zz}(\spt+ \tfrac{p h}{2}) -3a_{xx}(\spt+ \tfrac{p h}{2}) -2a_{yz}(\spt+ \tfrac{p h}{2})]h^2,
                \end{aligned}
                && p = (0,0,\pm 1),\\
                & \begin{aligned}
                & - \tfrac{a(\spt+ \tfrac{p h}{2})}{30} + \tfrac{a_x^2(\spt+ \tfrac{p h}{2})}{24a(\spt+ \tfrac{p h}{2})}h^2 - \tfrac{1}{24}[a_{xx}(\spt+ \tfrac{p h}{2}) \\
                & \quad +a_{yz}(\spt+ \tfrac{p h}{2}) -a_{yy}(\spt+ \tfrac{p h}{2}) \\
                & \quad -a_{xz}(\spt+ \tfrac{p h}{2})]h^2,
                \end{aligned}
                && p = (\pm 1,\mp 1,\mp 1),\\
                & \begin{aligned}
                & - \tfrac{a(\spt+ \tfrac{p h}{2})}{10} - \tfrac{1}{24a(\spt+ \tfrac{p h}{2})} [a_x^2(\spt+ \tfrac{p h}{2}) \\
                & \quad + a_x(\spt+ \tfrac{p h}{2}) a_y(\spt+ \tfrac{p h}{2}) \\
                & \quad -a_x(\spt+ \tfrac{p h}{2}) a_z(\spt+ \tfrac{p h}{2})]h^2 \\
                & \quad - \tfrac{1}{48}[3a_{yy}(\spt+ \tfrac{p h}{2}) -3a_{xx}(\spt+ \tfrac{p h}{2}) -a_{zz}(\spt+ \tfrac{p h}{2})\\
                & \quad +2a_{xz}(\spt+ \tfrac{p h}{2}) -2a_{xy}(\spt+ \tfrac{p h}{2}) -2a_{yz}(\spt+ \tfrac{p h}{2})]h^2,
                \end{aligned}
                && p=(\pm 1,\mp 1,0),\\
                & \begin{aligned}
                & - \tfrac{a(\spt+ \tfrac{p h}{2})}{10} - \tfrac{1}{24a(\spt+ \tfrac{p h}{2})} [a_x^2(\spt+ \tfrac{p h}{2}) - 2a_y^2(\spt+ \tfrac{p h}{2})\\
                & \quad -a_x(\spt+ \tfrac{p h}{2}) a_y(\spt+ \tfrac{p h}{2}) -a_x(\spt+ \tfrac{p h}{2}) a_z(\spt+ \tfrac{p h}{2})]h^2\\
                & \quad - \tfrac{1}{48}[a_{xx}(\spt+ \tfrac{p h}{2}) -a_{yy}(\spt+ \tfrac{p h}{2}) -a_{zz}(\spt+ \tfrac{p h}{2})\\
                & \quad +2a_{xy}(\spt+ \tfrac{p h}{2}) -2a_{xz}(\spt+ \tfrac{p h}{2}) -2a_{yz}(\spt+ \tfrac{p h}{2})]h^2,
                \end{aligned}
                && p =(\pm 1,\pm 1,0),\\
                & -\sum_{p \in \mathcal{S} \backslash \{(0,0,0)\}} C_p (\spt), && p = (0,0,0),
            \end{aligned} \right.\\
            \nonumber
            & f_h(\spt)
            := f(\spt) h^2 - \tfrac{1}{12 a^2(\spt)} [ a(\spt)(f(\spt) \Delta a(\spt) +\nabla a(\spt) \cdot \nabla f(\spt)) \hspace{2.2cm} \\
            \label{3d:explicit}
            & \quad -a^2(\spt) \Delta f(\spt) - \|\nabla a(\spt)\|^2_2 f(\spt)  ]h^4.
        \end{align}

        Define $\tilde{a} := - \ln a$. Finally, when $d=4$, the stencil coefficients $C_{p}(\spt)$ and the right-hand term $f_h(\spt)$ in \eqref{stencil:explicit:uf} take the following forms
        \begin{align}
        \nonumber
        C_p(\spt) & =
        \begin{cases}
            \tfrac{h^2}{24}a(|\nabla \tilde{a}|^2 + 2 \partial^{2\vec{e}_i}\tilde{a} - \Delta \tilde{a}) \left. \right|_{\spt + ph/2}, &
            p = \pm \vec{e}_i, \; 1\le i \le 4,\\
            a \left(-\tfrac{1}{6} + \tfrac{h^2}{24} \partial^{\vec{e}_i + \vec{e}_j}\tilde{a} \right) \left. \right|_{\spt + ph/2}, &
            p = \pm (\vec{e}_i+\vec{e}_j), \; 1\le i < j \le 4,\\
            a \left(-\tfrac{1}{6} - \tfrac{h^2}{24} \partial^{\vec{e}_i + \vec{e}_j}\tilde{a} \right) \left. \right|_{\spt + ph/2}, &
            p = \pm (\vec{e}_i-\vec{e}_j), \; 1\le i < j \le 4,\\
            - \sum_{p \in \SS \backslash\{\vec{0}\}} C_p(\spt), & p = \vec{0},\\
            0, & \text{otherwise},
        \end{cases}\\
        \label{Cp:fh:d4}
        f_h(\spt) &  = f(\spt)h^2 + \tfrac{h^4}{12} \nabla \cdot (a \nabla (f/a)) \left. \right|_{\spt}.
        \end{align}

       Note that each of the above finite difference stencils is only one possible example. Unlike the 1D case, their construction is not unique. We present a general construction framework in \cref{sec:construct:dD}, with the supporting details deferred to \cref{sec:details}.

    \subsection{Explicit stencils on $\Omega=(0,1)^d$ with $d \ge 2$ and constant $|\nabla a |^2 / a^2 - 2\Delta a / a$}
    The $d$-dimensional, 6th-order, and symmetric finite difference stencils $d \ge 2$ take the following form
    \[
        \mathcal{L}_h u (\spt) = \sum_{p \in \mathcal{S}} C_p (\spt) u(\spt + ph)  = f_h (\spt) + \bo_a(h^8),
    \]
    where
    \begin{align}
        \nonumber
        C_p (\spt) & := -\sqrt{a(\spt)} \widehat{C}_p(\spt) \sqrt{a(\spt + ph)}, \quad p \in \mathcal{S} \backslash \{\vec{0}\},\\
        \nonumber
        C_{\vec{0}}(\spt) & := - \sum_{p \in \mathcal{S} \backslash \{\vec{0}\}} C_p (\spt), \\
        \label{stencil:dD:constant}
        \widehat{C}_p(\spt) & :=
        \begin{cases}
            1-\tfrac{d-1}{3} + \tfrac{2}{15}\binom{d-1}{2} - \tfrac{(d-1)qh^2}{360}, & |p| = 1,\\
            \tfrac{9-2d}{30} + \tfrac{qh^2}{720}, & |p| = 2,\\
            \tfrac{1}{30}, & |p| = 3,\\
            0, & 4 \le |p| \le d,
        \end{cases}\\
        \nonumber
        f_h (\spt) & := \sqrt{a(\spt)} ( h^2 F(\spt) + \tfrac{h^4}{12} (\Delta F(\spt) - \tfrac{q}{4} F(\spt))
        + h^6 ( \tfrac{1}{360}\Delta^2F(\spt) \\
        \nonumber
        & \quad - \tfrac{q}{1440}\Delta F(\spt) + \tfrac{q^2}{5760} F(\spt)  +\tfrac{1}{180}\sum_{i<j}
       \partial^{2 (\vec{e}_i + \vec{e}_j)}F(\spt) )),\\
       \nonumber
       q & := |\nabla a|^2 / a^2 - 2\Delta a / a, \quad F:= f/\sqrt{a},
    \end{align}
    $|p| = |(p_1,\ldots,p_d)| := \sum_{i=1}^d |p_i|$, and by convention $ \binom{m}{n} :=0$ for all nonnegative integers $m<n$.

      \subsection{SPD property and convergence rates}
      \label{subsec:SPD:conv:thms}
      We introduce some notations used in the next theorems. Let $A_h \vec{u}_h = \vec{f}_h$ be the linear system arising from the discretization, where $A_h$ is a square matrix that depends on given stencil coefficients $C_p(\spt)$ with $p\in \mathcal{S}$, $\spt \in \Omega_h$, $\vec{u}_h$ is the vectorized form of the grid function $u_h$, and $\vec{f}_h$ is the right-hand-side vector obtained from the grid function $f_h$ after incorporating the known Dirichlet boundary condition $g$, under a given ordering of the grid points. Denote the vectorized form of $u|_{\Omega_h}$ under a given ordering of the grid points by $\overrightarrow{u|_{\Omega_h}}$.

      Below, we present our third main theoretical result concerning the SPD property of the linear systems arising from the discretization using the proposed $d$-dimensional compact, symmetric, high-order FDMs for $1 \le d \le 3$, as well as the convergence rates of the schemes.

    \begin{theorem}
        \label{SPD}
        Let $\Omega = (0,1)^d$ with $d=1,2,3$ and $h > 0$. Assume that $u$, $a$, $f$, and $g$ in \eqref{eq:PDE} are sufficiently smooth functions on $\overline{\Omega}$. Consider the compact finite difference scheme $\mathcal{L}_h u_h = f_h$ for \eqref{eq:PDE} on $\Omega_h$ given by $\mathcal{L}_h u_h (\spt) := \sum_{p \in \SS} C_p (\spt) u_h(\spt + ph)$ for all $\spt \in \Omega_h$, where $C_p (\spt)$ and $f_h(\spt)$ are given in \cref{thm:1D} for $d=1$, in \eqref{2d:explicit} for $d=2$, and in \eqref{3d:explicit} for $d=3$.
        Then, the following statements hold:
        \begin{itemize}
            \item[(a)] For all sufficiently small $h$ and $d=1,2,3$, the matrix $A_h$ is SPD.
            \item[(b)] For all sufficiently small $h$, the errors of the one-dimensional finite difference scheme satisfy
            \[
             \|\overrightarrow{u|_{\Omega_h}} - \vec{u}_h\|_{2,h} := h^{d/2} \|\overrightarrow{u|_{\Omega_h}} - \vec{u}_h\|_{2} \le C h^M \quad \text{and} \quad \|\overrightarrow{u|_{\Omega_h}} - \vec{u}_h\|_\infty \le C h^{M},
            \]
            where $M$ is any positive integer corresponding to the consistency order of the scheme, and $C$ is a positive constant independent of $h$. Meanwhile, for all sufficiently small $h$, the errors of two and three-dimensional finite difference schemes satisfy
            \[
             \|\overrightarrow{u|_{\Omega_h}} - \vec{u}_h\|_{2,h} := h^{d/2} \|\overrightarrow{u|_{\Omega_h}} - \vec{u}_h\|_{2} \le C h^{4} \quad \text{and} \quad \|\overrightarrow{u|_{\Omega_h}} - \vec{u}_h\|_\infty \le C h^{4},
            \]
            where $C$ is a positive constant independent of $h$.
        \end{itemize}
    \end{theorem}
    \begin{proof}
        See \cref{subsec:SPD:conv}.
    \end{proof}

    The next theorem is our fourth main result concerning the SPD property of the linear systems arising from the discretization using the proposed four-dimensional compact, symmetric, and 4th-order FDM as well as its convergence rates. The key difference here is we need an extra condition on the function $a$ in the form of $|\nabla \tilde{a}|^2 + 2\partial^{2\vec{e}_i}\tilde{a} - \Delta \tilde{a} \le 0$, where $\tilde{a}:=-\ln a$, on $\Omega$ for $i=1,\ldots,4$. An example of such a function is $a(\mathbf{x}) = \exp(-\sum_{i=1}^4 (x_i - \tfrac{1}{2})^2)$, where $\mathbf{x} = (x_1,\ldots,x_4) \in (0,1)^4$.

    \begin{theorem}
    \label{SPD:4}
    Let $\Omega = (0,1)^4$ and $h > 0$. Assume that $u$, $a$, $f$, and $g$ in \eqref{eq:PDE} are sufficiently smooth functions on $\overline{\Omega}$, and $|\nabla \tilde{a}|^2 + 2\partial^{2\vec{e}_i}\tilde{a} - \Delta \tilde{a} \le 0$ on $\Omega$, for $i=1,\ldots,4$, where $\tilde{a}:=-\ln a$. Consider the compact finite difference scheme $\mathcal{L}_h u_h = f_h$ for \eqref{eq:PDE} on $\Omega_h$ given by $\mathcal{L}_h u_h (\spt) := \sum_{p \in \SS} C_p (\spt) u_h(\spt + ph)$ for all $\spt \in \Omega_h$, where $C_p (\spt)$ and $f_h(\spt)$ are given in \eqref{Cp:fh:d4}.
    Then, for all sufficiently small $h$, the matrix $A_h$ is SPD, and the errors of the four-dimensional difference scheme satisfy
    \[
    \|\overrightarrow{u|_{\Omega_h}} - \vec{u}_h\|_{2,h} := h^{2} \|\overrightarrow{u|_{\Omega_h}} - \vec{u}_h\|_{2} \le C h^{4} \quad \text{and} \quad \|\overrightarrow{u|_{\Omega_h}} - \vec{u}_h\|_\infty \le C h^{4},
    \]
    where $C$ is a positive constant independent of $h$.
    \end{theorem}
    \begin{proof}
        See \cref{subsec:SPD:conv}.
    \end{proof}

    In what follows, we have our fifth and final main result concerning the SPD property of the linear systems arising from the discretization using the proposed $d$-dimensional for $d=2,3,4$, compact, symmetric, and 6th-order FDM assuming that $|\nabla a |^2 / a^2 - 2\Delta a / a$ is constant on $\Omega$. We also show their convergence rates.

    \begin{theorem}
        \label{SPD:sym:constant}
        Let $\Omega = (0,1)^d$ with $d=2,3,4$ and $h > 0$. Assume that $u$, $a$, $f$, and $g$ in \eqref{eq:PDE} are sufficiently smooth functions on $\overline{\Omega}$. Also, assume that $|\nabla a|^2 / a^2 - 2\Delta a / a$ is constant on $\Omega$. Consider the compact finite difference scheme $\mathcal{L}_h u_h = f_h$ for \eqref{eq:PDE} on $\Omega_h$ given by $\mathcal{L}_h u_h (\spt) := \sum_{p \in \SS} C_p (\spt) u_h(\spt + ph)$ for all $\spt \in \Omega_h$, where $C_p (\spt)$ and $f_h(\spt)$ are given in \eqref{stencil:dD:constant}.
        Then, the following statements hold:
        \begin{itemize}
            \item[(a)] For all sufficiently small $h$ and $d=2,3,4$, the matrix $A_h$ is SPD.
            \item[(b)] For all sufficiently small $h$, the errors of two, three, and four-dimensional finite difference schemes satisfy
            \[
             \|\overrightarrow{u|_{\Omega_h}} - \vec{u}_h\|_{2,h} := h^{d/2} \|\overrightarrow{u|_{\Omega_h}} - \vec{u}_h\|_{2} \le C h^{6} \quad \text{and} \quad \|\overrightarrow{u|_{\Omega_h}} - \vec{u}_h\|_\infty \le C h^{6},
            \]
            where $C$ is a positive constant independent of $h$.
        \end{itemize}
    \end{theorem}
    \begin{proof}
        The proof of this theorem is essentially identical to the proofs of \cref{SPD,SPD:4}.
    \end{proof}

    The construction of compact, symmetric, and 4th-order finite difference stencils for $d \ge 5$ and nonconstant $|\nabla a|^2 / a^2 - 2\Delta a / a$ is discussed in \cref{sec:any_dimension}. Note that for $d \ge 5$, the proposed schemes have off-diagonal stencil coefficients with opposite signs even for sufficiently small $h$, regardless of whether $|\nabla a|^2 / a^2 - 2\Delta a / a$ is constant on $\Omega$. Hence, we are unable to prove that the resulting linear systems are SPD and to establish theoretical convergence rates using the techniques employed in this paper.

\section{Numerical experiments}
    \label{sec:exp}
    Our goal now is to verify the convergence rates of our finite difference schemes through a series of numerical experiments. We use a uniform grid to discretize the domain $\Omega:=(0,1)^d$, where $d=1,2,3$. That is, we take $x_i = i h$, $y_j = j h$, and $z_k = k h$, where $i,j,k=1,\ldots,N-1$ and $h = 1/N$, as our grid. In addition, we define the following discrete errors for the 3D case
    \begin{align*}
    \|\overrightarrow{u|_{\Omega_h}} - \vec{u}_h\|_\infty & :=\max_{1\le i,j,k\le N-1} \left|(u_h)_{i,j,k}-u(x_i,y_j,z_k)\right|,\\
    \tfrac{\|\overrightarrow{u|_{\Omega_h}} - \vec{u}_h\|_2}{\|\overrightarrow{u|_{\Omega_h}}\|_2} & :=\sqrt{\tfrac{\sum_{i,j,k=1}^{N-1} \left|(u_h)_{i,j,k}-u(x_i,y_j,z_k)\right|^2}{\sum_{i,j,k=1}^{N-1}  \left|u(x_i,y_j,z_k)\right|^2}},
    \end{align*}
    where $(u_h)_{i,j,k}$ stands for the approximated solution on the grid point $(x_i, y_j, z_k)$ and $u$ is the exact solution. The errors for the 1D and 2D cases can be similarly computed. The convergence rates in the `order' columns are computed as $\log_2(E_h/E_{h/2})$, where $E_h$ denotes the corresponding error in the table.

    \subsection{Two examples in 1D}
    We use \cref{thm:1D} and \eqref{E12} in the discretization.
\begin{example}\label{Example:1:1D}
	\normalfont
	The functions  in \eqref{eq:PDE} with $\Omega:=(0,1)$ are given by
	\begin{align*}
		& a:= [2+\sin(1000x)]e^{x^2}[3+\cos(10x)],\qquad u:= e^{-x^3}[2+\cos(1000x)][4+\sin (13x)],
	\end{align*}
	and	$f,g$  are obtained by plugging the above functions into \eqref{eq:PDE}.
	The numerical results 
    are presented in \cref{Example:1:2:table:1D}. 
\end{example}

\begin{table}[htbp]
	\caption{Numerical results for \cref{Example:1:1D,Example:2:1D} using our 1D compact, symmetric, and 12th-order FDM.}
	\centering
	\resizebox{\textwidth}{!}{
			\begin{tabular}{c|c|c|c|c|c|c|c|c|c}
				\hline
				\multicolumn{10}{c}{ 1D compact, symmetric, and 12th-order FDM  } \\	
				\hline
				 \multicolumn{5}{c|}{\cref{Example:1:1D} with $\Omega=(0,1)$ } &
				 \multicolumn{5}{c}{\cref{Example:2:1D} with $\Omega=(0,1)$ } \\
				\cline{1-10}
				$h$&    $\frac{\|\vec{u}_{h}-\overrightarrow{u|_{\Omega_h}}\|_2}{\|\overrightarrow{u|_{\Omega_h}}\|_2}$    &  order &    $\|\vec{u}_{h}-\overrightarrow{u|_{\Omega_h}}\|_\infty$    &  order  & $h$&    $\frac{\|\vec{u}_{h}-\overrightarrow{u|_{\Omega_h}}\|_2}{\|\overrightarrow{u|_{\Omega_h}}\|_2}$      &order &  $\|\vec{u}_{h}-\overrightarrow{u|_{\Omega_h}}\|_\infty$   &order \\
				\hline
$\tfrac{1}{2^7}$  &  3.4520E+07  &    &  5.7515E+08  &    &    &    &    &    &  \\
$\tfrac{1}{2^8}$  &  2.3480E+04  &  10.52   &  3.5303E+05  &  10.67   &  $\tfrac{1}{2^9}$  &  3.2893E+04  &    &  3.6771E+04  &  \\
$\tfrac{1}{2^9}$  &  6.8676E-01  &  15.06   &  2.0243E+01  &  14.09   &  $\tfrac{1}{2^{10}}$  &  8.3446E+00  &  11.94   &  9.5414E+00  &  11.91 \\
$\tfrac{1}{2^{10}}$  &  3.8641E-04  &  10.80   &  8.0797E-03  &  11.29   &  $\tfrac{1}{2^{11}}$  &  8.2099E-04  &  13.31   &  8.8160E-04  &  13.40 \\
$\tfrac{1}{2^{11}}$  &  1.5819E-07  &  11.25   &  2.8859E-06  &  11.45   &  $\tfrac{1}{2^{12}}$  &  1.0228E-07  &  12.97   &  1.1290E-07  &  12.93 \\
$\tfrac{1}{2^{12}}$  &  1.4650E-11  &  13.40   &  1.6398E-10  &  14.10   &  $\tfrac{1}{2^{13}}$  &  4.0352E-10  &  7.99   &  4.2535E-10  &  8.05 \\
	\hline
&  \text{Average order:}  &  12.21   &  \text{Average order:}  &  12.32   &    &  \text{Average order:}  &  11.55   &  \text{Average order:}  &  11.57 \\
				\hline
	\end{tabular}}
	\label{Example:1:2:table:1D}
\end{table}

\begin{example}\label{Example:2:1D}
	\normalfont
	The functions  in \eqref{eq:PDE} with $\Omega:=(0,1)$ are given by
	\begin{align*}
		&a:= \ln(3x^3+5x^2+4)[3+\cos(1000x)],  \qquad  u:=\sin(40 \cdot 2^{x^2+2x+3}),
	\end{align*}
	and	$f,g$  are obtained by plugging the above functions into \eqref{eq:PDE}.
	The numerical results are presented in \cref{Example:1:2:table:1D}. 
\end{example}

In both 1D examples, the observed convergence rates are consistent with \cref{thm:1D} and \cref{SPD} item (b).

\subsection{Two examples in 2D with nonconstant $|\nabla a|^2/a^2-2\Delta a/a$ on $\Omega$}
We use \eqref{stencil:explicit:uf} and \eqref{2d:explicit} in the discretization.
\begin{example}\label{Example:1:2D}
	\normalfont
	The functions  in \eqref{eq:PDE} with $\Omega:=(0,1)^2$ are given by
	\begin{align*}
		& a:= 4+\cos \left( 5\pi \tanh (5x-3) \right)+\sin \left( 17.5 \tanh (4y-2) \right),\\
		& u:= e^{\sin(20\ln (3x^2+2y^2+1)) }\cos(20y),
	\end{align*}
	and	$f,g$  are obtained by plugging the above functions into \eqref{eq:PDE}.
	The numerical results are presented in \cref{Example:1:2:table:2D}.
    The function $a$ is obtained by modifying \cite[Example 7.5]{HWC999}.
\end{example}

\begin{example}\label{Example:2:2D}
	\normalfont
	The functions  in \eqref{eq:PDE} with $\Omega:=(0,1)^2$ are given by
	\begin{align*}
		& a:= 4+\sin(100x)\sin(100y),\qquad u:= \cos(100x)\cos(100y),
	\end{align*}
	and	$f,g$  are obtained by plugging the above functions into \eqref{eq:PDE}.
	The numerical results are presented in \cref{Example:1:2:table:2D}.
\end{example}	

In both 2D examples, the observed convergence rates are consistent with \cref{thm:max_order} item (a) and \cref{SPD} item (b). 

\begin{table}[htbp]
	\caption{Numerical results for \cref{Example:1:2D,Example:2:2D} using our 2D compact, symmetric, and 4th-order FDM.}
	\centering
	\resizebox{\textwidth}{!}{
			\begin{tabular}{c|c|c|c|c|c|c|c|c|c}
				\hline
				\multicolumn{10}{c}{ 2D compact, symmetric, and 4th-order FDM  } \\	
				\hline
				 \multicolumn{5}{c|}{\cref{Example:1:2D} with $\Omega=(0,1)^2$ } &
				 \multicolumn{5}{c}{\cref{Example:2:2D} with $\Omega=(0,1)^2$ } \\
				\cline{1-10}
				$h$&    $\frac{\|\vec{u}_{h}-\overrightarrow{u|_{\Omega_h}}\|_2}{\|\overrightarrow{u|_{\Omega_h}}\|_2}$    &  order &    $\|\vec{u}_{h}-\overrightarrow{u|_{\Omega_h}}\|_\infty$    &  order  & $h$&    $\frac{\|\vec{u}_{h}-\overrightarrow{u|_{\Omega_h}}\|_2}{\|\overrightarrow{u|_{\Omega_h}}\|_2}$      &order &  $\|\vec{u}_{h}-\overrightarrow{u|_{\Omega_h}}\|_\infty$   &order\\
				\hline
$\tfrac{1}{2^4}$ &  8.9368E-01 &   &  4.3616E+00 &   &  $\tfrac{1}{2^5}$ &  6.7687E+00 &   &  1.3146E+01 &  \\
$\tfrac{1}{2^5}$ &  2.7511E-02 &  5.02  &  1.6266E-01 &  4.74  &  $\tfrac{1}{2^6}$ &  9.2996E-02 &  6.19  &  1.1189E-01 &  6.88 \\
$\tfrac{1}{2^6}$ &  7.8214E-04 &  5.14  &  4.5310E-03 &  5.17  &  $\tfrac{1}{2^7}$ &  6.0184E-03 &  3.95  &  7.1631E-03 &  3.97 \\
$\tfrac{1}{2^7}$ &  4.7028E-05 &  4.06  &  2.6165E-04 &  4.11  &  $\tfrac{1}{2^8}$ &  3.8503E-04 &  3.97  &  4.5048E-04 &  3.99 \\
$\tfrac{1}{2^8}$ &  2.9450E-06 &  4.00  &  1.6378E-05 &  4.00  &  $\tfrac{1}{2^9}$ &  2.4257E-05 &  3.99  &  2.8574E-05 &  3.98 \\
$\tfrac{1}{2^9}$ &  1.8438E-07 &  4.00  &  1.0352E-06 &  3.98  &  $\tfrac{1}{2^{10}}$ &  1.5205E-06 &  4.00  &  1.7864E-06 &  4.00 \\
$\tfrac{1}{2^{10}}$ &  1.1535E-08 &  4.00  &  6.4784E-08 &  4.00  &  $\tfrac{1}{2^{11}}$ &  9.5137E-08 &  4.00  &  1.1166E-07 &  4.00 \\
				\hline
	\end{tabular}}
	\label{Example:1:2:table:2D}
\end{table}	

\subsection{Two examples in 3D with nonconstant $|\nabla a|^2/a^2-2\Delta a/a$ on $\Omega$}
We use \eqref{stencil:explicit:uf} and \eqref{3d:explicit} in the discretization.
\begin{example}\label{Example:2:3D}
	\normalfont
	The functions  in \eqref{eq:PDE} with $\Omega:=(0,1)^3$ are given by
	\begin{align*}
		& a:= 2+x^2+y^3+z^4,  \qquad u:=e^{x-2y+4z},
	\end{align*}
	and	$f,g$  are obtained by plugging the above functions into \eqref{eq:PDE}.
	The numerical results are presented in \cref{Example:1:2:table:3D}. 
\end{example}	

\begin{example}\label{Example:1:3D}
	\normalfont
	The functions  in \eqref{eq:PDE} with $\Omega:=(0,1)^3$ are given by
	\begin{align*}
		& a:= 2+\sin(5x-3y-3z),  \qquad u:=\cos(4x)\sin(4y)\cos(5z),
	\end{align*}
	and	$f,g$  are obtained by plugging the above functions into \eqref{eq:PDE}.
	The numerical results are presented in \cref{Example:1:2:table:3D}. 

\end{example}	

In both 3D examples, the observed convergence rates are consistent with \cref{thm:max_order} item (a) and \cref{SPD} item (b). 

\begin{table}[htbp]
	\caption{Numerical results for \cref{Example:1:3D,Example:2:3D} using our 3D compact, symmetric, and 4th-order FDM.}
	\centering
	\resizebox{\textwidth}{!}{
			\begin{tabular}{c|c|c|c|c|c|c|c|c|c}
				\hline
				\multicolumn{10}{c}{ 3D compact, symmetric, and 4th-order FDM  } \\	
				\hline
				 \multicolumn{5}{c|}{\cref{Example:2:3D} with $\Omega=(0,1)^3$ } &
				 \multicolumn{5}{c}{\cref{Example:1:3D} with $\Omega=(0,1)^3$ } \\
				\cline{1-10}
				$h$&    $\frac{\|\vec{u}_{h}-\overrightarrow{u|_{\Omega_h}}\|_2}{\|\overrightarrow{u|_{\Omega_h}}\|_2}$    &  order &    $\|\vec{u}_{h}-\overrightarrow{u|_{\Omega_h}}\|_\infty$    &  order  & $h$&    $\frac{\|\vec{u}_{h}-\overrightarrow{u|_{\Omega_h}}\|_2}{\|\overrightarrow{u|_{\Omega_h}}\|_2}$      &order &  $\|\vec{u}_{h}-\overrightarrow{u|_{\Omega_h}}\|_\infty$   &order \\
				\hline
                $\tfrac{1}{2}$ &  8.2602E-03 &   &  1.6294E+00 & &  $\tfrac{1}{2}$ &  1.6075E-01 &   &  3.1572E-01 &   \\
                $\tfrac{1}{2^2}$ &  1.2279E-03 &  2.75  &  1.4816E-01 &  3.46 & $\tfrac{1}{2^2}$ &  3.3301E-02 &  2.27  &  7.5349E-02 &  2.07 \\
                $\tfrac{1}{2^3}$ &  1.1451E-04 &  3.42  &  9.1878E-03 &  4.01 & $\tfrac{1}{2^3}$ &  3.5553E-03 &  3.23  &  5.6505E-03 &  3.74 \\
                $\tfrac{1}{2^4}$ &  8.8067E-06 &  3.70  &  5.7913E-04 &  3.99 & $\tfrac{1}{2^4}$ &  2.4658E-04 &  3.85  &  3.5319E-04 &  4.00 \\
                $\tfrac{1}{2^5}$ &  6.1362E-07 &  3.84  &  3.6321E-05 &  4.00 & $\tfrac{1}{2^5}$ &  1.6145E-05 &  3.93  &  2.2127E-05 &  4.00  \\
                $\tfrac{1}{2^6}$ &  4.0559E-08 &  3.92  &  2.2726E-06 &  4.00  & $\tfrac{1}{2^6}$ &  1.0330E-06 &  3.97  &  1.3805E-06 &  4.00 \\
				\hline
	\end{tabular}}
	\label{Example:1:2:table:3D}
\end{table}	

\cref{{ex:3d:cg}} reports the number of iterations required by the CGM to achieve a prescribed tolerance level. Due to the size of the linear system produced from discretizing this 3D problem, we compare the speed, accuracy, and convergence rates of solutions obtained from the built-in \verb|pcg| and \verb|mldivide| functions in MATLAB. For the former, we specify error tolerances according to the grid size so that the asymptotic 4th-order convergence rate remains unaffected. Computationally, we observe that the CGM is on average $6.4$ times faster than \verb|mldivide|. In particular, when $h=2^{-7}$, the CGM is around 29.9 times faster than \verb|mldivide|. These findings confirm well-known facts: the CGM performs well for SPD matrices, and \verb|mldivide|, which involves a direct solver, is computationally more expensive than the CGM for large sparse linear systems.

\begin{table}[htbp]
	\caption{CGM results for \cref{Example:1:3D}. The built-in \texttt{pcg} and \texttt{mldivide} MATLAB functions are used to solve the linear system $A_h \vec{u}_h=\vec{f}_h$ obtained from discretizing \eqref{eq:PDE} using a compact, symmetric, and 4th-order 3D FDM. Below, we define $R:=\frac{ \|A_h \vec{u}_h-\vec{f}_h\|_2 }{\|\vec{f}_h\|_2}$. Furthermore, we denote by $I$ the number of iterations of the conjugate gradient method, and $T_{CG}, T_{D}$ as the computational times (in seconds) for \texttt{pcg} and \texttt{mldivide} respectively to solve the linear system.}
	\centering
	\resizebox{\textwidth}{!}{
			\begin{tabular}{c|c|c|c|c|c|c|c|c|c|c}
				\hline
				\multicolumn{1}{c|}{} &
				\multicolumn{10}{c}{ 3D compact, symmetric, and 4th-order FDM } \\ 	
				\hline
				\multicolumn{1}{c|}{} &
				\multicolumn{6}{c|}{\texttt{pcg}} & 	 
				 \multicolumn{3}{c|}{\texttt{mldivide}} &
				\multicolumn{1}{c}{}  \\ 	
				\hline
				$h$  &   $\|\vec{u}_{h}-\overrightarrow{u|_{\Omega_h}}\|_\infty$    &order &  \text{tolerance} & $R$    & $I$ &  $T_{CG}$     & $\|\vec{u}_{h}-\overrightarrow{u|_{\Omega_h}}\|_\infty$    &order & $T_{D}$ & $T_{D}/T_{CG}$  \\
				\hline
                $\tfrac{1}{2}$ & 5.97E-01 &  & $10^{-2}$ & 9.2E-17 & 1  & 4.9E-03 & 5.97E-01 &  & 4.8E-04 & 0.1\\
                $\tfrac{1}{2^2}$ & 1.27E-01 & 2.2 & $10^{-2}$ & 4.5E-03 & 4  & 4.4E-03 & 1.21E-01 & 2.3 & 7.9E-03 & 1.8\\
                $\tfrac{1}{2^3}$ & 2.13E-02 & 2.6 & $10^{-2}$ & 8.3E-03 & 7  & 1.7E-03 & 6.11E-03 & 4.3 & 5.4E-03 & 3.3\\
                $\tfrac{1}{2^4}$ & 3.44E-04 & 6.0 & $10^{-4}$ & 7.6E-05 & 36  & 1.0E-02 & 3.59E-04 & 4.1 & 3.6E-02 & 3.6\\
                $\tfrac{1}{2^5}$ & 2.09E-05 & 4.0 & $10^{-6}$ & 9.8E-07 & 99  & 4.3E-01 & 2.22E-05 & 4.0 & 7.8E-01 & 1.8\\
                $\tfrac{1}{2^6}$ & 1.32E-06 & 4.0 & $10^{-8}$ & 9.6E-09 & 244  & 8.0E+00 & 1.38E-06 & 4.0 & 3.5E+01 & 4.4\\
                $\tfrac{1}{2^7}$ & 8.66E-08 & 3.9 & $10^{-10}$ & 9.9E-11 & 590  & 1.5E+02 & 8.62E-08 & 4.0 & 4.6E+03 & 29.9\\		
				\hline
	\end{tabular}}
	\label{ex:3d:cg}
\end{table}

\subsection{One example in 2D and one example in 3D with constant $|\nabla a|^2/a^2-2\Delta a/a$ on $\Omega$}
\label{sec:exp:constant:Qa}
We use \eqref{stencil:explicit:uf} and \eqref{stencil:dD:constant} in the discretization.

\begin{example}\label{Example:constant:2D}
    \normalfont
    The functions in \eqref{eq:PDE} with $\Omega:=(0,1)^2$ are given by
    \begin{align*}
        &a:=\left[\cosh(2x+y)+2\cosh(x-2y)+\cosh x\cosh(2y)\right]^2,
        \qquad u:=e^{3x-4y},
    \end{align*}
    and $f,g$ are obtained by plugging the above functions into \eqref{eq:PDE}.
    In this example, $|\nabla a|^2/a^2-2\Delta a/a=-20$ on $\Omega$.
    The numerical results are presented in \cref{Example:constant:table:2D:3D}.
\end{example}

\begin{example}\label{Example:constant:3D}
    \normalfont
    The functions in \eqref{eq:PDE} with $\Omega:=(0,1)^3$ are given by
    \[
        a:=\left(\frac{\sin(\pi r)}{\pi r}\right)^2,\quad
        u:=\frac{\pi r \sin(\pi x)\sin(2\pi y)\sin(3\pi z)}{\sin(\pi r)}, \quad 
        r:=\sqrt{(x-\tfrac12)^2+(y-\tfrac12)^2+(z-\tfrac12)^2},
    \]
    and $f,g$ are obtained by plugging the above functions into \eqref{eq:PDE}. In this example, $|\nabla a|^2/a^2-2\Delta a/a=4\pi^2$ on $\Omega$.
    The numerical results are presented in \cref{Example:constant:table:2D:3D}.
\end{example}

In the above 2D and 3D examples, the observed convergence rates are consistent with \cref{thm:max_order} item (b) and \cref{SPD:sym:constant} item (b). 

\begin{table}[htbp]
    \caption{Numerical results for \cref{Example:constant:2D,Example:constant:3D} using our 2D and 3D compact, symmetric, and 6th-order FDMs.}
    \centering
    \resizebox{\textwidth}{!}{
    \begin{tabular}{c|c|c|c|c|c|c|c|c|c}
        \hline
        \multicolumn{10}{c}{2D and 3D compact, symmetric, and 6th-order FDMs}\\
        \hline
        \multicolumn{5}{c|}{\cref{Example:constant:2D} with $\Omega=(0,1)^2$} &
        \multicolumn{5}{c}{\cref{Example:constant:3D} with $\Omega=(0,1)^3$}\\
        \cline{1-10}
        $h$ & $\frac{\|\vec{u}_{h}-\overrightarrow{u|_{\Omega_h}}\|_2}{\|\overrightarrow{u|_{\Omega_h}}\|_2}$ & order & $\|\vec{u}_{h}-\overrightarrow{u|_{\Omega_h}}\|_\infty$ & order &
        $h$ & $\frac{\|\vec{u}_{h}-\overrightarrow{u|_{\Omega_h}}\|_2}{\|\overrightarrow{u|_{\Omega_h}}\|_2}$ & order & $\|\vec{u}_{h}-\overrightarrow{u|_{\Omega_h}}\|_\infty$ & order\\
        \hline
$\tfrac{1}{2^{2}}$ & 1.8932E-03 &  & 5.3304E-03 &  & $\tfrac{1}{2^{2}}$ & 1.5813E-01 &  & 1.7564E-01 &  \\
$\tfrac{1}{2^{3}}$ & 2.1474E-05 & 6.46 & 9.5356E-05 & 5.80 & $\tfrac{1}{2^{3}}$ & 1.6494E-03 & 6.58 & 2.1837E-03 & 6.33 \\
$\tfrac{1}{2^{4}}$ & 2.7521E-07 & 6.29 & 1.6446E-06 & 5.86 & $\tfrac{1}{2^{4}}$ & 2.3174E-05 & 6.15 & 3.0682E-05 & 6.15 \\
$\tfrac{1}{2^{5}}$ & 3.8706E-09 & 6.15 & 2.5960E-08 & 5.99 & $\tfrac{1}{2^{5}}$ & 3.5252E-07 & 6.04 & 4.8759E-07 & 5.98 \\
$\tfrac{1}{2^{6}}$ & 5.7311E-11 & 6.08 & 4.0667E-10 & 6.00 & $\tfrac{1}{2^{6}}$ & 5.4712E-09 & 6.01 & 7.5676E-09 & 6.01 \\
$\tfrac{1}{2^{7}}$ & 8.6729E-13 & 6.05 & 6.3256E-12 & 6.01 & $\tfrac{1}{2^{7}}$ & 8.5344E-11 & 6.00 & 1.1812E-10 & 6.00 \\
        \hline
    \end{tabular}}
    \label{Example:constant:table:2D:3D}
\end{table}

In addition to the convergence summaries, \cref{tab:smallest-eigenvalues} reports the smallest eigenvalues of some coefficient matrices from Examples 1--8 assembled using the stencils in \cref{sec:stencils}. Each smallest eigenvalue is positive, and the ratios approach 4 as the mesh size is halved in each direction. Together with \eqref{sym:cond}, these results are consistent with \cref{thm:1D,SPD,SPD:sym:constant}, which establish the SPD property for all sufficiently small grid sizes.

\begin{table}[htbp]
	\caption{Smallest eigenvalues, $\lambda_{\min}$, and their consecutive ratios for Examples 1--8.}
	\centering
	\resizebox{\textwidth}{!}{
		\begin{tabular}{c|c|c|c|c|c|c|c|c|c|c|c}
			\hline
			\multicolumn{3}{c|}{\cref{Example:1:1D}} &
			\multicolumn{3}{c|}{\cref{Example:2:1D}} &
			\multicolumn{3}{c|}{\cref{Example:1:2D}} &
			\multicolumn{3}{c}{\cref{Example:2:2D}} \\
			\cline{1-12}
			$h$ & $\lambda_{\min}$ & Ratio &
			$h$ & $\lambda_{\min}$ & Ratio &
			$h$ & $\lambda_{\min}$ & Ratio &
			$h$ & $\lambda_{\min}$ & Ratio \\
			\hline
			$\tfrac{1}{2^{8}}$
			& 8.3133E-04
			& --
			& $\tfrac{1}{2^{9}}$
			& 1.9095E-04
			& --
			& $\tfrac{1}{2^{4}}$
			& 4.2337E-01
			& --
			& $\tfrac{1}{2^{5}}$
			& 1.1533E-01
			& -- \\

			$\tfrac{1}{2^{9}}$
			& 2.6983E-04
			& 3.08
			& $\tfrac{1}{2^{10}}$
			& 4.7738E-05
			& 4.00
			& $\tfrac{1}{2^{5}}$
			& 1.0623E-01
			& 3.99
			& $\tfrac{1}{2^{6}}$
			& 2.8688E-02
			& 4.02 \\

			$\tfrac{1}{2^{10}}$
			& 6.7460E-05
			& 4.00
			& $\tfrac{1}{2^{11}}$
			& 1.1934E-05
			& 4.00
			& $\tfrac{1}{2^{6}}$
			& 2.6590E-02
			& 4.00
			& $\tfrac{1}{2^{7}}$
			& 7.1728E-03
			& 4.00 \\

			$\tfrac{1}{2^{11}}$
			& 1.6865E-05
			& 4.00
			& $\tfrac{1}{2^{12}}$
			& 2.9836E-06
			& 4.00
			& $\tfrac{1}{2^{7}}$
			& 6.6496E-03
			& 4.00
			& $\tfrac{1}{2^{8}}$
			& 1.7934E-03
			& 4.00 \\

			$\tfrac{1}{2^{12}}$
			& 4.2162E-06
			& 4.00
			& $\tfrac{1}{2^{13}}$
			& 7.4590E-07
			& 4.00
			& $\tfrac{1}{2^{8}}$
			& 1.6625E-03
			& 4.00
			&
			&
			&  \\
			\hline
			            \multicolumn{3}{c|}{\cref{Example:2:3D}} &
            \multicolumn{3}{c|}{\cref{Example:1:3D}} &
            \multicolumn{3}{c|}{\cref{Example:constant:2D}} &
            \multicolumn{3}{c}{\cref{Example:constant:3D}}\\
            \cline{1-12}
            $h$ & $\lambda_{\min}$ & Ratio &
            $h$ & $\lambda_{\min}$ & Ratio &
            $h$ & $\lambda_{\min}$ & Ratio &
            $h$ & $\lambda_{\min}$ & Ratio\\
            \hline
$\tfrac{1}{2}$ & 2.5282E+00 & -- & $\tfrac{1}{2}$ & 1.6731E+00 & -- & $\tfrac{1}{2^{2}}$ & 4.9295E+01 & -- & $\tfrac{1}{2^{2}}$ & 7.0218E-01 & -- \\
$\tfrac{1}{2^{2}}$ & 9.7450E-01 & 2.59 & $\tfrac{1}{2^{2}}$ & 6.1153E-01 & 2.74 & $\tfrac{1}{2^{3}}$ & 1.3491E+01 & 3.65 & $\tfrac{1}{2^{3}}$ & 1.9777E-01 & 3.55 \\
$\tfrac{1}{2^{3}}$ & 2.7489E-01 & 3.55 & $\tfrac{1}{2^{3}}$ & 1.8195E-01 & 3.36 & $\tfrac{1}{2^{4}}$ & 3.4454E+00 & 3.92 & $\tfrac{1}{2^{4}}$ & 5.1358E-02 & 3.85 \\
$\tfrac{1}{2^{4}}$ & 7.0873E-02 & 3.88 & $\tfrac{1}{2^{4}}$ & 4.7328E-02 & 3.84 & $\tfrac{1}{2^{5}}$ & 8.6590E-01 & 3.98 & $\tfrac{1}{2^{5}}$ & 1.2969E-02 & 3.96 \\
$\tfrac{1}{2^{5}}$ & 1.7856E-02 & 3.97 & $\tfrac{1}{2^{5}}$ & 1.1947E-02 & 3.96 & $\tfrac{1}{2^{6}}$ & 2.1676E-01 & 3.99 & $\tfrac{1}{2^{6}}$ & 3.2507E-03 & 3.99 \\
 &  &  &  &  &  & $\tfrac{1}{2^{7}}$ & 5.4208E-02 & 4.00 & $\tfrac{1}{2^{7}}$ & 8.1318E-04 & 4.00 \\
            \hline
		\end{tabular}}
	\label{tab:smallest-eigenvalues}
\end{table}

\section{Compact, symmetric, and fourth-order $d$-dimensional FDMs, $d\ge 2$}
    \label{sec:construct:dD}

    In this section, we provide the main steps for constructing finite difference schemes satisfying the properties described in \Cref{thm:max_order}, thereby recovering the specific examples provided in \cref{sec:stencils}. We shall elaborate on an earlier remark on the non-uniqueness of the finite difference schemes in dimensions greater than 1. The supporting technical details are deferred to \Cref{sec:details}.

    In the second half of this section, we shall show that $d$-dimensional compact, symmetric, and 4th-order finite difference schemes with $d \ge 3$ can be constructed by combining the one and two-dimensional finite difference schemes thanks to the bilinear structure of the grid function $f_h$ in \eqref{eq:f_h_any_dim}.

    \subsection{The 2D and 3D cases}
    \label{sec:basis}

    One key observation that we utilize is
    \begin{equation}
        \label{eq:u_taylor_int}
        u(\spt+ ph)
        = \sum_{\bm{\ell} \in \LN_{M + 1}} \sum_{k = |\bm{\ell}|}^{M + 1}
        A^k_{\bm{\ell}} (p) h^k \cdot \partial^{\bm{\ell}} u(\spt) + F(p) + \bo_{a,u}(h^{M + 2}), \quad \forall\spt \in \Omega_h,
    \end{equation}
    where $A^k_{\bm{\ell}} (p)$, $F(p)$ are explicitly known quantities obtained from Taylor expansions that depend on $a$ and its partial derivatives at $\spt$, $F(p)$ also depends on $f$ and its partial derivatives at $\spt$, and
    \begin{equation} \label{PiM}
        \LN_M := \{ \bm{\ell} \in \NN^d: \ell_1 \le 1, \, |\bm{\ell}| \le M \}, \ M \in \NN
        \quad \text{with} \quad
        \bm{\ell} := (\ell_1,\ldots,\ell_d).
    \end{equation}
    All the coefficients $c_{p,k}$ in \Cref{thm:max_order} can be determined by solving a linear system, which has special cascade structures and involves the quantities $A^k_{\bm{\ell}} (p)$ in \eqref{eq:u_taylor_int}. Moreover, with the aid of symbolic computation, we observe that there are free parameters available. As a result, for all $\spt \in \Omega_h$, we can write
    \begin{align} \label{eq:C_p_2D:0}
    C_p(\spt) & = \sum_{m = -1}^{M} \sum_{k = 1}^{K_m} \kappa_{m, k} C_p^{[m], k}(\spt) h^{M - m}, \quad \forall  p \in \SS,
    \end{align}
    where $K_m$ is the number of available basis stencil coefficients obtained via symbolic computation, $\kappa_{m, k}$ are some constants for all $\spt \in \Omega_h$, $C_p^{[m], k}(\spt)$ are basis stencil coefficients, and $\mathcal{S} = \{-1,0,1\}^d$ with $d=2,3$. Furthermore, we can compute the right-hand basis grid functions, $f_h^{[m], k}(\spt)$, which are defined as
    \[
    f_h^{[m], k}(\spt) := \sum_{p \in \SS} C_p^{[m], k}(\spt) F(p), \quad \forall \spt \in \Omega_h,
    \]
    where $F(p)$ comes from Taylor expansions done in \eqref{eq:u_taylor_int}.

    Now, the $M$th-order compact finite difference schemes take the following form
    \begin{equation}
    \label{eq:assemble}
    \begin{aligned}
    & \sum_{p \in \mathcal{S}} \sum_{m = -1}^{M} \sum_{k = 1}^{K_m} \kappa_{m, k} C_p^{[m], k}(\spt) h^{M - m} u(\spt + ph) \\
    & \hspace{2cm} = \sum_{m = -1}^{M} \sum_{k = 1}^{K_m} \kappa_{m, k} f_h^{[m], k}(\spt) h^{M - m} 
    + \bo_{a,u}(h^{M+2}), \quad \forall \spt \in \Omega_h.
    \end{aligned}
    \end{equation}
    By \eqref{sym:cond}, to guarantee that the linear system produced by the discretization is symmetric, we further enforce the following conditions on the stencil coefficients in \eqref{eq:assemble}
    \be \label{cpmk:sym}
    \sum_{m = -1}^{M} \sum_{k = 1}^{K_m} \kappa_{m, k} C_p^{[m], k}(\spt) h^{M - m}
    = \sum_{m = -1}^{M} \sum_{k = 1}^{K_m} \kappa_{m, k} C_{-p}^{[m], k}(\spt + ph) h^{M - m},
    \ee
    for all $\spt \in \Omega_h$ and $p \in \mathcal{S}$, where each $C_p^{[m], k}$ is evaluated at the point $\spt + ph / 2$.

    Next, we shall present the explicit forms of all possible $C_p^{[m], k}(\spt)$ for $p \in \SS \bs \{\vec{0}\}$ with $C_{\vec{0}}^{[m], k} (\spt)$ being uniquely determined by $C_{\vec{0}}^{[m], k} (\spt) = -\sum_{p \in \SS \bs \{\vec{0}\}} C_p^{[m], k} (\spt)$ according to \Cref{lem:C_p_constraint_int}(b), and all possible $f_h^{[m], k}(\spt)$. Afterwards, we shall discuss some restrictions on $\kappa_{m,k}$. In what follows, we define $\tilde{a}:=-\ln a$ and $\tilde{f}:=-f/a$.

    For $d = 2$, we have $\SS= \{(0,0), \pm (1,0), \pm(0,1), \pm(1,1), \pm (1,-1)\}$, and
    \be \label{KM:2D}
    \begin{aligned}
    & K_{-1}=8, \quad K_{0}=6, \quad K_{1}=4, \quad K_{2}=2, \quad K_{3} = K_{4} = K_{5} = K_{6} =1, \quad \\ 
    & K_{M} = 0 \quad \forall M \ge 7.
    \end{aligned}
    \ee
    The basis stencil coefficients can be written as
    \begin{equation}
        \label{eq:C_p_2D}
        C_p^{[m],k} (\spt) =
        \begin{cases}
            (a \Phi^\mathsf{H} \pm a \widetilde{\Phi}^\mathsf{H}) (\spt + ph/2), & p = \pm (1, 0), \\
            (a \Phi^\mathsf{V} \pm a \widetilde{\Phi}^\mathsf{V}) (\spt + ph/2), & p = \pm (0, 1), \\
            (a \Phi^\mathsf{D} \pm a \widetilde{\Phi}^\mathsf{D}) (\spt + ph/2), & p = \pm (1, 1), \\
            (a \Phi^\mathsf{A} \pm a \widetilde{\Phi}^\mathsf{A}) (\spt + ph/2), & p = \pm (1, -1), 
        \end{cases}
    \end{equation}
    where superscripts $``\mathsf{H}"$ (horizontal), $``\mathsf{V}"$ (vertical), $``\mathsf{D}"$ (diagonal), and $``\mathsf{A}"$ (antidiagonal) emphasize the stencil's location. The sign in between $\Phi$ and $\tilde{\Phi}$ matches the sign of $p$. For example, if $p = (1,0)$, then $C_p^{[m],k} (\spt) = (a \Phi^\mathsf{H} + a \widetilde{\Phi}^\mathsf{H}) (\spt + ph/2)$.

    Below, we list the values that these $\Phi$ and $\tilde{\Phi}$ with various superscripts as well as the right-hand basis grid functions take depending on $m$ and $k$:
    \begin{align*}
        & (\Phi^\mathsf{H}, \Phi^\mathsf{V}, \Phi^\mathsf{D}, \Phi^\mathsf{A}, \widetilde{\Phi}^\mathsf{H}, \widetilde{\Phi}^\mathsf{V}, \widetilde{\Phi}^\mathsf{D}, \widetilde{\Phi}^\mathsf{A},f^{[m],k}_h) = \\
        & \quad
        \begin{cases}
        -\vec{e}_k, & \text{if} \; m=-1,0,\; \text{and} \; k=1,\ldots,4,\\
        \vec{e}_k, & \text{if} \; m=-1,\; \text{and} \; k=5,\ldots,8,\\
        2\vec{e}_5-\vec{e}_7+\vec{e}_8, & \text{if} \; m = 0, \, k = 5 \mbox{; or, } m = 1, \, k = 3,\\
        2\vec{e}_6-\vec{e}_7-\vec{e}_8, & \text{if} \; m = 0, \, k = 6 \mbox{; or, } m = 1, \, k = 4,\\
        -\vec{e}_{2k-1}-\vec{e}_{2k} + k h^2 f \vec{e}_9, & \text{if} \; m = 1, 2,\, \mbox{ and } \, k = 1,2,
        \end{cases}
    \end{align*}
    where $\vec{e}_k$, $k=1,\ldots,9$ are the standard unit vector basis of $\mathbb{R}^9$.

    If $m=3,4$ and $k=1$, then
    \begin{align*}
        \Phi^\mathsf{H} &= -4 + \frac{1}{4} h^2 \left( |\nabla \tilde{a}|^2 + \tilde{a}_{xx} - \tilde{a}_{yy} \right), \quad
        \Phi^\mathsf{V} = -4 + \frac{1}{4} h^2 \left( |\nabla \tilde{a}|^2 - \tilde{a}_{xx} + \tilde{a}_{yy} \right), \\
        \Phi^\mathsf{D} & = -1 + \frac{1}{4} h^2 \tilde{a}_{xy}, \hspace{3cm} \Phi^\mathsf{A}  = -1 - \frac{1}{4} h^2 \tilde{a}_{xy}, \\
        \widetilde{\Phi}^\mathsf{H} & = \widetilde{\Phi}^\mathsf{V} = \widetilde{\Phi}^\mathsf{D} = \widetilde{\Phi}^\mathsf{A} = 0,
        \hspace{2cm}  f^{[m],k}_h = -a \left( 6 h^2 \tilde{f} + \frac{1}{2} h^4 \left( \Delta \tilde{f} - \nabla \tilde{a} \cdot \nabla \tilde{f} \right) \right).
    \end{align*}

    Finally, if $m=5,6$ and $k=1$, then
    \begin{align*}
        \Phi^\mathsf{H} & = -4 + h^2 \left( - \frac{11}{60} (\eta_1 + \eta_2) + \frac{1}{2} \tilde{a}_{xx} \right) + \frac{1}{960} h^4 \Big( \! - \! 8 |\nabla \tilde{a}|^4 + 26 \nabla \tilde{a} \cdot \nabla \Delta \tilde{a}   \\
        & \quad - 16 \Delta^2 \tilde{a} + 14 (\Delta \tilde{a})^2 + 7 |\nabla \tilde{a}|^2 \Delta \tilde{a} + 10 \tilde{a}_{xxyy} - 22 \tilde{a}_{xx} \tilde{a}_{yy} + 5 \Delta \tilde{a} \eta_5 - 7 \Delta \eta_5 \\
        & \quad - 11 |\nabla \tilde{a}|^2 \eta_5 - 2 \nabla \tilde{a} \cdot \nabla \eta_5 \Big), \\
        \Phi^\mathsf{V} & = -4 + h^2 \left( - \frac{11}{60} (\eta_1 + \eta_2) + \frac{1}{2} \tilde{a}_{yy} \right) + \frac{1}{960} h^4 \Big( \! - \! 8 |\nabla \tilde{a}|^4 + 26 \nabla \tilde{a} \cdot \nabla \Delta \tilde{a}   \\
        & \quad - 16 \Delta^2 \tilde{a} + 14 (\Delta \tilde{a})^2 + 7 |\nabla \tilde{a}|^2 \Delta \tilde{a} + 10 \tilde{a}_{xxyy} - 22 \tilde{a}_{xx} \tilde{a}_{yy} - 5 \Delta \tilde{a} \eta_5 + 7 \Delta \eta_5  \\
        & \quad + 11 |\nabla \tilde{a}|^2 \eta_5 + 2 \nabla \tilde{a} \cdot \nabla \eta_5 \Big), \\
        \Phi^\mathsf{D} & = -1 + h^2 \left( \frac{1}{120} \eta_4 + \frac{1}{4} \tilde{a}_{xy} \right) + \frac{1}{1440} h^4 \big( \partial_{xy} (4(\eta_1 + \eta_2) + \eta_4) - 3 \tilde{a}_{xy} \eta_4 \big), \\
        \Phi^\mathsf{A} & = -1 + h^2 \left( \frac{1}{120} \eta_4 - \frac{1}{4} \tilde{a}_{xy} \right) - \frac{1}{1440} h^4 \big( \partial_{xy} (4(\eta_1 + \eta_2) + \eta_4) - 3 \tilde{a}_{xy} \eta_4 \big), \\
        \widetilde{\Phi}^\mathsf{H} & = \frac{1}{40} h^3 \partial_x (\eta_1 + \eta_2) + \frac{1}{960} h^5 \big( \partial_{xxx} (\eta_1 + \eta_2) - 3 \tilde{a}_{xx} \partial_x (\eta_1 + \eta_2) \big), \\
        \widetilde{\Phi}^\mathsf{V} & = \frac{1}{40} h^3 \partial_y (\eta_1 + \eta_2) + \frac{1}{960} h^5 \big( \partial_{yyy} (\eta_1 + \eta_2) - 3 \tilde{a}_{yy} \partial_y (\eta_1 + \eta_2) \big), \\
        \widetilde{\Phi}^\mathsf{D} & = \widetilde{\Phi}^\mathsf{A} = 0, \\
        f^{[m],k}_h & = -a \bigg[ 6 h^2 \tilde{f} + \frac{1}{2} h^4 \left( \Delta \tilde{f} - \nabla \tilde{a} \cdot \nabla \tilde{f} + \frac{1}{10} (\eta_1 + \eta_2) \tilde{f} \right) - \frac{1}{240} h^6 \bigg( \! - \! 4 (\tilde{f}_{xxxx}  \\
        & \quad + 4 \tilde{f}_{xxyy} + \tilde{f}_{yyyy}) + 8 (\tilde{a}_x \tilde{f}_{xxx} + 2 \tilde{a}_y \tilde{f}_{xxy} + 2 \tilde{a}_x \tilde{f}_{xyy} + \tilde{a}_y \tilde{f}_{yyy})  \\
        & \quad + \big( (\eta_1 - \eta_2) (\tilde{f}_{xx} - \tilde{f}_{yy}) + 16 \eta_3 \tilde{f}_{xy} \big) + \big( \tilde{a}_x (3\eta_1 + \eta_2) - 2 \partial_x (\eta_1 - \eta_2) \big) \tilde{f}_x  \\
        & \quad + \big( \tilde{a}_y (\eta_1 + 3\eta_2) + 2 \partial_y (\eta_1 - \eta_2) \big) \tilde{f}_y + \Big( \frac{1}{4} (\eta_1^2 + \eta_2^2) + \tilde{a}_x \partial_x (2\eta_1 + \eta_2) \\
        & \quad + \tilde{a}_y \partial_y (\eta_1 + 2\eta_2) - \partial_{xx} (2\eta_1 + \eta_2) - \partial_{yy} (\eta_1 + 2\eta_2) \Big) \tilde{f} \bigg) \bigg],
    \end{align*}
    where $\eta_1 = 2\tilde{a}_{xx} - \tilde{a}_x^2$, $\eta_2 = 2\tilde{a}_{yy} - \tilde{a}_y^2$, $\eta_3 = 2 \tilde{a}_{xy} - \tilde{a}_x \tilde{a}_y$, $\eta_4 = \Delta \tilde{a} + 7 |\nabla \tilde{a}|^2$, $\eta_5 = \tilde{a}_{xx} - \tilde{a}_{yy}$.

    A few observations are in order. The basis stencil coefficients $C_p^{[m],k} (\spt)$ are symmetric if $\widetilde{\Phi}^\mathsf{H}$, $\widetilde{\Phi}^\mathsf{V}$, $\widetilde{\Phi}^\mathsf{D}$, and $\widetilde{\Phi}^\mathsf{A}$ are all equal to zeros. On the other hand, the basis stencil coefficients $C_p^{[m],k} (\spt)$ are anti-symmetric if $\Phi^\mathsf{H}$, $\Phi^\mathsf{V}$, $\Phi^\mathsf{D}$, and $\Phi^\mathsf{A}$ are all equal to zeros; i.e., $C_p^{[m], k}(\spt) = -C_{-p}^{[m], k}(\spt + ph)$ for all $\spt \in \Omega_h$ and $p \in \mathcal{S}$. 
    The basis stencil coefficients $C^{[5],1}_p, C^{[6],1}_p$ are neither symmetric nor anti-symmetric, but they can be decomposed into symmetric and anti-symmetric parts by taking even and odd powers of $h$.

    The free parameters $\kappa_{m, k}$ in \eqref{eq:C_p_2D:0} can be chosen to help us satisfy \eqref{cpmk:sym} given the symmetric/anti-symmetric property of the basis stencil coefficients. One natural choice is to set $\kappa_{m, k} = 0$ for all basis stencil coefficients that are not symmetric (this also includes the anti-symmetric ones). At the same time, this freedom allows us to potentially minimize the leading truncation error term, which is expected to produce a smaller error for a given grid size. Moreover, if the free parameters $\kappa_{m, k}$ are also chosen so that $- C_p(\spt) \ge c_a >0$ for all $ p \in \mathcal{S} \backslash \{\vec{0}\}$ with a fixed positive constant $c_a$ that only depends on the diffusion term $a$, then we can prove that the scheme produces SPD linear systems and achieves 4th-order convergence for all sufficiently small grid sizes by using the arguments in the proof of \cref{SPD}.

    For a general $a$, we can take $M \le 4$ in \eqref{eq:C_p_2D} to ensure that the resulting finite difference scheme is symmetric. Meanwhile, if we further assume that $2 \Delta \tilde{a} - |\nabla \tilde{a}|^2$ is constant, then we can take $M \le 6$ in \eqref{eq:C_p_2D} to guarantee the symmetry. Under this assumption, the basis stencil coefficients become symmetric, since $\tilde{\Phi}^\mathsf{A}$, $\tilde{\Phi}^\mathsf{D}$ are already zeros to begin with, and $\tilde{\Phi}^\mathsf{H}$, $\tilde{\Phi}^\mathsf{V}$ are now also zeros for $M = 5$, $6$. In particular, if $a$ is constant in \eqref{eq:PDE} and $M=6$, we recover the 6th-order 2D finite difference scheme, which is widely known in the literature \cite{kyei2010family,zhai2013family,zhai2014new}.  If symmetry is not required, we can take $M \le 6$ in \eqref{eq:C_p_2D} and the only restriction imposed on $\kappa_{m, k}$ in \eqref{eq:C_p_2D} is given in \Cref{lem:C_p_constraint_int}.

    For $d = 3$, we have
    \begin{align*}
        & \SS = \{(0,0,0), \pm (1,0,0), \pm(0,1,0), \pm(0,0,1), \pm(1,1,0), \pm(1,-1,0),  \\
        \nonumber
        & \qquad \pm(1,0,1), \pm(1,0,-1), \pm (0,1,1), \pm (0,1,-1), \pm (1,1,1), \pm (-1,1,1),\\
        & \qquad \pm (1,-1,1), \pm (1,1,-1)\},
    \end{align*}
    and the number of available basis stencil coefficients are
    \begin{align}
        \label{KM:3D:1}
        & K_{-1}=26, \quad K_{0}=23, \quad K_{1}=18, \quad K_{2}=11, \quad K_{3} = 5, \quad K_{4} = 2,\\
        \nonumber
        & K_{5} = K_{6} =1, \quad K_{M} = 0 \quad \forall M \ge 7;
    \end{align}
    meanwhile, if we further exclude 8 corner points from $\SS$, then we have
    \be \label{KM:3D:2}
    \begin{aligned}
    & K_{-1}=18, \quad K_{0}=15, \quad K_{1}=10, \quad K_{2}=4, \quad K_{3} = K_{4} = 1, \quad K_{5} = K_{6} = 0,\\
    & K_{M} = 0 \quad \forall M \ge 7,
    \end{aligned}
    \ee
    and no additional points can be removed from $\SS$ while keeping $K_4 > 0$. Similar to \eqref{eq:C_p_2D}, the stencil coefficients can be written as
    \begin{equation}
        \label{eq:C_p_3D}
        C_p^{[m],k} (\spt) =
        \begin{cases}
            (a \Phi^\mathsf{F} \pm a \widetilde{\Phi}^\mathsf{F}) (\spt + ph/2),^{\phantom{\mathsf{ABCD}}} & p = \pm (1, 0, 0), \\
            (a \Phi^\mathsf{R} \pm a \widetilde{\Phi}^\mathsf{R}) (\spt + ph/2), & p = \pm (0, 1, 0), \\
            (a \Phi^\mathsf{U} \pm a \widetilde{\Phi}^\mathsf{U}) (\spt + ph/2), & p = \pm (0, 0, 1), \\
            (a \Phi^\mathsf{FR} \pm a \widetilde{\Phi}^\mathsf{FR}) (\spt + ph/2), & p = \pm (1, 1, 0), \\
            (a \Phi^\mathsf{BR} \pm a \widetilde{\Phi}^\mathsf{BR}) (\spt + ph/2), & p = \pm (-1, 1, 0), \\
            (a \Phi^\mathsf{FU} \pm a \widetilde{\Phi}^\mathsf{FU}) (\spt + ph/2), & p = \pm (1, 0, 1), \\
            (a \Phi^\mathsf{BU} \pm a \widetilde{\Phi}^\mathsf{BU}) (\spt + ph/2), & p = \pm (-1, 0, 1), \\
            (a \Phi^\mathsf{RU} \pm a \widetilde{\Phi}^\mathsf{RU}) (\spt + ph/2), & p = \pm (0, 1, 1), \\
            (a \Phi^\mathsf{LU} \pm a \widetilde{\Phi}^\mathsf{LU}) (\spt + ph/2), & p = \pm (0, -1, 1), \\
            (a \Phi^\mathsf{FRU} \pm a \widetilde{\Phi}^\mathsf{FRU}) (\spt + ph/2), & p = \pm (1, 1, 1), \\
            (a \Phi^\mathsf{FRD} \pm a \widetilde{\Phi}^\mathsf{FRD}) (\spt + ph/2), & p = \pm (1, 1, -1), \\
            (a \Phi^\mathsf{FLU} \pm a \widetilde{\Phi}^\mathsf{FLU}) (\spt + ph/2), & p = \pm (1, -1, 1), \\
            (a \Phi^\mathsf{BRU} \pm a \widetilde{\Phi}^\mathsf{BRU}) (\spt + ph/2), & p = \pm (-1, 1, 1), 
        \end{cases}
    \end{equation}
    where superscripts $``\mathsf{F}"$ (front), $``\mathsf{B}"$ (back), $``\mathsf{R}"$ (right),
    $``\mathsf{L}"$ (left), $``\mathsf{U}"$ (up), $``\mathsf{D}"$ (down), or their combinations indicate the stencil's location. Below, we list the values that these $\Phi$ and $\tilde{\Phi}$ with various superscripts as well as the right-hand basis grid functions take depending on $m$ and $k$:
    \begin{align*}
    & (\Phi^\mathsf{F},\Phi^\mathsf{R},\Phi^\mathsf{U},\Phi^\mathsf{FR},\Phi^\mathsf{BR},\Phi^\mathsf{FU},\Phi^\mathsf{BU},\Phi^\mathsf{RU},\Phi^\mathsf{LU},\Phi^\mathsf{FRU},\Phi^\mathsf{FRD},\Phi^\mathsf{FLU},\Phi^\mathsf{BRU},\\
    & \quad \tilde{\Phi}^\mathsf{F},\tilde{\Phi}^\mathsf{R},\tilde{\Phi}^\mathsf{U},\tilde{\Phi}^\mathsf{FR},\tilde{\Phi}^\mathsf{BR},\tilde{\Phi}^\mathsf{FU},\tilde{\Phi}^\mathsf{BU},\tilde{\Phi}^\mathsf{RU},\tilde{\Phi}^\mathsf{LU},\tilde{\Phi}^\mathsf{FRU},\tilde{\Phi}^\mathsf{FRD},\tilde{\Phi}^\mathsf{FLU},\tilde{\Phi}^\mathsf{BRU}, f^{[m],k}_h) = \\
    & \qquad
    \left[
        \begin{aligned}
            & -\vec{e}_{k}, && \text{if} \; m=-1,0,\; \text{and} \; k=1,\ldots,13,\\
            & -\vec{e}_k, && \text{if} \; m=-1,\; \text{and} \; k=14,\ldots,26,\\
            & \vec{e}_{14} + \vec{e}_{15} - \vec{e}_{17}, && \text{if} \; m=0, k=14 \mbox{; or, } m=1, k=9,\\
            & -\vec{e}_{14} + \vec{e}_{15} - \vec{e}_{18}, && \text{if} \; m = 0, k = 15 \mbox{; or, } m = 1, k = 10,\\
            & \vec{e}_{14} + \vec{e}_{16} - \vec{e}_{19}, && \text{if} \; m = 0, k = 16 \mbox{; or, } m = 1, k = 11,\\
            & -\vec{e}_{14} + \vec{e}_{16} - \vec{e}_{20}, && \text{if} \; m = 0, k = 17 \mbox{; or, } m = 1, k = 12,\\
            & \vec{e}_{15} + \vec{e}_{16} - \vec{e}_{21}, && \text{if} \; m = 0, k = 18 \mbox{; or, } m = 1, k = 13,\\
            & -\vec{e}_{15} + \vec{e}_{16} - \vec{e}_{22}, && \text{if} \; m = 0, k = 19 \mbox{; or, } m = 1, k = 14,\\
            & \vec{e}_{14} + \vec{e}_{15} + \vec{e}_{16} - \vec{e}_{23}, && \text{if} \; m = 0, k = 20 \mbox{; or, } m = 1, k = 15,\\
            & \vec{e}_{14} + \vec{e}_{15} - \vec{e}_{16} - \vec{e}_{24}, && \text{if} \; m = 0, k = 21 \mbox{; or, } m = 1, k = 16,\\
            & \vec{e}_{14} - \vec{e}_{15} + \vec{e}_{16} - \vec{e}_{25}, && \text{if} \; m = 0, k = 22 \mbox{; or, } m = 1, k = 17,\\
            & -\vec{e}_{14} + \vec{e}_{15} + \vec{e}_{16} - \vec{e}_{26}, && \text{if} \; m = 0, k = 23 \mbox{; or, } m = 1, k = 18,\\
            & -\vec{e}_1 - \vec{e}_2 - \vec{e}_3 + h^2 f \vec{e}_{27}, && \text{if} \; m=1,2, \mbox{ and } k=1,\\
            & -2\vec{e}_{k - 1} - \vec{e}_{2k} - \vec{e}_{2k + 1} + 2 h^2 f \vec{e}_{27}, && \text{if} \; m=1,2, \mbox{ and } k = 2, 3, 4,\\
            & -\vec{e}_5 - \vec{e}_7 - \vec{e}_9 - \vec{e}_{10} + 3 h^2 f \vec{e}_{27}, && \text{if} \; m=1,2, \mbox{ and } k=5,\\
            & -\vec{e}_5 - \vec{e}_6 - \vec{e}_8 - \vec{e}_{11} + 3 h^2 f \vec{e}_{27}, && \text{if} \; m=1,2, \mbox{ and } k=6,\\
            & -\vec{e}_4 - \vec{e}_7 - \vec{e}_8 - \vec{e}_{12} + 3 h^2 f \vec{e}_{27}, && \text{if} \; m=1,2, \mbox{ and } k=7,\\
            & -\vec{e}_4 - \vec{e}_6 - \vec{e}_9 - \vec{e}_{13} + 3 h^2 f \vec{e}_{27}, && \text{if} \; m=1,2, \mbox{ and } k=8,\\
            & \begin{aligned}
            & -4\vec{e}_{14}+2(\vec{e}_{17}-\vec{e}_{18}+\vec{e}_{19}-\vec{e}_{20})\\
            & \quad -\vec{e}_{23}-\vec{e}_{24}-\vec{e}_{25}+\vec{e}_{26},
            \end{aligned} && \text{if} \; m=2,k=9 \mbox{; or, } m=3,k=3,\\
            & \begin{aligned}
            & -4\vec{e}_{15}+2(\vec{e}_{17}+\vec{e}_{18}+\vec{e}_{21}-\vec{e}_{22})\\
            & \quad -\vec{e}_{23}-\vec{e}_{24}+\vec{e}_{25}-\vec{e}_{26},
            \end{aligned} && \text{if} \; m=2,k=10 \mbox{; or, } m=3,k=4,\\
            & \begin{aligned}
            & -4\vec{e}_{16}+2(\vec{e}_{19}+\vec{e}_{20}+\vec{e}_{21}+\vec{e}_{22})\\
            & \quad -\vec{e}_{23}+\vec{e}_{24}-\vec{e}_{25}-\vec{e}_{26},
            \end{aligned} && \text{if} \; m=2,k=11 \mbox{; or, } m=3,k=5,
        \end{aligned}
        \right.
    \end{align*}
    where $\vec{e}_k$, $k=1,\ldots,27$ are the standard unit vector basis of $\mathbb{R}^{27}$.

    If $m=3,4$ and $k=1$, then
    \begin{align*}
        \Phi^\mathsf{F} & = -2 + \tfrac{1}{4} h^2 \left( |\nabla \tilde{a}|^2 + \tilde{a}_{xx} - \tilde{a}_{yy} - \tilde{a}_{zz} \right),\\
        \Phi^\mathsf{R} & = -2 + \tfrac{1}{4} h^2 \left( |\nabla \tilde{a}|^2 - \tilde{a}_{xx} + \tilde{a}_{yy} - \tilde{a}_{zz} \right),\\
        \Phi^\mathsf{U} & = -2 + \tfrac{1}{4} h^2 \left( |\nabla \tilde{a}|^2 - \tilde{a}_{xx} - \tilde{a}_{yy} + \tilde{a}_{zz} \right),\\
        \Phi^\mathsf{FR} & = -1 + \tfrac{1}{4} h^2 \tilde{a}_{xy},
        \quad \Phi^\mathsf{FU} = -1 + \tfrac{1}{4} h^2 \tilde{a}_{xz},
        \quad \Phi^\mathsf{RU} = -1 + \tfrac{1}{4} h^2 \tilde{a}_{yz},\\
        \Phi^\mathsf{BR} & = -1 - \tfrac{1}{4} h^2 \tilde{a}_{xy},
        \quad \Phi^\mathsf{BU} = -1 - \tfrac{1}{4} h^2 \tilde{a}_{xz},
        \quad \Phi^\mathsf{LU}  = -1 - \tfrac{1}{4} h^2 \tilde{a}_{yz},\\
        f^{[m],k}_h & = -a \left( 6 h^2 \tilde{f} + \tfrac{1}{2} h^4 \left( \Delta \tilde{f} - \nabla \tilde{a} \cdot \nabla \tilde{f} \right) \right),
    \end{align*}
    and $\Phi^\mathsf{FRU}$, $\Phi^\mathsf{FRD}$, $\Phi^\mathsf{FLU}$, $\Phi^\mathsf{BRU}$, $\tilde{\Phi}^\mathsf{F}$, $\tilde{\Phi}^\mathsf{R}$, $\tilde{\Phi}^\mathsf{U}$, $\tilde{\Phi}^\mathsf{FR}$, $\tilde{\Phi}^\mathsf{BR}$, $\tilde{\Phi}^\mathsf{FU}$, $\tilde{\Phi}^\mathsf{BU}$, $\tilde{\Phi}^\mathsf{RU}$, $\tilde{\Phi}^\mathsf{LU}$, $\tilde{\Phi}^\mathsf{FRU}$, $\tilde{\Phi}^\mathsf{FRD}$, $\tilde{\Phi}^\mathsf{FLU}$, $\tilde{\Phi}^\mathsf{BRU}$ are all zeros.

    Finally, if $m=3,4$ and $k=2$, then
    \begin{align*}
        & \makebox[120pt][l]{$\Phi^\mathsf{F} = -8 + h^2 \tilde{a}_{xx},$}
        \makebox[120pt][l]{$\Phi^\mathsf{R} = -8 + h^2 \tilde{a}_{yy},$}
        \makebox[110pt][l]{$\Phi^\mathsf{U} = -8 + h^2 \tilde{a}_{zz},$} \\
        & \makebox[120pt][l]{$\widetilde{\Phi}^\mathsf{F} = \frac{1}{4} h^3 \partial_x (|\nabla \tilde{a}|^2 - \Delta \tilde{a}),$}
        \makebox[120pt][l]{$\widetilde{\Phi}^\mathsf{R} = \frac{1}{4} h^3 \partial_y (|\nabla \tilde{a}|^2 - \Delta \tilde{a}),$}
        \makebox[110pt][l]{$\widetilde{\Phi}^\mathsf{U} = \frac{1}{4} h^3 \partial_z (|\nabla \tilde{a}|^2 - \Delta \tilde{a}),$} \\
        & \makebox[185pt][l]{$\Phi^\mathsf{FRU} = -1 + \frac{1}{4} h^2 \big( \tilde{a}_{xy} + \tilde{a}_{xz} + \tilde{a}_{yz} \big),$}
        \makebox[185pt][l]{$\Phi^\mathsf{FRD} = -1 + \frac{1}{4} h^2 \big( \tilde{a}_{xy} - \tilde{a}_{xz} - \tilde{a}_{yz} \big),$} \\
        & \makebox[185pt][l]{$\Phi^\mathsf{FLU} = -1 + \frac{1}{4} h^2 \big( -\tilde{a}_{xy} + \tilde{a}_{xz} - \tilde{a}_{yz} \big),$}
        \makebox[185pt][l]{$\Phi^\mathsf{BRU} = -1 + \frac{1}{4} h^2 \big( -\tilde{a}_{xy} - \tilde{a}_{xz} + \tilde{a}_{yz} \big),$} \\
        & f_h^{[m], k} = -a \left( 12 h^2 \tilde{f} + h^4 \left( \tfrac{1}{2} \big( \Delta \tilde{f} - \nabla \tilde{a} \cdot \nabla \tilde{f} \big) + \big( |\nabla \tilde{a}|^2 - \Delta \tilde{a} \big) \tilde{f} \right) \right),
    \end{align*}
    and $\Phi^\mathsf{FR}$, $\Phi^\mathsf{BR}$, $\Phi^\mathsf{FU}$, $\Phi^\mathsf{BU}$, $\Phi^\mathsf{RU}$, $\Phi^\mathsf{LU}$, $\tilde{\Phi}^\mathsf{FR}$, $\tilde{\Phi}^\mathsf{BR}$, $\tilde{\Phi}^\mathsf{FU}$, $\tilde{\Phi}^\mathsf{BU}$, $\tilde{\Phi}^\mathsf{RU}$, $\tilde{\Phi}^\mathsf{LU}$, $\tilde{\Phi}^\mathsf{FRU}$, $\tilde{\Phi}^\mathsf{FRD}$, $\tilde{\Phi}^\mathsf{FLU}$, $\tilde{\Phi}^\mathsf{BRU}$ are all zeros.

    Same as before, the basis stencil coefficients $C_p^{[m],k} (\spt)$ are symmetric if all $\tilde{\Phi}$ are zeros. Meanwhile, if all $\Phi$ are zeros, then the basis stencil coefficients $C_p^{[m],k} (\spt)$ are anti-symmetric. Additionally, the basis stencil coefficient $C_p^{[4], 2}$ is neither symmetric nor anti-symmetric, but it can be decomposed into symmetric and anti-symmetric parts by taking even and odd powers of $h$. The freedom in the choices of $\kappa_{m, k}$ in \eqref{eq:C_p_3D} can be exploited to help us achieve \eqref{cpmk:sym}, minimize the leading truncation error, and prove the SPD property as well as the convergence rates. For a general $a$, we can take $M \le 4$ in \eqref{eq:C_p_3D} to ensure that the resulting finite difference scheme is symmetric. If symmetry is not required, or symmetry is required but $2 \Delta \tilde{a} - |\nabla \tilde{a}|^2$ is constant, we can take $M \le 6$ in \eqref{eq:C_p_3D}. For $M=5,6$, there is a single basis stencil coefficient, but it is difficult to compute symbolically.

    \subsection{The general $d$-dimensional case}
    \label{sec:any_dimension}


    We begin by adapting the lower dimensional schemes into higher dimensions $d \ge 3$, where they become schemes that involve partial derivatives along one and two axes. To this end, we define $\mathcal{L}_i' := \partial^{2 \vec{e}_i} - \partial^{\vec{e}_i} \tilde{a} \cdot \partial^{\vec{e}_i}$ and $\tilde{f}_i = \mathcal{L}_i' u$, $1 \le i \le d$. Observe that the original PDE in \eqref{eq:PDE} can be written as $\sum_{i = 1}^d \mathcal{L}_i' u = \tilde{f} = \sum_{i = 1}^d \tilde{f}_i$. Let $\SS = \{-1, 0, 1\}^d$ as usual. Then, the scheme along the $i$-th direction is given by $\sum_{p \in \SS} C_{p, i} (\spt) u_h(\spt + ph) = f_{h, i} (\spt)$, where
    \begin{align*}
        C_{p, i} (\spt) & =
        \begin{cases}
            a \left( -1 + \frac{1}{24} h^2 \left( \partial^{2 \vec{e}_i} \tilde{a} + (\partial^{\vec{e}_i} \tilde{a})^2 \right) \right) \big|_{\spt + ph/2}, & \text{if} \quad p = \pm \vec{e}_i,\\
            -(C_{-\vec{e}_i, i} (\spt) + C_{\vec{e}_i, i} (\spt)), & \text{if} \quad p = \vec{0}, \\
            0, & \text{otherwise},
        \end{cases} \\
        f_{h, i} & = -a \big( h^2 \tilde{f}_i + \tfrac{1}{12} h^4 \mathcal{L}_i' \tilde{f}_i \big)
        = -a \left( h^2 \tilde{f}_i + \tfrac{1}{12} h^4 (\partial^{2 \vec{e}_i} \tilde{f}_i - \partial^{\vec{e}_i} \tilde{a} \partial^{\vec{e}_i} \tilde{f}_i) \right).
    \end{align*}
    The scheme along the $i$-th and $j$-th direction is given by $\frac{1}{6}$ times the 2D basis coefficient $C_{p}^{[4], 1}$ and the corresponding right-hand basis grid function $f_h^{[4], 1}$ in \Cref{sec:basis}. More explicitly, the scheme along two directions $1 \le i < j \le d$ is given by $\sum_{p \in \SS} C_{p, i, j} (\spt) u_h(\spt + ph) = f_{h, i, j} (\spt)$, where
    \begin{align*}
        & C_{p,i,j}(\spt) = \\
        & \quad
        \begin{cases}
            a(-\frac{2}{3} + \frac{1}{24} h^2 \left( (\partial^{\vec{e}_i} \tilde{a})^2 + (\partial^{\vec{e}_j}\tilde{a})^2 + \partial^{2\vec{e}_i} \tilde{a} - \partial^{2\vec{e}_j}\tilde{a} \right))\big|_{\spt + ph/2}, & \text{if} \quad p = \pm \vec{e}_i,\\
            a(-\frac{2}{3} + \frac{1}{24} h^2 \left( (\partial^{\vec{e}_i} \tilde{a})^2 + (\partial^{\vec{e}_j}\tilde{a})^2 - \partial^{2\vec{e}_i} \tilde{a} + \partial^{2\vec{e}_j}\tilde{a} \right))\big|_{\spt + ph/2}, & \text{if} \quad p = \pm \vec{e}_j,\\
            a(-\frac{1}{6} + \frac{1}{24} h^2 \partial^{\vec{e}_i + \vec{e}_j}\tilde{a})\big|_{\spt + ph/2}, & \text{if} \quad p = \pm (\vec{e}_i+\vec{e}_j), \\
            a(-\frac{1}{6} - \frac{1}{24} h^2 \partial^{\vec{e}_i + \vec{e}_j}\tilde{a})\big|_{\spt + ph/2}, & \text{if} \quad p = \pm (\vec{e}_i-\vec{e}_j), \\
            -\sum_{p \in \SS \bs \{\vec{0}\}} C_{p,i,j} (\spt), & \text{if} \quad p = \vec{0},\\
            0, & \text{otherwise},
        \end{cases}\\
        & f_{h,i,j} = -a \left( h^2 (\tilde{f}_i + \tilde{f}_j) + \tfrac{1}{12} h^4 (\mathcal{L}_i' + \mathcal{L}_j') (\tilde{f}_i + \tilde{f}_j) \right).
    \end{align*}
    Note that these two schemes are 4th-order consistent with the PDEs $\mathcal{L}_i' u = \tilde{f}_i$ and $(\mathcal{L}_i' + \mathcal{L}_j') u = \tilde{f}_i + \tilde{f}_j$ respectively.

    Let us turn to the $d$-dimensional scheme. The observation is that, setting
    \begin{equation}
        \label{eq:f_h_any_dim}
        f_h := -a \left( h^2 \tilde{f} + \tfrac{1}{12} h^4 (\Delta \tilde{f} - \nabla \tilde{a} \cdot \nabla \tilde{f}) \right),
    \end{equation}
    we have
    \begin{align*}
        f_h &= -a h^2 \sum_{i = 1}^d \tilde{f}_i - \tfrac{1}{12} a h^4 \sum_{i = 1}^d \sum_{j = 1}^d \mathcal{L}_i' \tilde{f}_j \\
        &= -a \Bigg[ (d - 1) h^2 \sum_{i = 1}^d \tilde{f}_i + \tfrac{1}{12} h^4 \bigg( (d - 1) \sum_{i = 1}^d \mathcal{L}_i' \tilde{f}_i + \sum_{1 \le i < j \le d} (\mathcal{L}_i' \tilde{f}_j + \mathcal{L}_j' \tilde{f}_i) \bigg) \Bigg] \\
        & \quad \ + (d - 2) a \Bigg[ h^2 \sum_{i = 1}^d \tilde{f}_i + \tfrac{1}{12} h^4 \sum_{i = 1}^d \mathcal{L}_i' \tilde{f}_i \Bigg] \\
        &= \sum_{1 \le i < j \le d} f_{h, i, j} - (d - 2) \sum_{i = 1}^d f_{h, i}.
    \end{align*}
    Now, for each $p \in \SS$, if we set
    \begin{equation}
        \label{eq:C_p_any_dim}
        C_p(\spt) = \sum_{1 \le i < j \le d} C_{p, i, j} (\spt) - (d - 2) \sum_{i = 1}^d C_{p, i} (\spt),
    \end{equation}
    then the finite difference scheme for \eqref{eq:PDE} given by $\sum_{p \in \SS} C_p (\spt) u_h(\spt + ph) = f_h (\spt)$, where \eqref{eq:f_h_any_dim} and \eqref{eq:C_p_any_dim} hold, is compact and 4th-order consistent. Moreover, this finite difference scheme is symmetric, since it is a linear combination of a scheme involving derivatives along a single axis and another scheme involving derivatives along two axes, both of which are inherently symmetric. Note that $C_p(\spt) = 0$ except for $p \in \SS$ with up to two nonzero coordinates.

    When $d = 3$, the scheme in this section coincides with the basis coefficients $C_p^{[4], 1}$ and the corresponding right-hand basis grid function. When $a$ is constant, this scheme coincides with the $4$th-order ones proposed by \cite{kyei2010family,spotz1996high,zhai2013family}.

    As a final remark for this section, we note that similar calculation and analysis can be done for the non-uniform case, where the grid size in each axis is different. However, our calculation suggests that the maximum consistency order for a compact symmetric $d$-dimensional FDM with $d \ge 2$ on a non-uniform grid is only 2. Therefore, we omit this case from the present paper.

\section{Derivations of compact, symmetric, and high-order FDMs}
    \label{sec:details}
    In this section, we present the complete derivations and proofs of the main results in \cref{sec:main:thms,subsec:SPD:conv:thms}, together with the supporting details for \cref{sec:construct:dD}.

    \subsection{The 1D case}
    \label{sec:pf:1D}
    The proof for $d=1$ is easier than the $d\ge 2$ case, so we treat it separately.
    \begin{proof}[Proof of \cref{thm:1D}]
        We recall a few basic facts. By \cite[Proposition 3.1]{HanMichelleWong}, for any $b^{*} \in [0,1]$,
        \be \label{expandu}
        \begin{aligned}
        u(b^* + h) & = u(b^*) + u'(b^*) h \left( 1 + \sum_{j=2}^{M+1} \tfrac{E_{j,1}(b^*)}{j!} h^{j-1} \right) \\
        & \quad + \sum_{\ell=0}^{M-1} h^{\ell+2} f^{(\ell)} (b^*) \left( \sum_{j=\ell+2}^{M+1} \tfrac{F_{j,\ell}(b^*)}{j!}h^{j-\ell-2} \right) + \bo_{a,f}(h^{M+2}),
        \end{aligned}
        \ee
        where the quantities $E_{j,1}$, $F_{j,k}$ are computed using the above recursions. Furthermore, if we define $v:=a u'$, take the Taylor expansion of $v(b^* + h)$ about the point $b^*$, and use the fact that $-(a u')' =f$, we have
        \be \label{auprime}
        \begin{aligned}
        a(b^* + h) u'(b^* + h) & = a(b^*) u'(b^*) + \sum_{j=1}^{M} \tfrac{(a u')^{(j)} (b^*)}{j!} h^j + \bo_{a,f}(h^{M+1})\\
        & = a(b^*) u'(b^*) - \sum_{j=1}^{M} \tfrac{f^{(j-1)} (b^*)}{j!} h^j + \bo_{a,f}(h^{M+1}).
        \end{aligned}
        \ee

        Now, we expand $u(c^*) - u(c^*-h)$ about the point $c^* - h/2$ by using \eqref{expandu}. More specifically, to expand $u(c^*)$, we let $b^{*} = c^* - h/2$ and replace $h$ with $h/2$ in \eqref{expandu}. On the other hand, to expand $u(c^* - h)$, we let $b^{*} = c^* - h/2$ and replace $h$ with $-h/2$ in \eqref{expandu}. We do a similar expansion for $u(c^*+h) - u(c^*)$ about the point $c^* + h/2$ to get the following results
        \begin{align*}
        u(c^*) - u(c^* -h) & = u'(c^* - \tfrac{h}{2})
        \left( 2 + \sum_{j=2}^{M+1} \tfrac{E_{j,1}\left(c^* - h/2\right)}{j!} (\tfrac{h}{2})^{j-1} (1+(-1)^{j-1})\right)
        \tfrac{h}{2} \\
        & \quad + \sum_{\ell = 0}^{M-1} f^{(\ell)}(c^* - \tfrac{h}{2}) \mathcal{G}_{\ell}(c^* - \tfrac{h}{2}) (\tfrac{h}{2})^{\ell+2} + \bo_{a,f}(h^{M+2}),\\
        u(c^*+h) - u(c^*) &  = u'(c^* + \tfrac{h}{2}) \left( 2 + \sum_{j=2}^{M+1} \tfrac{E_{j,1}\left(c^* + h/2\right)}{j!} (\tfrac{h}{2})^{j-1} (1+(-1)^{j-1})\right) \tfrac{h}{2} \\
        & \quad + \sum_{\ell = 0}^{M-1} f^{(\ell)}(c^* + \tfrac{h}{2}) \mathcal{G}_{\ell}(c^{*} + \tfrac{h}{2})(\tfrac{h}{2})^{\ell+2} + \bo_{a,f}(h^{M+2}),
        \end{align*}
        where $\mathcal{G}_{\ell}$ is defined in \eqref{Gk}. Using $C_{-1}, C_0, C_1$ in \eqref{c:1D}, we have
        \begin{align*}
            & C_{-1}(c^*) u(c^{*}-h) + C_0(c^*) u(c^{*}) + C_1(c^*) u(c^{*}+h) \\
            & \quad = - C_{-1}(c^*) (u(c^{*})-u(c^{*}-h)) + C_{1}(c^*) (u(c^* +h) - u(c^{*})) \\
            & \quad = a(c^* - \tfrac{h}{2}) u'(c^* - \tfrac{h}{2}) h - C_{-1}(c^*) \sum_{\ell = 0}^{M-1} f^{(\ell)}(c^* - \tfrac{h}{2}) \mathcal{G}_{\ell}(c^{*}-\tfrac{h}{2})(\tfrac{h}{2})^{\ell+2} \\
            & \qquad - a(c^{*} + \tfrac{h}{2}) u'(c^{*} + \tfrac{h}{2}) h + C_{1}(c^*) \sum_{\ell=0}^{M-1} f^{(\ell)}(c^* + \tfrac{h}{2}) \mathcal{G}_{\ell}(c^{*}+\tfrac{h}{2})(\tfrac{h}{2})^{\ell+2} \\
            & \qquad + \bo_{a,f}(h^{M+2}) \\
            & \quad = -2 \sum_{j=1}^{M} \tfrac{f^{(j-1)}(c^{*})}{j!} ((-1)^{j}-1) (\tfrac{h}{2})^{j+1} - C_{-1}(c^*) \sum_{\ell = 0}^{M-1} f^{(\ell)}(c^* - \tfrac{h}{2}) \mathcal{G}_{\ell}(c^{*}-\tfrac{h}{2}) (\tfrac{h}{2})^{\ell+2} \\
            & \qquad + C_{1}(c^*) \sum_{\ell=0}^{M-1} f^{(\ell)}(c^* + \tfrac{h}{2}) \mathcal{G}_{\ell}(c^{*}+\tfrac{h}{2}) (\tfrac{h}{2})^{\ell+2}  + \bo_{a,f}(h^{M+2}) \\
            & \quad = \sum_{\ell = 0}^{M-1} \sum_{j=0}^{M-1-\ell} \tfrac{f^{(\ell+j)}(c^*)}{j!} \left( - C_{-1}(c^*)  (-1)^j \mathcal{G}_{\ell}(c^{*}-\tfrac{h}{2}) + C_{1}(c^*) \mathcal{G}_{\ell}(c^{*}+\tfrac{h}{2}) \right) (\tfrac{h}{2})^{\ell+j+2} \\
            & \qquad - 2\sum_{j=0}^{M-1} \tfrac{f^{(j)}(c^{*})}{(j+1)!} ((-1)^{j+1}-1) (\tfrac{h}{2})^{j+2} + \bo_{a,f}(h^{M+2}) \\
            & \quad = \sum_{\ell = 0}^{M-1} \tfrac{f^{(\ell)}(c^*)}{2^{\ell + 2}} \left(- \tfrac{2((-1)^{\ell+1}-1)}{(\ell+1)!} \right.\\
            & \qquad \left. + \sum_{k = 0}^{\ell} \tfrac{C_{-1}(c^*) (-1)^{\ell - k+1} \mathcal{G}_{k}(c^{*}-\tfrac{h}{2}) + C_{1}(c^*) \mathcal{G}_{k}(c^{*}+\tfrac{h}{2})}{(\ell-k)!}\right) h^{\ell+2} + \bo_{a,f}(h^{M+2}),
        \end{align*}
        where we expanded $f^{(\ell)}(c^{*}-\tfrac{h}{2})$ and $f^{(\ell)}(c^{*}+\tfrac{h}{2})$ about the point $c^{*}$ to arrive at the second last line, and rearrange the indices to obtain the final line. The above calculation implies that
        \begin{align*}
        C_{-1}(c^{*}) u(c^{*}-h) + C_0(c^{*}) u(c^{*}) + C_1(c^{*}) u(c^{*}+h) & = \sum_{\ell=0}^{M-1} d_{\ell}(h) h^{\ell+2} f^{(\ell)}(c^*) \\
        & \quad + \bo_{a,f}(h^{M+2}),
        \end{align*}
        where $d_{\ell}(h)$ is defined in $\eqref{c:1D}$. Setting $f_h(c^*) = \sum_{\ell=0}^{M-1} d_{\ell}(h) h^{\ell+2} f^{(\ell)}(c^*)$, we have $f_h(c^*) = d_0(0) h^2 f(c^*) + \bo_{a, f} (h^3) = f(c^*)h^2 + \frac{1}{4} (\mathcal{E}^{M}(c^*) - \mathcal{E}^{M}(c^*)) \mathcal{G}_0 (c^*) f(c^*)h^2 + \bo_{a, f} (h^3) = f(c^*)h^2 + \bo_{a, f} (h^3)$. Therefore, the finite difference scheme $\mathcal{L}_h u_h = f_h$ is $M$th-order consistent in the sense of \eqref{eq:consistency_def} with $\Lambda(h) = h^2$. The proof is completed.
    \end{proof}

    In the remaining part of this section, we adopt the approach of \cite{feng2024sixth,han2025convergent} to present a general framework for the proposed $d$-dimensional finite difference schemes with $d \ge 2$. We present the structure of all possible stencil coefficients and prove that the maximum consistency order is generally 4 for compact symmetric finite difference schemes in higher dimensions.

    \subsection{Expansion of the solution for $d \ge 2$}
    \label{sec:nD_expansion}
    This subsection contains key relations used in the upcoming lemmas, propositions, and ultimately our main result, \Cref{thm:max_order}. We first rewrite the PDE in \eqref{eq:PDE} as
    \begin{equation}
        \label{eq:PDE_simp}
        \partial^{2 \vec{e}_1} u = \nabla \tilde{a} \cdot \nabla u - \sum_{j = 2}^d \partial^{2 \vec{e}_j} u + \tilde{f}
        \quad \mbox{ with }\quad \tilde{a} := -\ln a \quad \mbox{ and }\quad \tilde f := -\frac{f}{a}.
    \end{equation}
    Taking partial derivatives of both sides of \cref{eq:PDE_simp}, applying the Leibniz rule with $\bm{k} \ge 2 \vec{e}_1$ to the first term, and rearranging the indices, we get
    \begin{equation}
        \label{eq:PDE_deriv}
        \begin{aligned}
            \partial^{\bm{k}} u
            &= \sum_{\substack{\vec{0} \le \bm{\ell} \le \bm{k} - 2\vec{e}_1}} \binom{\bm{k} - 2\vec{e}_1}{\bm{\ell}} \sum_{j = 1}^d \partial^{\bm{k} - \bm{\ell} - 2 \vec{e}_1 + \vec{e}_j} \tilde{a} \, \partial^{\bm{\ell} + \vec{e}_j} u  - \sum_{j = 2}^d \partial^{\bm{k} - 2 \vec{e}_1 + 2 \vec{e}_j} u + \partial^{\bm{k} - 2 \vec{e}_1} \tilde{f} \\
            &= \sum_{\substack{\vec{0} \le \bm{\ell} \le \bm{k} \\ \ell_1 < k_1}} \tilde{a}^{\bm{k}}_{\bm{\ell}} \partial^{\bm{\ell}} u + \partial^{\bm{k} - 2 \vec{e}_1} \tilde{f},
        \end{aligned}
    \end{equation}
    where
    \begin{equation}
        \label{eq:a_tilde_kl}
        \tilde{a}^{\bm{k}}_{\bm{\ell}} =
        \begin{cases}
            \displaystyle
            - \sum_{j = 2}^d \td(|\bm{k} - \bm{\ell} - 2 \vec{e}_1 + 2 \vec{e}_j|), & \mbox{if} \ |\bm{\ell}| = |\bm{k}|, \\
            \displaystyle
            \sum_{j = 1}^d (1 - \td(\ell_j)) \binom{\bm{k} - 2\vec{e}_1}{\bm{\ell} - \vec{e}_j}  \partial^{\bm{k} - \bm{\ell} - 2 \vec{e}_1 + 2 \vec{e}_j} \tilde{a}, & \mbox{if} \ |\bm{\ell}| < |\bm{k}|.
        \end{cases}
    \end{equation}
    Now, if we use \eqref{eq:PDE_deriv} repeatedly, we can represent $\partial^{\bm{k}} u$ in terms of $\partial^{\bm{\ell}} u$ for $\bm{\ell} \in \LN_{|\bm{k}|}$, where $\LN_M$ is defined in \eqref{PiM}. More concisely, we have
    \begin{equation}
        \label{eq:u_k_repr}
        \partial^{\bm{k}} u
        = \sum_{\bm{\ell} \in \LN_{|\bm{k}|}} \tilde{A}^{\bm{k}}_{\bm{\ell}}
        \partial^{\bm{\ell}} u + \tilde{F}_{\bm{k}},
    \end{equation}
    where the coefficients $\tilde{A}^{\bm{k}}_{\bm{\ell}}$ and $\tilde{F}_{\bm{k}}$ are uniquely determined through the following recursive formulas:
    \begin{equation}
        \label{eq:AF_tilde}
        \tilde{A}^{\bm{k}}_{\bm{\ell}} :=
        \begin{cases}
            \td(|\bm{k} - \bm{\ell}|), & \mbox{if} \ k_1 \le 1, \\[3pt]
            \sum\limits_{\substack{|\bm{\ell}| \le |\bm{j}| \le |\bm{k}| \\ j_1 < k_1}} \tilde{a}^{\bm{k}}_{\bm{j}} \tilde{A}^{\bm{j}}_{\bm{\ell}}, & \mbox{if} \ k_1 \ge 2,
        \end{cases}
        \quad
        \tilde{F}_{\bm{k}} :=
        \begin{cases}
            0, & \mbox{if} \ k_1 \le 1, \\[3pt]
            \sum\limits_{\substack{|\bm{\ell}| \le |\bm{k}| \\ \ell_1 < k_1}} \tilde{a}^{\bm{k}}_{\bm{\ell}} \tilde{F}_{\bm{\ell}} + \partial^{\bm{k} - 2 \vec{e}_1} \tilde{f}, & \mbox{if} \ k_1 \ge 2.
        \end{cases}
    \end{equation}
    For $M \in \NN$, a base point $\bpt \in \overline{\Omega}$ and $p = (p_1, \ldots, p_d) \in \Z^d$, we have the Taylor expansion
    \begin{equation*}
        u(\bpt + ph) = \sum_{0 \le |\bm{k}| \le M + 1} \left( \prod_{j = 1}^d \frac{p_j^{k_j}}{k_j!} \right) h^{|\bm{k}|} \partial^{\bm{k}} u(\bpt) + \bo_u (h^{M + 2}).
    \end{equation*}
    Now we incorporate \eqref{eq:u_k_repr} into the Taylor expansion of $u$ to obtain \eqref{eq:u_taylor_int}, where $A^k_{\bm{\ell}} (p)$ and $F(p)$ are explicitly defined below:
    \begin{equation}
        \label{eq:A^k_l(p)}
        A^k_{\bm{\ell}} (p) := \sum_{|\bm{k}| = k} \left( \prod_{j = 1}^d \frac{p_j^{k_j}}{k_j!} \right) \tilde{A}^{\bm{k}}_{\bm{\ell}}
        \quad \mbox{for} \quad \bm{\ell} \in \LN_{M + 1} \quad \mbox{and} \quad k = 0, \ldots, M + 1,
    \end{equation}
    and
    \begin{equation}
        \label{eq:F(p)}
        F(p) := \sum_{|\bm{k}| \le M + 1} \left( \prod_{j = 1}^d \frac{p_j^{k_j}}{k_j!} \right) \tilde{F}_{\bm{k}} h^{|\bm{k}|}.
    \end{equation}

    We also consider the converse to guarantee that our FDMs are exhaustive. It will be used in the proof of \Cref{lem:C_p_constraint_int}. For any real numbers $( u_{\bm{\ell}} )_{ \bm{\ell} \in \LN_{M + 1} }$, we construct a function $v$ with
    \begin{equation*}
        v(\bpt + ph) = \sum_{0 \le |\bm{k}| \le M + 1} \left( \prod_{j = 1}^d \frac{p_j^{k_j}}{k_j!} \right) h^{|\bm{k}|} v_{\bm{k}},
        \quad 
        v_{\bm{k}} =
        \begin{cases}
            u_{\bm{k}}, & \mbox{if} \ \bm{k} \in \LN_{M + 1}, \\
            \sum\limits_{\bm{\ell} \in \LN_{|\bm{k}|}} \tilde{A}^{\bm{k}}_{\bm{\ell}} v_{\bm{\ell}} + \tilde{F}_{\bm{k}}, & \mbox{if} \ \bm{k} \not\in \LN_{M + 1},
        \end{cases}
    \end{equation*}
    (compare with \eqref{eq:u_k_repr}), then $\partial^{\bm{\ell}} v(\bpt) = u_{\bm{\ell}}$ for $\bm{\ell} \in \LN_{M + 1}$ and $\partial^{\bm{\ell}} \nabla \cdot (a \nabla v) = \partial^{\bm{\ell}} f$ at $\bpt$ for $|\bm{\ell}| \le M - 1$.

    In preparation for \Cref{sec:structure} (particularly for the proof of \Cref{lem:C_p_constraint_int}), we mention some special values for the quantities defined above. When $k_1 \ge 2$, we have $\tilde{a}^{\bm{k}}_{\vec{0}} = 0$ by \eqref{eq:a_tilde_kl}. Now we can deduce recursively from \eqref{eq:AF_tilde} that $\tilde{A}^{\bm{k}}_{\vec{0}} = \td(|\bm{k}|)$ for all $\bm{k} \in \NN^d$. It follows from \eqref{eq:A^k_l(p)} that $A^k_{\vec{0}}(p) = \td(k)$ for $k \in \NN$. Besides, we can calculate from \eqref{eq:AF_tilde} and \eqref{eq:F(p)} that
    \begin{equation}
        \label{eq:F(p)_leading}
        F(p) = \frac{1}{2} p_1^2 h^2 \tilde{f} (\spt) + \bo_{\tilde{a}, \tilde{f}} (h^3)
        = -\frac{1}{2a(\spt)} p_1^2 h^2 f(\spt) + \bo_{\tilde{a}, \tilde{f}} (h^3),
    \end{equation}
    where $p_1$ is the first entry of $p \in \Z^d$ and we recall that $\tilde{f}=-f/a$ as stated in \eqref{eq:PDE_simp}. The coefficient $\tilde{A}^{\bm{k}}_{\bm{\ell}}$ is of particular interest when $|\bm{k}| = |\bm{\ell}|$ and $\bm{\ell} \in \LN_{|\bm{k}|}$. In this case,  \eqref{eq:a_tilde_kl} and \eqref{eq:AF_tilde} imply
    \begin{equation*}
        \tilde{A}^{\bm{k}}_{\bm{\ell}} =
        \begin{cases}
            \td(|\bm{k} - \bm{\ell}|), & \mbox{if } k_1 \le 1, \\
            -\sum_{j = 2}^d \tilde{A}^{\bm{k} - 2 \vec{e}_1 + 2 \vec{e}_j}_{\bm{\ell}}, & \mbox{if } k_1 \ge 2.
        \end{cases}
    \end{equation*}
    Combined with \eqref{eq:A^k_l(p)}, we can see that $A^k_{\bm{\ell}}(p)$ is a constant independent of the functions $a$, $f$, and $g$ in the PDE \eqref{eq:PDE}. In the particular case of $d = 2$, we can deduce that
    \begin{equation}
        \label{eq:A_tilde_special}
        \tilde{A}^{\bm{k}}_{\bm{\ell}} = 0
        \quad \mbox{if} \quad \frac{1}{2} (\bm{k} - \bm{\ell}) \not \in \NN^d \quad \mbox{and} \quad
        \tilde{A}^{\bm{k}}_{\bm{\ell}} = (-1)^{(k_1 - l_1) / 2}
        \quad \mbox{if} \quad \frac{1}{2} (\bm{k} - \bm{\ell}) \in \NN^d.
    \end{equation}
    This yields the following closed form of $A^k_{\bm{\ell}} (p)$ when $d = 2$:
    \begin{equation}
        \label{eq:A^k_0k}
        A^k_{(0, k)}(p) = \sum_{\substack{j \le k \\ j \text{ is even}}} \frac{p_1^j p_2^{k - j}}{j! (k - j)!} (-1)^{j/2}
        = \frac{1}{k!} \re \left( \ii p_1 + p_2 \right)^k.
    \end{equation}
    Similarly, $A^k_{(1, k - 1)}(p) = \frac{1}{k!} \im \left( \ii p_1 + p_2 \right)^k$ for $k \in \N$.

    \subsection{Structure of the stencil coefficients for $d\ge 2$}
    \label{sec:structure}

    Let $\spt \in \Omega_h$ and $\SS = \{-1,0,1\}^d$. We expand the solution $u$ at each point $\spt + ph$ according to \eqref{eq:u_taylor_int} with the base point $\bpt = \spt$. In view of this, we aim to construct a finite difference scheme in the form of \eqref{eq:FDM_form} that satisfies \eqref{eq:consistency_def} with $f_h = \sum_{p \in \SS} C_p (\spt) F(p)$. The conditions on $C_p$ are given by the following lemma.
    \begin{lemma}
        \label{lem:C_p_constraint_int}
        Let $M \ge -1$ be an integer. Consider the finite difference scheme $\mathcal{L}_h u_h = f_h$ for \eqref{eq:PDE} on $\Omega_h$ satisfying the assumptions of \Cref{thm:max_order}. Then, the following statements are equivalent:
        \begin{itemize}
            \item[(a)] The discretization operator $\mathcal{L}_h$ satisfies the following relation
            \be \label{eq:Lh_int}
                \mathcal{L}_h u (\spt) = \sum_{p \in \SS} C_p (\spt) F(p) + \bo_{\tilde{a}, u} (h^{M + 2}), \quad \forall \spt \in \Omega_h,
            \ee
            where $F(p)$ is defined in \eqref{eq:F(p)}.

            \item[(b)] For all $\spt \in \Omega_h$, we have
            \begin{equation}
                \label{eq:C_p_constraint_int}
                \sum_{p \in \SS} A^{|\bm{\ell}|}_{\bm{\ell}}(p) c_{p, j}
                = -\sum_{k = 0}^{j - 1} \sum_{p \in \SS} A^{|\bm{\ell}| + j - k}_{\bm{\ell}}(p) c_{p, k}, \; \forall \, j = 0, \ldots, M + 1, \ \bm{\ell} \in \LN_{M + 1 - j},
            \end{equation}
            where $A^k_{\bm{\ell}}(p)$ is defined in \eqref{eq:A^k_l(p)}. In particular, $\sum_{p \in \SS} C_p(\spt) = 0$.

            \item[(c)] $\mathcal{L}_h u (\spt)$ is independent of $\partial^{\bm{\ell}} u(\spt)$, $\bm{\ell} \in \LN_{M + 1}$ for any $u \in C^\infty (\overline{\Omega})$ and any $\spt \in \Omega_h$.
        \end{itemize}
        In addition, consider the statement
        \begin{itemize}
            \item[(d)] There exists a grid function $f_h$ such that the scheme $\mathcal{L}_h u_h = f_h$ is $M$th-order consistent.
        \end{itemize}
        Statement (d) implies any of the statements (a)-(c), which further implies statement (d) with the extra assumption of $\sum_{p = (p_1, \ldots, p_d) \in \SS} p_1^2 c_{p, 0} \neq 0$.
    \end{lemma}

    \begin{proof}
        We observe that expanding $u(\spt + ph)$ at the base point $\bpt = \spt$ via \eqref{eq:u_taylor_int} yields
        \begin{equation}
            \label{eq:Lh_int_expansion}
            \mathcal{L}_h u(\spt)
            = \sum_{\bm{\ell} \in \LN_{M + 1}} P_{\bm{\ell}}(h) \partial^{\bm{\ell}} u(\spt) +
            \sum_{p \in \SS} C_p (\spt) F(p) + \bo_{\tilde{a}, u} (h^{M + 2}),
        \end{equation}
        where
        \begin{equation*}
            P_{\bm{\ell}} (h) := \sum_{k = |\bm{\ell}|}^{M + 1} \sum_{p \in \SS} A^k_{\bm{\ell}}(p) C_p(\spt) h^k, \quad \bm{\ell} \in \LN_{M + 1}.
        \end{equation*}
        We will prove that the statements (a)-(c) in the lemma are further equivalent to
        \begin{equation*}
            \text{\textit{(e)}} \ \ P_{\bm{\ell}} (h) = \bo_{\tilde{a}, u} (h^{M + 2}), \, \forall \bm{\ell} \in \LN_{M + 1}.
        \end{equation*}
        In fact, the equivalence of statements (b) and (e) is treated in \cite[Lemma 3.1]{han2025convergent} with a slight modification of the summation indices. Besides, the implications (e) $\Rightarrow$ (a) and (a) $\Rightarrow$ (c) are trivial.

        (c) $\Rightarrow$ (e): We want to show that if $P_{\bm{\ell}}(h)$ is not simultaneously zero for all $\bm{\ell} \in \LN_{M + 1}$ at some point $\spt \in \Omega_h$, then $\mathcal{L}_h u(\spt)$ depends on some $\partial^{\bm{\ell}} u(\spt)$, $\bm{\ell} \in \LN_{M + 1}$. From the discussion in \Cref{sec:nD_expansion}, for any given values $( u_{\bm{\ell}} )_{ \bm{\ell} \in \LN_{M + 1} }$, there exists a smooth function $v$ such that $\partial^{\bm{\ell}} v(\spt) = u_{\bm{\ell}}$ for $\bm{\ell} \in \LN_{M + 1}$ and $\partial^{\bm{\ell}} \nabla \cdot (a \nabla v)(\spt) = \partial^{\bm{\ell}} f(\spt)$ for $|\bm{\ell}| \le M - 1$. From \eqref{eq:AF_tilde} and \eqref{eq:F(p)}, we can see that $F(p) \in \Span \{ \partial^{\bm{\ell}} f(\spt): |\bm{\ell}| \le M - 1 \}$, which implies that $\sum_{p \in \SS} C_p (\spt) F(p)$ remains constant when we vary $( u_{\bm{\ell}} )_{ \bm{\ell} \in \LN_{M + 1} }$, whereas $\sum_{\bm{\ell} \in \LN_{M + 1}} P_{\bm{\ell}}(h) \partial^{\bm{\ell}} v(\spt)$ can be equal to any number. Thus, neglecting the $\bo_{\tilde{a}, u} (h^{M + 2})$ term, the right-hand side of \eqref{eq:Lh_int_expansion} must be dependent on some $\partial^{\bm{\ell}} u(\spt)$, $\bm{\ell} \in \LN_{M + 1}$. Thus, (c) implies (e).

        (a) $+ \sum_{p \in \SS} p_1^2 c_{p, 0} \neq 0 \Rightarrow$ (d): Suppose that $\sum_{p \in \SS} p_1^2 c_{p, 0} \neq 0$ and set $f_h := \sum_{p \in \SS} C_p F(p)$. Then, \eqref{eq:F(p)_leading} implies that $f_h(\spt) = \alpha a(\spt)^{-1} h^2 (f(\spt) + \so_{\tilde{a}, \tilde{f}} (1))$ for some $\alpha \neq 0$. This confirms the $M$th-order consistency of the scheme according to \Cref{def}, and proves item (d).

        (d) $\Rightarrow$ (e): The same proof of (c) $\Rightarrow$ (e) shows that, if $P_{\bm{\ell}}(h)$ is not simultaneously zero for all $\bm{\ell} \in \LN_{M + 1}$ at some point $\spt \in \Omega_h$, the right-hand side of \eqref{eq:Lh_int_expansion} must be dependent on $u$ in addition to $a$ and $f$. This proves that we cannot find $f_h$ that forms an $M$th-order consistent scheme.
    \end{proof}

    \Cref{eq:C_p_constraint_int} provides an efficient way to find the stencil coefficients. We view \eqref{eq:C_p_constraint_int} as a cascade of linear systems $\mathbb{A}_j \vec{c}_j = \vec{b}_j$ with unknowns $c_{p, j}$ for $j = 0, \ldots, M + 1$, where $\mathbb{A}_j$ is an $\# \LN_{M + 1 - j} \times \# \SS$ matrix, $\vec{b}_j\in \mathbb{R}^{\# \LN_{M + 1 - j}}$, and $\vec{c}_j \in \mathbb{R}^{\# \SS}$. Note that $\mathbb{A}_j$ is a constant matrix according to the end of \Cref{sec:nD_expansion} and it consists of merely the first $\# \LN_{M + 1 - j}$ rows of $\mathbb{A}_0$. Moreover, when $d = 2$, the matrix $\mathbb{A}_j$ in this article and \cite{han2025convergent} only differ by a reordering of the columns.

    Following the result of \Cref{lem:C_p_constraint_int}, the next proposition shows that the set of all possible stencil coefficients has a cascade structure. A more explicit explanation of this proposition is provided after the proof.

    \begin{proposition}
        \label{prop:decomposition}
        For any integer $M \ge -1$ and $\spt \in \Omega_h$, denote $\mathcal{N}_M$ to be the solution space of \eqref{eq:C_p_constraint_int} viewed as a single linear system for $\{ c_{p, j}: 0 \le j \le M + 1, p \in \SS \}$. Moreover, writing $\vec{c}_j = (c_{p, j})_{p \in \SS}$, we define a shifting operator $T: \R^{(M + 1) \# \SS} \to \R^{(M + 2) \# \SS}$, $(\vec{c}_0, \ldots, \vec{c}_M) \mapsto (\vec{0}, \vec{c}_0, \ldots, \vec{c}_M)$, where $\vec{0} \in \R^{\# \SS}$. Then $\mathcal{N}_M$ is a linear space with a decomposition
        \begin{equation}
            \label{eq:decomposition}
            \mathcal{N}_M = T\mathcal{N}_{M - 1} + \mathcal{N}_M^* = \sum_{m = -1}^M T^{M - m} \mathcal{N}_m^*,
        \end{equation}
        where $\mathcal{N}_m^*$ is a linear space and $\dim \mathcal{N}_M^*$ is equal to the dimension of the solution space of \eqref{eq:C_p_constraint_int} with $j = 0$.
    \end{proposition}

    \begin{proof}
        It is clear that $\mathcal{N}_M$ is a linear space from \eqref{eq:C_p_constraint_int}. We can equivalently write \eqref{eq:C_p_constraint_int} as
        \begin{equation}
            \label{eq:proof_decomp_1}
            \sum_{p \in \SS} A^{|\bm{\ell}|}_{\bm{\ell}}(p) c_{p, 0} = 0, \quad \forall \bm{\ell} \in \LN_{M + 1},
        \end{equation}
        when $j=0$ and
        \begin{equation}
            \label{eq:proof_decomp_2}
            \sum_{k = 1}^j \sum_{p \in \SS} A^{|\bm{\ell}| + j - k}_{\bm{\ell}}(p) c_{p, k} = -\sum_{p \in \SS} A^{|\bm{\ell}| + j}_{\bm{\ell}}(p) c_{p, 0}, \quad \forall j =1, \ldots, M + 1, \ \bm{\ell} \in \LN_{M + 1 - j}.
        \end{equation}
        Viewing \eqref{eq:proof_decomp_2} as a single linear system for $(\vec{c}_1, \ldots, \vec{c}_{M+1})$, we know that its solution space $\mathcal{N} (\vec{c}_0)$ can be decomposed into $\mathcal{N}^0 + (\vec{c}_1^{\, *} (\vec{c}_0), \ldots, \vec{c}_{M+1}^{\, *} (\vec{c}_0))$, where $\mathcal{N}^0$ is the solution space of the homogeneous equation and $(\vec{c}_1^{\, *} (\vec{c}_0), \ldots, \vec{c}_{M+1}^{\, *} (\vec{c}_0))$ is a particular solution that is dependent on $\vec{c}_0$. It follows that
        \begin{equation*}
            \mathcal{N}_M = \left\{ (\vec{c}_0, \ldots, \vec{c}_{M+1}): \vec{c}_0 \mbox{ satisfies \eqref{eq:proof_decomp_1}}, \, (\vec{c}_1, \ldots, \vec{c}_{M+1}) \in \mathcal{N} (\vec{c}_0) \right\}.
        \end{equation*}
        Denote
        \begin{equation*}
            \mathcal{N}_M^* := \left\{ \big( \vec{c}_0, \vec{c}_1^{\, *} (\vec{c}_0), \ldots, \vec{c}_{M+1}^{\, *} (\vec{c}_0) \big): \vec{c}_0 \mbox{ satisfies \eqref{eq:proof_decomp_1}} \right\},
        \end{equation*}
        then clearly $\mathcal{N}_M^*$ is a subspace of $\mathcal{N}_M$. Thus, $\mathcal{N}_M = T \mathcal{N}^0 \oplus \mathcal{N}_M^*$. On the other hand, by shifting the index from $j$ to $j + 1$, we see that $(\vec{c}_1, \ldots, \vec{c}_{M+1}) \in \mathcal{N}^0$ if and only if
        \begin{equation*}
            \sum_{k = 0}^j \sum_{p \in \SS} A^{|\bm{\ell}| + j - k}_{\bm{\ell}}(p) c_{p, k + 1} = 0, \quad \forall j=0,\ldots,M, \ \bm{\ell} \in \LN_{M - j}.
        \end{equation*}
        The above equation holds if and only if $(\vec{c}_1, \ldots, \vec{c}_{M+1}) \in \mathcal{N}_{M - 1}$. Therefore, $\mathcal{N}^0 = \mathcal{N}_{M - 1}$ and we get $\mathcal{N}_M = T \mathcal{N}_{M - 1} \oplus \mathcal{N}_M^*$. The second equality in \eqref{eq:decomposition} follows directly from the first equality. Lastly, we can see from the above proof that $\dim \mathcal{N}_M^*$ is equal to the dimension of the solution space of \eqref{eq:proof_decomp_1}, which is the same as \eqref{eq:C_p_constraint_int} with $j = 0$.
    \end{proof}

    Now we correspond the solution $\{ c_{p, j}: 0 \le j \le M + 1, p \in \SS \}$ of \eqref{eq:C_p_constraint_int} to the stencil coefficients $\{ C_p = \sum_{j = 0}^{M + 1} c_{p, j} h^j: p \in \SS \}$, then the shifting operator $T$ sends $C_p$ into $hC_p$. Moreover, let $\{ (C_p^{[M], k})_{p \in \SS}: 1 \le k \le K_M := \dim \mathcal{N}_M^* \}$ be a basis of $\mathcal{N}_M^*$. According to \Cref{prop:decomposition}, any stencil coefficients $C_p$ satisfying \eqref{eq:Lh_int} can be decomposed as
    \begin{equation}
        \label{eq:C_p_structure}
        C_p = \sum_{k = 1}^{K_M} \kappa_{M, k} C_p^{[M], k} + h C_p'
        = \sum_{m = -1}^{M} \sum_{k = 1}^{K_m} \kappa_{m, k} C_p^{[m], k} h^{M - m},
    \end{equation}
    where $\kappa_{m, k} \in \R$ for all $k = 1, \ldots, K_m$, $m= -1, \ldots, M$, and $C_p'$ are some stencil coefficients that satisfy \eqref{eq:Lh_int} with $M$ replaced by $M - 1$.

    Now we consider the compact stencil, that is, $\SS = \{ -1, 0, 1 \}^d$. Since $\mathbb{A}_j$, $0 \le j \le M + 1$ are constant matrices, it is easy to perform symbolic calculation and obtain the dimension $K_M$. The results for 2D and 3D are recorded in \eqref{KM:2D}, \eqref{KM:3D:1}, and \eqref{KM:3D:2}. Using these explicit expressions, it is easy to verify whether the stencil coefficients in \eqref{eq:C_p_structure} satisfy the condition $\sum_{p \in \SS} p_1^2 c_{p, 0} \neq 0$ required by \Cref{lem:C_p_constraint_int}, which leads to a compact, symmetric, and $M$th-order consistent FDM for \eqref{eq:PDE}.

    \subsection{Proofs of maximum consistency order of $d$-dimensional compact FDMs for $d\ge 2$}
    \label{sec:2D_3D_symmetric}

    We start this subsection by proving the maximum consistency order of compact nonsymmetric FDMs for \eqref{eq:PDE}. Recall from the previous subsection that the consistency order of a compact $d$-dimensional finite difference scheme is closely related to the dimension $K_M$. For dimensions $2$ and $3$, we already know that the maximum consistency order is $6$ based on \Cref{lem:C_p_constraint_int}, but it is challenging to directly determine the dimension $K_M$ for $d>3$. Hence, to circumvent this issue, we shall employ a key dimensional reduction technique in the following lemma.

    \begin{lemma}
        \label{lem:max_order_nonsymmetric}
        Let $d \ge 2$. Then any compact finite difference scheme $\mathcal{L}_h u_h = f_h$ for \eqref{eq:PDE} on $\Omega_h$ satisfying the assumptions of \Cref{thm:max_order} is at most 6th-order consistent.
    \end{lemma}
    \begin{proof}
        %
        Let $\mathcal{L}_h u_h = f_h$ be an $M$-th-order consistent FDM with $\mathcal{L}_h u(\spt) := \sum_{p \in \SS} C_p(h) u(\spt + ph)$. Let $\widehat{u}$ be an arbitrary 2-dimensional smooth function and define a restriction operator $R: \R^d \to \R^{2}$, $\mathbf{x} = (x_1, \ldots, x_d) \mapsto (x_1, x_2)$. According to \Cref{lem:C_p_constraint_int}(d), we know that $\mathcal{L}_h u$ is independent of $\partial^{\bm{k}} u$, $\bm{k} \in \LN_{M + 1}^d$. The superscript on $\LN_{M + 1}$ indicates the dimension. In particular, for the constant extension $u$ of $\widehat{u}$, that is, $u(\mathbf{x}) = \widehat{u} (R\mathbf{x})$, the term $\mathcal{L}_h \widehat{u}$ does not depend on $\partial^{\bm{k}} \widehat{u}$, $\bm{k} \in \LN_{M + 1}^2$. Besides, setting $\widehat{\spt} = R \spt$ and $\widehat{\SS} = R \SS$, we get
        \begin{equation*}
            \mathcal{L}_h u (\spt) = \sum_{\widehat{p} \in \widehat{\SS}} \bigg( \sum_{\{p \setsp \spt = R \widehat{p}\}} C_p (\spt) \bigg) \widehat{u} (\widehat{\spt} + \widehat{p}h).
        \end{equation*}
        This naturally induces an operator $\widehat{\mathcal{L}}_h$ on $\widehat{u}$ with stencil coefficients $\widehat{C}_{\widehat{p}} (\widehat{\spt}) := \sum$ $_{\{p \setsp \spt = R \widehat{p}\}} C_p (\spt)$ for some $\spt$ satisfying $\widehat{\spt} = R\spt$. Due to the independence of $\widehat{\mathcal{L}}_h \widehat{u}$ on $\partial^{\bm{k}} \widehat{u}$, $\bm{k} \in \LN_{M + 1}^2$, \Cref{lem:C_p_constraint_int} again shows that the linear system \eqref{eq:C_p_constraint_int} corresponding to $\widehat{C}_{\widehat{p}} (\widehat{\spt})$ has a solution. However, in two dimensions this is true if and only if $M \leq 6$. Therefore, any finite difference scheme in $d$ dimensions can be at most sixth-order consistent.
    \end{proof}

    The following lemma shows that the symmetry of the finite difference stencil implies formal self-adjointness of each differential operator in its Taylor expansion.

    \begin{lemma} \label{propn:Tr:Tr*}
        Let $M \in \mathbb{N}$. Suppose that the stencil coefficients of a finite difference scheme satisfy the following relations
        \begin{align*}
        \wh{C}_p(\spt) & = \wh{C}_{-p}(\spt + ph), \quad \text{for all } p \in \mathcal{S}, \; \spt, \spt+ph \in \Omega_h,\\
        \wh{C}_p(\spt) & = \sum_{k=0}^{M+1} h^k \wh{c}_{p,k}(\spt) + \bo_a(h^{M+2}), \quad
        \text{for all } p \in \mathcal{S},
        \end{align*}
        where $\wh{c}_{p,k}$ are some smooth functions. Define the differential operators $T_r$ and $T_r^*$ by
        \be \label{Tr:Trs}
        \begin{aligned}
         (T_r v)(\bm x)
         & := \sum_{p\in\SS}\sum_{|\bm\ell|\le r}
           \frac{p^{\bm\ell}}{\bm\ell!}
           \widehat c_{p,r-|\bm\ell|}(\bm x)\partial^{\bm\ell}v(\bm x),\\
         (T_r^*v)(\bm x)
         & := \sum_{p\in\SS}\sum_{|\bm\ell|\le r}
           \frac{(-1)^{|\bm\ell|}p^{\bm\ell}}{\bm\ell!}
           \partial^{\bm\ell}\!\left(\widehat c_{p,r-|\bm\ell|}v\right)(\bm x),
       \end{aligned}
        \ee
        for $v \in C^{\infty}(\overline{\Omega})$, $\bm x \in \Omega$, and $r=0,\ldots,M+1$. Then, $(T_r v)(\spt) = (T_r^* v)(\spt)$ for $v \in C^{\infty}(\overline{\Omega})$, $\spt \in \Omega_h$, and $r=0,\ldots,M+1$.
    \end{lemma}
    \begin{proof}
        By the stencil's symmetry and the fact that $\mathcal{S}= -\mathcal{S}$, we have
        \be \label{adj:stencil}
            \wh{\mathcal{L}}_h v(\spt) = \sum_{p\in\SS} \wh{C}_p(\spt) v(\spt+hp)
            = \sum_{p\in\SS} \wh{C}_{-p}(\spt) v(\spt-hp)
            = \sum_{p\in\SS} \wh{C}_{p}(\spt - hp) v(\spt-hp).
        \ee
        Taking the Taylor expansion about any fixed point $\spt \in \Omega_h$, we have
        \begin{align*}
             \sum_{p\in\SS}\widehat C_p(\spt) v(\spt+hp)
             & =\sum_{r=0}^{M+1}h^r
             (T_r v) (\spt) +\bo_a(h^{M+2}), \\
             \sum_{p\in\SS}\widehat C_p(\spt-hp)v(\spt-hp)
             & =\sum_{r=0}^{M+1}h^r
             (T_r^* v)(\spt) + \bo_a(h^{M+2}).
        \end{align*}
        By \eqref{adj:stencil} and the fact that $(T_r v)(\spt)$, $(T_r^* v)(\spt)$ are independent of $h$, it follows that $(T_r v)(\spt) = (T_r^* v)(\spt)$ for $v \in C^{\infty}(\overline{\Omega})$, $\spt \in \Omega_h$, and $r=0,\ldots,M+1$.
    \end{proof}

    We now state and prove two lemmas concerning the highest attainable consistency order for compact symmetric FDMs. The next one describes the nonexistence of a compact, symmetric, $M$-th order consistent FDM for \eqref{eq:PDE} with $M \ge 5$, when $|\nabla a|^2/a^2 - 2\Delta a/a$ is nonconstant on $\Omega$.

    \begin{lemma} \label{Qa:noncnst}
        Let $\Omega = (0,1)^d$ with $d \ge 2$, and $a$ be a sufficiently smooth variable function, where $a(\mathbf{x}) \ge a_0 >0$. If there exists a compact symmetric finite difference scheme for \eqref{eq:PDE} satisfying the assumptions of \Cref{thm:max_order} with consistency order $M \ge 5$, then $|\nabla a|^2/a^2 - 2\Delta a/a$ is constant on $\Omega$.
    \end{lemma}

    \begin{proof}
        Let $M \ge 5$. By assumption, we have the following scheme for \eqref{eq:PDE}
        \begin{align*}
        \mathcal{L}_h u(\spt) & := \sum_{p \in \SS} C_p (\spt) u(\spt + ph) = f_h(\spt) + \bo_a(h^{M+2}), \quad \text{for all } \spt \in \Omega_h,
        \quad \text{and}\\
        C_p (\spt) & = C_{-p}(\spt + ph), \quad \text{for all } p \in \mathcal{S}, \; \spt, \spt+ph \in \Omega_h,
        \end{align*}

        We introduce the change of variables $\rho := \sqrt{a}$, $v = \rho u$, and $F=f/\rho$.
        Define  $C(\mathbf{x}) := \frac{|\nabla a(\mathbf x)|^2}{a(\mathbf x)^2} -\frac{2\Delta a(\mathbf x)}{a(\mathbf x)}$ for all $\mathbf{x} \in \Omega$. Then, the transformed problem becomes
        \be \label{v:PDE:trans}
        \wh{\mathcal{L}}v:=- \Delta v - \tfrac{1}{4}C(\mathbf{x}) v = F \text{ in } \Omega, \quad v= \rho g \text{ on } \partial\Omega.
        \ee
        We also observe that
        \be \label{arhov}
             -\nabla\cdot\left(a\nabla\frac v\rho\right)
             =-\rho\Delta v+v\Delta\rho
             =\rho\left(-\Delta v-\tfrac14 Cv\right)
             = \rho F.
        \ee
        Now, suppose that there is a compact, symmetric, and $M$th-order finite difference scheme for \eqref{arhov} with $M \ge 5$, which we denote by $\mathcal{L}_h u_h(\spt) = f_h (\spt)$, $\spt \in \Omega_h$. Define $\wh{\mathcal{L}}_h v(\spt)
        :=-\frac{1}{\rho(\spt)} \mathcal{L}_h(v/\rho)(\spt)$. It follows that
        \[
        \wh{\mathcal{L}}_h v(\spt) = \sum_{p \in \mathcal{S}} \wh{C}_p(\spt) v(\spt + hp), \quad
        \wh{C}_p(\spt) := -\frac{C_{p}(\spt)}{\rho(\spt) \rho(\spt +ph)},
        \]
        where we can further write
        \[
        \wh{C}_p(\mathbf{x}) = \sum_{r=0}^{M+1} h^r \wh{c}_{p,r}(\mathbf{x}) + \bo_a(h^{M+2}), \quad \text{for all } \mathbf{x} \in \Omega, \; p \in \mathcal{S},
        \]
        and for some smooth functions $\wh{c}_{p,r}$. Additionally, we can also write
        \[
        \frac{f_h(\mathbf{x})}{\sqrt{a(\mathbf{x})}} = h^2 F(\mathbf{x}) + \sum_{j=1}^{M-1} h^{j+2} (\mathcal{P}_j F)(\mathbf{x}) + \bo_a(h^{M+2}),
        \quad \text{for all } \mathbf{x} \in \Omega,
        \]
        where $F:=f/\sqrt{a}$, and each $\mathcal{P}_j$ is a finite-order linear differential operator independent of $h$ and $F$.
        Hence, we have the following compact, symmetric, and $M$-th order consistent, with $M \ge 5$, finite difference scheme for \eqref{v:PDE:trans}
        \begin{align*}
            \wh{\mathcal{L}}_h v(\spt) & = - \frac{f_h(\spt)}{\sqrt{a(\spt)}} + \bo_a(h^{M+2}) = - h^2 F(\spt) - \sum_{j=1}^{M-1} h^{j+2} (\mathcal{P}_j F)(\spt) + \bo_a(h^{M+2})\\
            & = h^2 (\Delta + \tfrac{1}{4}C I) v(\spt) + \sum_{j=1}^{M-1} h^{j+2} \mathcal{P}_j (\Delta + \tfrac{1}{4}C I) v(\spt) + \bo_a(h^{M+2}),
        \end{align*}
        for all $\spt \in \Omega_h$.

        Using the following Taylor expansion
        \[
             v(\spt+ph)
             =\sum_{|\bm\ell|\le M+1}
               \frac{p^{\bm\ell}h^{|\bm\ell|}}{\bm\ell!}
               \partial^{\bm\ell}v(\spt)+\bo_a(h^{M+2}),
        \]
        and recalling the expression for $\wh{C}_p(\spt)$, we have
        \[
             \wh{\mathcal{L}}_h v(\spt)
             =\sum_{r=0}^{M+1}h^r
               (T_r v)(\spt)+\bo_a(h^{M+2})
             = \sum_{r=0}^{M+1}h^r
               \sum_{|\bm\ell|\le r}
               t_{r,\bm\ell} (\spt)
               \partial^{\bm\ell}v(\spt)+\bo_a(h^{M+2}),
        \]
        where $T_r v$ is defined in \eqref{Tr:Trs} and
        \[
            t_{r,\bm\ell}(\spt):=\frac{1}{\bm\ell!}\sum_{p\in\SS}
                      p^{\bm\ell}\widehat c_{p,r-|\bm\ell|}(\spt).
        \]
        Let $I$ be the identity operator. Our calculation, so far, implies that
        \be \label{terms:trl:Pj}
        \begin{aligned}
        \sum_{r=0}^{M+1}h^r
               \sum_{|\bm\ell|\le r}
               t_{r,\bm\ell} (\spt)
               \partial^{\bm\ell}v(\spt)
       & =
       h^2 (\Delta + \tfrac{1}{4}C I) v(\spt)  \\
       & \quad + \sum_{j=1}^{M-1} h^{j+2} \mathcal{P}_j (\Delta + \tfrac{1}{4}C I) v(\spt) + \bo_a(h^{M+2}).
       \end{aligned}
        \ee

        Our next major task is to do term-by-term comparisons, list some equations that must hold given that \eqref{terms:trl:Pj} holds, and show how these equations yield the claim of this lemma. The selected necessary equations discussed below play a direct role in proving that $C(\mathbf{x})$ is constant on $\Omega$. We heavily rely on symbolic calculation to carry out this task.

        We start by performing a term-by-term comparison on the following equation
        \[
            \sum_{|\bm\ell|\le 2}
               t_{2,\bm\ell} (\spt)
               \partial^{\bm\ell}v(\spt)
            =
            (\Delta + \tfrac{1}{4}C I) v(\spt).
        \]
        Then, we must have
        \[
        t_{2,2\vec e_i}=1; \quad
        t_{2,\vec e_i+\vec e_j} =0, \; \text{for all } i<j; \quad
        t_{2,\vec e_i}=0; \quad
        t_{2,\vec{0}}=C/4.
        \]

        Next, we perform a term-by-term comparison on the following equation
        \be \label{P2:eqn}
        \sum_{|\bm\ell|\le4} t_{4,\bm\ell}(\spt)
           \partial^{\bm\ell}v(\spt)
        = \mathcal{P}_2 (\Delta +\tfrac{1}{4} C I) v (\spt).
        \ee
        Since the highest partial derivatives on the left-hand side take the form of $\partial^{\bm\ell}v$ with $|\bm \ell | =4$, the linear differential operator $\mathcal{P}_2$ must have the following form
        \[
        \mathcal P_2=\sum_{1 \le i\le j \le d}b_{i,j}^{(2)}\partial^{\vec e_i+\vec e_j}
        +\sum_{i=1}^d \gamma_i^{(2)}\partial^{\vec e_i}+w^{(2)}I
        \]
        for some smooth real-valued functions $b_{i,j}^{(2)}, \gamma_i^{(2)}, w^{(2)}$ on $\Omega$ such that \eqref{P2:eqn} holds and $T_4 v= T_4^* v$ in \cref{propn:Tr:Tr*}. This forces
        \be \label{t4:soln}
        \begin{aligned}
            & b_{i,i}^{(2)} = \tfrac{1}{12}; \quad
            b_{i,j}^{(2)}=0 \; \text{for all } i<j; \quad \gamma_i^{(2)} = 0; \quad w^{(2)} = \tfrac{1}{48}C + \kappa, \\
            & t_{4,2\vec{e}_i+2\vec e_j} = \tfrac{1}{6} \; \text{for all } i \ne j; \quad
            t_{4,2\vec{e}_i}= \tfrac{1}{24} C +\kappa; \quad
            t_{4,\vec{e}_i}=\tfrac{1}{24} \partial^{\vec{e}_i}C
        \end{aligned}
        \ee
        for some constant $\kappa \in \R$.

        Next, we perform a term-by-term comparison on the following equation
        \be \label{P4:eqn}
         \sum_{|\bm\ell|\le6}
         t_{6,\bm\ell}(\spt)
               \partial^{\bm\ell}v(\spt)
         = \mathcal{P}_4 (\Delta +\tfrac{1}{4} C I) v (\spt).
        \ee
        Since the highest partial derivatives on the left-hand side take the form of $\partial^{\bm\ell}v$ with $|\bm \ell | = 6$, the linear differential operator $\mathcal{P}_4$ must have the following form
        \[
        \mathcal{P}_4 =
         \sum_{|\bm\alpha|=4}\eta_{\bm\alpha}\partial^{\bm\alpha}
         +\sum_{|\bm\alpha|=3}\zeta_{\bm\alpha}\partial^{\bm\alpha}
         +\sum_{1 \le i\le j \le d}b_{i,j}^{(4)}\partial^{\vec{e}_i+\vec{e}_j}
         +\sum_{i=1}^{d} \gamma_i^{(4)}\partial^{\vec e_i}+w^{(4)}I.
        \]
        for some smooth real-valued functions $\eta_{\bm\alpha}, \zeta_{\bm\alpha}, b_{i,j}^{(4)}, \gamma_i^{(4)}, w^{(4)}$ on $\Omega$ such that \eqref{P4:eqn} holds and $T_6 v= T_6^* v$ in \cref{propn:Tr:Tr*}. This forces
        \be \label{t6:soln}
        \begin{aligned}
            & \eta_{4\vec e_i} = \tfrac{1}{360}; \quad
            \eta_{2\vec e_i+2\vec e_j} = \tfrac{1}{90} \; \text{for all } i < j; \quad
            \eta_{\bm\alpha} =0 \; \text{for all other } |\bm\alpha|=4; \\
            & \zeta_{\bm\alpha}=0, \;  |\bm\alpha|=3; \quad b_{i,i}^{(4)} =\tfrac{1}{360} C + \tfrac{1}{12} \kappa; \quad
            b_{i,j}^{(4)} = 0 \; \text{for all } i<j.
        \end{aligned}
        \ee

        Our final task is to analyze equations involving $\gamma_i^{(4)}$ for all $1\le i \le d$. To this end, we first compare both sides of \eqref{P4:eqn} involving $\partial^{3\vec{e}_i} v$. Observe that the corresponding coefficient on the left-hand sides of \eqref{P4:eqn} is
        \[
        t_{6,3\vec{e}_i}
        = \tfrac{1}{6} \sum_{p\in\SS} p_i^{3} \widehat c_{p,3}(\spt)
        = \tfrac{1}{6} \sum_{p\in\SS} p_i \widehat c_{p,3}(\spt)
        = \tfrac{1}{6} t_{4, \vec{e}_i}
        = \tfrac{1}{144} \partial^{\vec{e}_i}C,
        \]
        where we have used \eqref{t4:soln} to arrive at the last equality. Meanwhile, the coefficient of $\partial^{3\vec{e}_i} v$ on the right-hand side is $\gamma_i^{(4)} + \tfrac{1}{360} \partial^{\vec{e}_i}C$. The following is our first equation involving $\gamma_i^{(4)}$
        \be \label{gam:1}
        \tfrac{1}{144} \partial^{\vec{e}_i}C = \gamma_i^{(4)} + \tfrac{1}{360} \partial^{\vec{e}_i}C \quad \text{or} \quad
        \gamma^{(4)}_i = \tfrac{1}{240} \partial^{\vec{e}_i}C.
        \ee
        Let us find the second equation involving $\gamma_i^{(4)}$. Comparing the terms involving $\partial^{\vec{e}_i+2\vec{e}_j} v$, where $i \neq j$, on both sides of \eqref{P4:eqn}, we have
        \[
        t_{6,\vec{e}_i + 2\vec{e}_j} = \gamma^{(4)}_i + \tfrac{1}{180} \partial^{\vec{e}_i}C.
        \]
        On the other hand, comparing the terms involving $\partial^{2\vec{e}_i+2\vec{e}_j} v$ on both sides of \eqref{P4:eqn}, we have
        \[
        t_{6,2\vec{e}_i + 2\vec{e}_j} = b^{(4)}_{i,i} + b^{(4)}_{j,j} + \tfrac{1}{4}C \eta_{2\vec{e}_i + 2\vec{e}_j} = 2 \left(\tfrac{1}{360} C + \tfrac{1}{12} \kappa\right) + \tfrac{1}{360} C = \tfrac{1}{120} C + \tfrac{1}{6} \kappa
        \]
        by \eqref{t6:soln}. Now, comparing the terms involving $\partial^{\vec{e}_i+2\vec{e}_j} v$ in the relation $T_6 v = T_6^* v$ in \cref{propn:Tr:Tr*}, we have
        \[
        t_{6,\vec{e}_i + 2\vec{e}_j} = - t_{6,\vec{e}_i + 2\vec{e}_j} + 2 \partial^{\vec{e}_i} (t_{6,2\vec{e}_i + 2\vec{e}_j})
        \quad
        \text{or}
        \quad
        t_{6,\vec{e}_i + 2\vec{e}_j} = \tfrac{1}{120} \partial^{\vec{e}_i}C,
        \]
        The second equation involving $\gamma_i^{(4)}$ is as follows
        \be \label{gam:2}
        \gamma^{(4)}_i + \tfrac{1}{180} \partial^{\vec{e}_i}C = \tfrac{1}{120} \partial^{\vec{e}_i}C
        \quad \text{or} \quad
        \gamma^{(4)}_i = \tfrac{1}{360} \partial^{\vec{e}_i}C.
        \ee
        Solving \eqref{gam:1} and \eqref{gam:2} gives $\partial^{\vec{e}_i}C = 0$ for all $1\le i \le d$ and $\nabla C = \vec{0}$ for any $\spt \in \Omega_h$. This implies that $C = \frac{|\nabla a |^2}{a^2} -\frac{2\Delta a }{a }$ is constant on $\Omega$. The proof is completed.
    \end{proof}

     The following lemma complements the previous one and shows that if $\frac{|\nabla a|^2}{a^2} - \frac{2\Delta a}{a}$ is constant on $\Omega = (0,1)^d$ with $d\ge 2$ and symmetry is required, 6th-order consistency can be achieved. This order matches the highest possible consistency order attainable by a nonsymmetric compact FDM for \eqref{eq:PDE}.

    \begin{lemma}
         \label{lem:sym:6th}
        Let $\frac{|\nabla a|^2}{a^2} - \frac{2\Delta a}{a}$ be constant on $\Omega = (0,1)^d$, $d \ge 2$. Then there exists a compact, symmetric, and 6th-order consistent finite difference scheme $\mathcal{L}_h u_h = f_h$ for \eqref{eq:PDE} on $\Omega_h$ satisfying the assumptions of \Cref{thm:max_order}.
    \end{lemma}

    \begin{proof}
        Consider the model problem \eqref{eq:PDE} with $d\ge2$. Since $a(\mathbf{x}) \ge a_0 >0$ for all $\mathbf{x} \in \Omega$, we can take a positive function $\rho(\mathbf{x})$ such that $a(\mathbf{x})=\rho^2(\mathbf{x})$. By assumption, $\frac{|\nabla a|^2}{a^2} - \frac{2\Delta a}{a} = C$ for some constant $C$. If $v:=\rho u$ and $F:=f/\rho$, applying change-of-variable yields
        \begin{equation} \label{model:helm}
        \widehat{\mathcal{L}}v:= - \Delta v - \tfrac{1}{4} C v = F \quad \text{in} \; \Omega; \quad v = \rho g \quad \text{on} \; \partial \Omega.
        \end{equation}
        Next, we recall the definition of the centered second difference. For a smooth function $v$, define its centered second difference along the $i$-th axis
        \[
        \Delta_{h,i}v := v(\mathbf{x} + h\vec{e}_i) - 2v(\mathbf{x}) + v(\mathbf{x} - h\vec{e}_i).
        \]
        The following relations can be verified by elementary calculations using Taylor expansions:
        \begin{align*}
            \sum_{i=1}^{d} \Delta_{h,i} v & = h^2 \Delta v + \tfrac{h^4}{12} \sum_{i=1}^{d} \partial^{4 \vec{e}_i} v + \tfrac{h^6}{360} \sum_{i=1}^{d} \partial^{6 \vec{e}_i} v + \bo_a(h^8), \\
            \sum_{1\le i < j \le d} \Delta_{h,i} \Delta_{h,j} v  & = h^4 \sum_{1\le i < j \le d} \partial^{2(\vec{e}_i+\vec{e}_j)} v + \tfrac{h^6}{12} \sum_{1\le i < j \le d} (\partial^{4\vec{e}_i+2\vec{e}_j} v + \partial^{2\vec{e}_i+4\vec{e}_j} v)\\
            & \quad + \bo_a(h^8), \\
            \sum_{1\le i < j < \ell \le d} \Delta_{h,i} \Delta_{h,j} \Delta_{h,\ell} v  & = h^6 \sum_{1\le i < j < \ell \le d} \partial^{2(\vec{e}_i + \vec{e}_j + \vec{e}_\ell)} v + \bo_a(h^8).
        \end{align*}
        Define
        \begin{align*}
        \widehat{\mathcal{L}}_h & := \sum_{i=1}^{d} \Delta_{h,i} + \left(\tfrac{1}{6} + \tfrac{C h^2}{720}\right) \sum_{1\le i < j \le d} \Delta_{h,i} \Delta_{h,j} \\
        & \quad + \tfrac{1}{30}  \sum_{1\le i < j < \ell \le d} \Delta_{h,i} \Delta_{h,j} \Delta_{h,\ell} + \tfrac{1}{4} C h^2 \left(1-\tfrac{C h^2}{48} + \tfrac{C^2 h^4}{5760}\right)I,\\
        \mathcal{P}_h & := h^2 I + \tfrac{h^4}{12}(\Delta - \tfrac{1}{4}C I) + \tfrac{h^6}{360} (\Delta^2 - \tfrac{1}{4}C \Delta + \tfrac{1}{16}C^2 I) + \tfrac{h^6}{180} \sum_{1\le i < j \le d} \partial^{2(\vec{e}_i + \vec{e}_j)},
        \end{align*}
        where $I$ is the identity operator. It follows that for a sufficiently smooth function $v$ and any $\spt \in \Omega_h$, we have
        \begin{align*}
        \widehat{\mathcal{L}}_h v(\spt) & = (\Delta + \tfrac{1}{4}C I) \mathcal{P}_h v(\spt) + \bo_a(h^8) = \mathcal{P}_h (\Delta + \tfrac{1}{4}C I)  v(\spt) + \bo_a(h^8) \\
        & = - \mathcal{P}_hF(\spt) + \bo_a(h^8),
        \end{align*}
        where we used \eqref{model:helm} to obtain the final line. Since $\widehat{\mathcal{L}}_h$ consists of sums and products of centered second differences which are symmetric, we have a compact, symmetric, and 6th-order FDM for \eqref{model:helm}. That is,
        \[
        \widehat{\mathcal{L}}_h v(\spt) = \sum_{p \in \mathcal{S}} \widehat{C}_p(\spt) v(\spt + ph) = - \mathcal{P}_hF(\spt) + \bo_a(h^8), \quad \spt \in\Omega_h,
        \]
        where $\widehat{C}_p(\spt) = \widehat{C}_{-p}(\spt+ph)$ for all $p \in \mathcal{S}$.
        Recall that $v = \rho u$. Define
        \[
        \mathcal{L}_h u (\spt) := -\rho(\spt)\widehat{\mathcal{L}}_h(\rho u) (\spt) \quad \text{and} \quad f_h(\spt) := \rho(\spt) \mathcal{P}_h(\tfrac{f}{\rho})(\spt).
        \]
        Then, for any $\spt \in \Omega_h$, we have
        \[
        \mathcal{L}_h u (\spt) = -\rho(\spt)\widehat{\mathcal{L}}_h(\rho u) (\spt) = \rho(\spt)\mathcal{P}_h(\tfrac{f}{\rho})(\spt) + \bo_a(h^8) = f_h (\spt) + \bo_a(h^8),
        \]
        which implies
        \[
        \mathcal{L}_h u (\spt) = \sum_{p \in \mathcal{S}} C_p (\spt) u(\spt + ph)  = f_h (\spt) + \bo_a(h^8),
        \]
        where $C_p (\spt) := -\sqrt{a(\spt)} \widehat{C}_p(\spt) \sqrt{a(\spt + ph)}$ and $C_p (\spt) = C_{-p} (\spt+ph)$ for all $p \in \mathcal{S}$.
        This gives us a compact, symmetric, and 6th-order FDM for \eqref{eq:PDE} under the assumption that $\frac{|\nabla a|^2}{a^2} - \frac{2\Delta a}{a}$ is constant on $\Omega = (0,1)^d$, $d \ge 2$.
    \end{proof}

    \subsection{Proofs of SPD and convergence rates for $1\le d \le 4$}
    \label{subsec:SPD:conv}
    Next, we prove \cref{SPD}, which says that the discretization matrix is SPD, and the convergence rates coincide with the consistency orders.

    \begin{proof}[Proof of \cref{SPD}]
         Our goal is to show that for all nonzero vector $\vec{v}$, we have $\vec{v}^\tp A_h \vec{v} > 0$. Recall that $\mathcal{S} = \{-1,0,1\}^d$, where $d=1,2,3$, and $\dot{\iota}(\spt)$, $\spt \in \Omega_h$, is a map from the point $\spt$ to its row index in the matrix $A_h$. The $\dot{\iota}(\spt)$-th entry of the vector $A_h \vec{v}$ is given by
        \[
        (A_h \vec{v})_{\dot{\iota}(\spt)} = C_{\vec{0}} (\spt) \vec{v}_{\dot{\iota}(\spt)} + \sum_{\substack{p \in \mathcal{S} \backslash \{\vec{0}\},\\ \spt + ph \in \Omega_h}} C_p(\spt) \vec{v}_{\dot{\iota}(\spt+ph)}, \quad \spt \in \Omega_h.
        \]
        As a result, we have
        \begin{align*}
        \vec{v}^\tp A_h \vec{v} & = \sum_{\spt \in \Omega_h} C_{\vec{0}}(\spt) \vec{v}^2_{\dot{\iota}(\spt)} + \sum_{\spt \in \Omega_h} \sum_{\substack{p \in \mathcal{S} \backslash \{\vec{0}\},\\ \spt + ph \in \Omega_h}} C_{p}(\spt) \vec{v}_{\dot{\iota}(\spt)} \vec{v}_{\dot{\iota}(\spt + ph)}\\
        & = -\sum_{\spt \in \Omega_h} \sum_{\substack{p \in \mathcal{S} \backslash \{\vec{0}\},\\ \spt + ph \in \Omega_h}} C_p(\spt) \vec{v}^2_{\dot{\iota}(\spt)} + \sum_{\spt \in \Omega_h} \sum_{\substack{p \in \mathcal{S} \backslash \{\vec{0}\},\\ \spt + ph \in \Omega_h}} C_{p}(\spt) \vec{v}_{\dot{\iota}(\spt)} \vec{v}_{\dot{\iota}(\spt + ph)} \\
        & \quad -\sum_{\spt \in \Omega_h} \sum_{\substack{p \in \mathcal{S} \backslash \{\vec{0}\},\\ \spt + ph \in \partial \Omega_h}} C_p(\spt) \vec{v}^2_{\dot{\iota}(\spt)},
        \end{align*}
        where $\partial \Omega_h := \partial \Omega \cap (h \mathbb{Z}^{d})$, $d=1,2,3$. By \eqref{sym:cond}, we have $C_p(\spt) = C_{-p}(\spt + ph)$ for $\spt,\spt + ph \in \Omega_h$, and $p \in \mathcal{S}\backslash \{\vec{0}\}$. Moreover, if $p \in \mathcal{S}$, then $-p \in \mathcal{S}$ as well. This gives
        \begin{align}
        \nonumber
        \vec{v}^\tp A_h \vec{v}
        & = - \frac{1}{2}\sum_{\spt \in \Omega_h} \sum_{\substack{p \in \mathcal{S} \backslash \{\vec{0}\},\\ \spt + ph \in \Omega_h}} C_p(\spt) (\vec{v}_{\dot{\iota}(\spt)} - \vec{v}_{\dot{\iota}(\spt + ph)})^2 \\
        \label{vTAv:start}
        & \quad -\sum_{\spt \in \Omega_h} \sum_{\substack{p \in \mathcal{S} \backslash \{\vec{0}\},\\ \spt + ph \in \partial \Omega_h}} C_p(\spt) \vec{v}^2_{\dot{\iota}(\spt)}.
        \end{align}
        Notice from \cref{thm:1D} for $d=1$, \eqref{2d:explicit} for $d=2$, and \eqref{3d:explicit} for $d=3$ that
        \[
            C_p(\spt) = \theta_p a(\spt + ph/2) + \bo_a (h^2), \quad \theta_p < 0, \quad \text{for all } \spt \in \Omega_h.
        \]
        Hence, for any sufficiently small $h$, we have
        \begin{align}
        \nonumber
        \vec{v}^\tp A_h \vec{v}
        & \ge \sum_{\spt \in \Omega_h} \sum_{\substack{p \in \mathcal{S} \backslash \{\vec{0}\},\\ \spt + ph \in \Omega_h}} \left(- \tfrac{1}{4} \theta_p a_0\right) (\vec{v}_{\dot{\iota}(\spt)} - \vec{v}_{\dot{\iota}(\spt + ph)})^2 \\
        \label{vTAv}
        & \quad + \sum_{\spt \in \Omega_h} \sum_{\substack{p \in \mathcal{S} \backslash \{\vec{0}\},\\ \spt + ph \in \partial \Omega_h}} \left(- \tfrac{1}{2} \theta_p a_0\right) \vec{v}^2_{\dot{\iota}(\spt)},
        \end{align}
        where we use the fact that $a(\mathbf{x}) \ge a_0 >0$ for all $\mathbf{x} \in \Omega$. I.e., $ \vec{v}^\tp A_h \vec{v}  \ge 0$ for any $\vec{v}$. If $ \vec{v}^\tp A_h \vec{v} =0$, then
        \begin{align*}
        & \vec{v}_{\dot{\iota}(\spt)} = \vec{v}_{\dot{\iota}(\spt + ph)}, \quad \text{for all } \spt \in \Omega_h, \; p \in \mathcal{S} \backslash \{\vec{0}\}, \; \spt + ph \in \Omega_h, \quad \text{and}\\
        & \vec{v}_{\dot{\iota}(\spt)} = 0, \quad \text{for all } \spt \in \Omega_h, \; p \in \mathcal{S} \backslash \{\vec{0}\}, \; \spt + ph \in \partial \Omega_h,
        \end{align*}
        because $- \tfrac{1}{4} \theta_p a_0 >0$ and $- \tfrac{1}{2} \theta_p a_0 >0$. This implies that $\vec{v} = \vec{0}$. Hence, for any nonzero $\vec{v}$, we have $\vec{v}^\tp A_h \vec{v} >0$, which proves the matrix $A_h$ obtained from discretization is SPD for all sufficiently small $h$. This proves item (a).

        Since $A_h$ is SPD and its off-diagonal entries are nonpositive, it follows that $A_h$ is an M-matrix \cite[Definition 10.3.3 and Theorem 10.3.3]{G97book}.

        Next, we prove the convergence in the 2-norm. Define $\overline{\Omega}_h := \Omega_h \cup \partial \Omega_h$. For notational convenience, we extend $\dot{\iota}$ from $\Omega_h$ to $\overline{\Omega}_h$ by assigning an auxiliary index to every point in $\partial\Omega_h$, while keeping the original indices for points in $\Omega_h$. Define
        \[
        \vec{\tilde{v}}_{\dot{\iota}(\spt)} :=
        \begin{cases}
            \vec{v}_{\dot{\iota}(\spt)}, & \spt \in \Omega_h,\\
            0, & \spt \in \partial \Omega_h,
        \end{cases}
        \quad
        \text{and}
        \quad
        \theta_0 :=\min_{p \in \mathcal{S}\backslash\{\vec{0}\}} |\theta_p|.
        \]
        Recall that we use a uniform grid with mesh size $h = N^{-1}$. From \eqref{vTAv}, we have
        \begin{align}
        \nonumber
        \vec{v}^\tp A_h \vec{v} & \ge \frac{1}{4} \theta_0 a_0 \sum_{\spt \in \Omega_h} \sum_{\substack{p \in \mathcal{S} \backslash \{\vec{0}\}, \\ \spt + ph \in \overline{\Omega}_h}} (\vec{\tilde{v}}_{\dot{\iota}(\spt)} - \vec{\tilde{v}}_{\dot{\iota}(\spt+ph)})^2\\
        \label{coercive:1}
         & \ge \frac{1}{4} \theta_0 a_0 \sum_{\spt \in \Gamma_j}\sum_{k=1}^{N} (\vec{\tilde{v}}_{\dot{\iota}(\spt+k\vec{e}_j h)} - \vec{\tilde{v}}_{\dot{\iota}(\spt+(k-1)\vec{e}_j h)})^2,
        \end{align}
        for each $j=1,\ldots,d$, where $\Gamma_j := \{h(i_1,\ldots,i_{j-1},0,i_{j+1},\ldots,i_d):1\le i_m \le N-1, m \neq j\}$. To further find a lower bound, we observe that
        \[
        \vec{\tilde{v}}_{\dot{\iota}(\spt+\ell \vec{e}_j h)} = \sum_{k=1}^{\ell} (\vec{\tilde{v}}_{\dot{\iota}(\spt+k\vec{e}_j h)} - \vec{\tilde{v}}_{\dot{\iota}(\spt+(k-1)\vec{e}_j h)}),
        \]
        which by the Cauchy-Schwarz inequality implies that
        \[
            |\vec{\tilde{v}}_{\dot{\iota}(\spt+\ell \vec{e}_j h)}|^2 \le N  \sum_{k=1}^{N} (\vec{\tilde{v}}_{\dot{\iota}(\spt+k\vec{e}_j h)} - \vec{\tilde{v}}_{\dot{\iota}(\spt+(k-1)\vec{e}_j h)})^2.
        \]
        It follows that
        \begin{equation} \label{coercive:2}
             N (N-1) \sum_{\spt \in \Gamma_j} \sum_{k=1}^{N} (\vec{\tilde{v}}_{\dot{\iota}(\spt+k\vec{e}_j h)} - \vec{\tilde{v}}_{\dot{\iota}(\spt+(k-1)\vec{e}_j h)})^2 \ge \sum_{\spt \in \Gamma_j} \sum_{\ell=1}^{N-1} |\vec{\tilde{v}}_{\dot{\iota}(\spt+\ell \vec{e}_j h)}|^2 = \|\vec{v}\|^2_2.
        \end{equation}
        Combining \eqref{coercive:1} and \eqref{coercive:2}, and using the fact that $d=1,2,3$, we have
        \[
        \vec{v}^\tp A_h \vec{v} \ge \frac{1}{4} \theta_0 a_0  h^2 \|\vec{v}\|^2_2
        \quad
        \text{and}
        \quad
        \|A_h^{-1}\|_2 \le 4 \theta_0^{-1} a_0^{-1} h^{-2}.
        \]
        Our finite difference scheme gives $A_h \overrightarrow{u|_{\Omega_h}} = \vec{f}_h + \vec{\tau}_h$, where the vector $\vec{\tau}_h$ stores truncation errors with
        \[
        \|\vec{\tau}_h\|_{\infty} \le
        \begin{cases}
             C_{\tau,\infty} h^{M+2}, & \text{if } d=1,\\
             C_{\tau,\infty} h^{6}, & \text{if } d=2,3,
        \end{cases}
        \quad
        \text{and}
        \quad
        \|\vec{\tau}_h\|_{2,h} \le
        \begin{cases}
             C_{\tau,2} h^{M+2}, & \text{if } d=1,\\
             C_{\tau,2} h^{6}, & \text{if } d=2,3,
        \end{cases}
        \]
        for some positive constants $C_{\tau,\infty}$ and $C_{\tau,2}$.
        At the same time, $A_h \vec{u}_h = \vec{f}_h$. This gives $A_h (\overrightarrow{u|_{\Omega_h}}-\vec{u}_h) = \vec{\tau}_h$. By calculation, we have
        \[
        \|\overrightarrow{u|_{\Omega_h}}-\vec{u}_h\|_{2,h} = h^{d/2}\|A_h^{-1}  \vec{\tau}_h\|_{2} \le \|A_h^{-1}\|_2 \|\tau_h\|_{2,h} \le
        \begin{cases}
             \frac{4C_{\tau,2}}{\theta_0 a_0} h^{M}, & \text{if } d=1,\\
             \frac{4C_{\tau,2}}{\theta_0 a_0} h^{4}, & \text{if } d=2,3,
        \end{cases}
        \]
        for all sufficiently small $h$.

        We now prove the convergence in the $\infty$-norm. Consider the comparison function
        \[
        \phi(x_1,\ldots,x_d):= e^{\mu} - e^{\mu x_1}, \quad \mu := 1 + \left\| \partial_{x_1}a/a \right\|_{\infty}.
        \]
        It is clear that $\phi \ge 0$ on $\overline{\Omega}$ and
        \[
        \mathcal{L}\phi = -\nabla \cdot (a \nabla \phi) = \mu (\partial_{x_1}a + a \mu) e^{\mu x_1} \ge \mu a_0 =:c_\phi.
        \]
        For all sufficiently small $h$, we have
        \[
        \mathcal{L}_h \phi(\spt) = h^2 \mathcal{L}\phi(\spt) + \bo_{a}(h^3) \ge \tfrac{1}{2} c_{\phi} h^2, \quad \text{for all } \spt \in \Omega_h.
        \]
        Since $\phi \ge 0$ on $\overline{\Omega}$ and $C_p(\spt) \le 0$ for all $\spt \in \Omega_h$ and $p \in \mathcal{S}\backslash \{\vec{0}\}$, we have $A_h \overrightarrow{\phi|_{\Omega_h}} \ge \tfrac{1}{2}c_{\phi}h^2 \vec{1}$,
        where $\vec{1}$ is a vector of ones. Since $A_h$ is a nonsingular M-matrix, we have $A_h^{-1} \ge 0$ and $A_h^{-1} \vec{1} \le 2c^{-1}_{\phi} h^{-2} \overrightarrow{\phi|_{\Omega_h}}$. Since all elements are positive, we have $\|A_h^{-1}\|_{\infty} \le 2c^{-1}_{\phi} e^{\mu} h^{-2}$. Because $A_h (\overrightarrow{u|_{\Omega_h}}-\vec{u}_h) = \vec{\tau}_h$, it follows that
        \[
        \|\overrightarrow{u|_{\Omega_h}}-\vec{u}_h\|_{\infty} \le \|A_h^{-1}\|_\infty \|\vec{\tau}_h\|_{\infty} \le \begin{cases}
             C_{\tau,\infty} (2c^{-1}_{\phi} e^{\mu}) h^{M}, & \text{if } d=1,\\
             C_{\tau,\infty} (2c^{-1}_{\phi} e^{\mu}) h^{4}, & \text{if } d=2,3,
        \end{cases}
        \]
        for all sufficiently small $h$. This proves item (b).
    \end{proof}

    In what follows, we provide the analogous SPD and convergence proofs for $d=4$.

    \begin{proof}[Proof of \cref{SPD:4}]
        For $d=4$, recall that the stencil coefficients $C_{p}(\spt)$ and the right-hand term $f_h(\spt)$ in \eqref{stencil:explicit:uf} take the forms in \eqref{Cp:fh:d4}.
        Under the assumption that $|\nabla \tilde{a}|^2 + 2\partial^{2\vec{e}_i}\tilde{a} - \Delta \tilde{a} \le 0$ on $\Omega$, for $i=1,\ldots,4$, we have $C_{p}(\spt) \le 0$ for $p = \pm \vec{e}_i$ and $1\le i \le 4$. Meanwhile, for all sufficiently small $h$, we have
        \[
            C_{p}(\spt) = -\frac{1}{6} a(\spt + ph/2) + \bo_{a}(h^{2}), \quad  \pm (\vec{e}_i\pm\vec{e}_j), \; 1\le i < j \le 4.
        \]
        The inequality in \eqref{vTAv:start} yields
        \begin{align*}
        \vec{v}^\tp A_h \vec{v}
        & \ge \sum_{\spt \in \Omega_h} \sum_{\substack{p = \pm(\vec{e}_1 + \vec{e}_2),\\ \spt + ph \in \Omega_h}} \left(\tfrac{a_0}{24}\right) (\vec{v}_{\dot{\iota}(\spt)} - \vec{v}_{\dot{\iota}(\spt + ph)})^2 \\
        & \quad + \sum_{\spt \in \Omega_h} \sum_{\substack{p = \pm(\vec{e}_1 + \vec{e}_2),\\ \spt + ph \in \partial \Omega_h}} \left(\tfrac{a_0}{24}\right) \vec{v}^2_{\dot{\iota}(\spt)}.
        \end{align*}
        I.e., $\vec{v}^\tp A_h \vec{v} \ge 0$ for any $\vec{v}$. If $\vec{v}^\tp A_h\vec{v} = 0$, then
        \begin{align*}
        & \vec{v}_{\dot{\iota}(\spt)} = \vec{v}_{\dot{\iota}(\spt \pm(\vec{e}_1 + \vec{e}_2)h)}, \quad \text{for all } \spt \in \Omega_h, \; \spt \pm(\vec{e}_1 + \vec{e}_2) h \in \Omega_h, \quad \text{and}\\
        & \vec{v}_{\dot{\iota}(\spt)} = 0, \quad \text{for all } \spt \in \Omega_h, \; \spt \pm(\vec{e}_1 + \vec{e}_2) h \in \partial \Omega_h,
        \end{align*}
        because $a_0/24 > 0$. This implies that $\vec{v} = \vec{0}$. This proves that the discretization matrix $A_h$ is SPD for all sufficiently small $h$. To prove the coercivity, we replace $\vec{e}_j$ in \eqref{coercive:1} with $\vec{e}_1 + \vec{e}_2$. The convergence rates for $d=4$ can now be proved by using the same ideas in the proof of \cref{SPD}.
    \end{proof}

    \section{Conclusion}
    \label{sec:conclusion}
    We characterize the highest consistency orders attainable by compact symmetric FDMs for the variable-coefficient Poisson equation on uniform grids in a $d$-dimensional hypercube. This study is motivated by the limited existing theory for high-order compact symmetric FDMs with general variable coefficients, despite the well-known 2nd-order construction. Constructing higher-order schemes is challenging because compactness limits the available degrees of freedom, high-order consistency imposes numerous local Taylor conditions, and matrix symmetry requires compatibility between neighboring stencils.

    For $d=1$, compact symmetric FDMs can attain arbitrarily high consistency orders. For $d\ge2$, the highest attainable consistency order is 4 when $\frac{|\nabla a|^2}{a^2}-\frac{2\Delta a}{a}$ is nonconstant on $\Omega$, and 6 when it is constant. The latter matches the maximum order attainable without requiring matrix symmetry. Furthermore, under the stated assumptions and for sufficiently small grid sizes, the proposed schemes in dimensions $1\le d\le4$ produce SPD matrices and converge at their respective consistency orders in both the weighted $2$-norm and the $\infty$-norm. We provide explicit stencil coefficients, a general construction procedure, and numerical experiments demonstrating the convergence and performance of the schemes.

    As future work, we plan to extend our approach to other PDEs and more general domains. Although some techniques developed here may carry over, new ideas will be needed to construct interior and boundary stencils that jointly satisfy high-order consistency and matrix symmetry, particularly when using larger stencils or treating irregular domains.

    \section{Statement on the use of AI}
    OpenAI GPT-5.6 Sol was used to assist with the literature review, explore the possibility of obtaining the SPD property, review existing links between the SPD property of FDM matrices and scheme convergence, and simplify algebraic expressions in the proof of Lemma 7.5. The authors manually carried out the entire rigorous proofs of the SPD property, the scheme’s convergence, and Lemma 7.5. They independently reviewed and verified all AI-assisted outputs and take full responsibility for the entire content of the manuscript.
\bibliographystyle{siamplain}
\bibliography{ref}
\end{document}